\documentclass{amsart}

\usepackage[colorlinks=true, urlcolor=black, citecolor=black, linkcolor=black, hyperfootnotes=true]{hyperref}
\usepackage{amssymb}
\usepackage{aliascnt}%for theorem counters
\usepackage{accents}

\numberwithin{equation}{section}

\newtheorem{thm}{Theorem}[section]

\newaliascnt{prp}{thm}
\newtheorem{prp}[prp]{Proposition}
\aliascntresetthe{prp}

\newaliascnt{cor}{thm}
\newtheorem{cor}[cor]{Corollary}
\aliascntresetthe{cor}

\newaliascnt{lem}{thm}
\newtheorem{lem}[lem]{Lemma}
\aliascntresetthe{lem}

\theoremstyle{definition}

\newaliascnt{dfn}{thm}
\newtheorem{dfn}[dfn]{Definition}
\aliascntresetthe{dfn}

\newaliascnt{xpl}{thm}
\newtheorem{xpl}[xpl]{Example}
\aliascntresetthe{xpl}

\newaliascnt{rmk}{thm}
\newtheorem{rmk}[rmk]{Remark}
\aliascntresetthe{rmk}

\author{Tristan Bice}
\thanks{Tristan Bice is supported by the GA\v{C}R project 22-07833K and RVO: 67985840.}
\address{Institute of Mathematics of the Czech Academy of Sciences, \v{Z}itn\'a 25, Prague}
\email{bice@math.cas.cz}

\makeatletter
\@namedef{subjclassname@2020}{\textup{2020} Mathematics Subject Classification}
\makeatother

\makeatletter
\protected\def\xvcenter{%
  \hbox\bgroup$\everyvbox{\everyvbox{}\aftergroup\m@th\aftergroup$\aftergroup\egroup}%
  \vcenter
}
\DeclareRobustCommand{\midscript}[1]{
  \mathchoice{\mid@script\scriptstyle{#1}}
    {\mid@script\scriptstyle{#1}}
    {\mid@script\scriptscriptstyle{#1}}
    {\mid@script\scriptscriptstyle{#1}}
}
\newcommand{\mid@script}[2]{
  \vcenter{\hbox{$\m@th#1#2$}}
}

\makeatletter

\title[Noncommutative Pierce Duality]{Noncommutative Pierce Duality between Steinberg Rings and Ample Ringoid Bundles}

\begin{document}

\begin{abstract}
Classic work of Pierce and Dauns-Hofmann shows that biregular rings are dual to simple ring bundles over Stone spaces.  We extend this duality to Steinberg rings, a purely algebraic generalisation of Steinberg algebras, and ringoid bundles over ample groupoids.  We base this largely on an even more general extension of Lawson's noncommutative Stone duality, specifically between Steinberg semigroups, a generalisation of Boolean inverse semigroups, and category bundles over ample groupoids.
\end{abstract}

\maketitle

\section*{Introduction}

\subsection*{Background}

Rings are commonly analysed by representing them as continuous sections of certain kinds of ring bundles.  For example, classic results in algebraic geometry say that a commutative ring can always be represented on a bundle of local rings over a spectral space, namely its Zariski spectrum of prime ideals.  Other classic results due to Pierce (see \cite{Pierce1967}) and Dauns-Hofmann (\cite{DaunsHofmann1966}) say that a biregular ring can always be represented on a bundle of simple rings over a Stone space, namely its spectrum of maximal ideals.  However, in the noncommutative case, biregular rings may have few maximal ideals, which makes the resulting Pierce representation somewhat trivial.  To gain any insight into such noncommutative rings, it is natural to consider more general ringoid bundles over \'etale groupoids.

This idea first appeared in work of Feldman and Moore \cite{FeldmanMoore1977} on Cartan subalgebras of von Neumann algebras, which was later extended to C*-algebras by Kumjian \cite{Kumjian1986} and Renault \cite{Renault2008}.  In the past decade, attention has been directed to more algebraic analogs of groupoid C*-algebras known as Steinberg algebras, which generalise the Leavitt path algebras associated to directed graphs.  These were introduced in \cite{Steinberg2010} and \cite{ClarkFarthingSimsTomforde2014} and consist of locally constant functions on ample groupoids taking values/coefficients in $\mathbb{C}$ or some other ring.

Recently, a number of people have shown how to reconstruct the groupoid from the Steinberg algebra, together with some additional data coming from a Cartan-like subalgebra or its normaliser semigroup \textendash\, see \cite{BrownClarkanHuef2017}, \cite{CarlsenRout2018}, \cite{Steinberg2019} and \cite{BiceClark2021}.  What has been lacking so far is an abstract algebraic characterisation of Steinberg algebras, together with a generalised Pierce representation analogous to the Kumjian-Renault Weyl groupoid representation for C*-algebras.

The goal of the present paper is to develop an algebraic version of the Kumjian-Renault Weyl groupoid representation which does indeed apply to Steinberg algebras and even more general rings consisting of sections of ringoid bundles over ample groupoids.  These are characterised as abstract `Steinberg rings', namely rings together with some extra data coming from appropriate subsemigroups and an expectation.  We further build this up to a duality between appropriate categories which extends the duality above due to Pierce and Dauns-Hofmann.

In fact, we are able to do most of the heavy lifting already at the semigroup level.  Indeed, there has been a parallel series of papers in the world of semigroups/lattices starting with the classic Stone duality \cite{Stone1936} between Boolean algebras and Stone spaces, leading to Lawson's noncommutative Stone duality in \cite{Lawson2010} and further extensions with Kudryavtseva in \cite{KudryavtsevaLawson2016} and \cite{KudryavtsevaLawson2017}, as well as another extension by the author and Starling in \cite{BiceStarling2018}.  Here we extend Lawson's work in a different direction, where we consider semigroups arising from sections of category bundles over ample groupoids.  These are characterised as abstract `Steinberg semigroups' likewise consisting of semigroups with a distinguished subsemigroup and an expectation.  This leads to a generalisation of Lawson's duality which is then easily extended to the desired generalisation of Pierce and Dauns-Hofmann's duality.

\subsection*{Related Work}
In the later stages of preparing the first draft of this paper, we learnt about some similar work being done independently at the same time in \cite{ArmstrongCastroClarkCourtneyLinMcCormickRamaggeSimsSteinberg2021}.  We would like to thank the authors for sharing their preprint and discussing how it relates to the present paper, as we outline here.

On the topological side, \cite{ArmstrongCastroClarkCourtneyLinMcCormickRamaggeSimsSteinberg2021} deals with groupoid twists, which can be viewed as groupoid bundles with an additional group action -- see \autoref{TwistedRemark} below.  These can be built up to `ample line bundles', i.e. ample ringoid bundles $\pi:C\rightarrow G$ where $\pi^{-1}[G^0]=R\times G^0$, for some fixed unital ring $R$, which in \cite{ArmstrongCastroClarkCourtneyLinMcCormickRamaggeSimsSteinberg2021} is also commutative.  On the other hand, here we consider more general category bundles and ringoid bundles, where the unit fibres can even be noncommutative.  Steinberg's `local bisection hypothesis' is also not required in our work.

On the algebraic side, we consider Steinberg rings, namely quadruples $(A,S,Z,\Phi)$ where $A$ is a ring with subsemigroups $S$ and $Z$ and $\Phi:A\rightarrow S$ is an expectation, satisfying certain algebraic conditions.  On the other hand, the structures considered in \cite{ArmstrongCastroClarkCourtneyLinMcCormickRamaggeSimsSteinberg2021} are algebraic quasi-Cartan pairs $(A,B)$, where $A$ is an $R$-algebra with subalgebra $B$, again satisfying certain somewhat different algebraic conditions.  However, quasi-Cartan pairs can again be viewed as special kinds of Steinberg rings when we take $\Phi$ to be the unique expectation onto $B$, $Z$ to be the idempotents in $B$ and $S$ to be the orthospan of the $Z$-invertible normalisers of $B$ -- see \autoref{QuasiCartanPairs} below.

Consequently, our more general framework could be viewed as subsuming the duality in \cite{ArmstrongCastroClarkCourtneyLinMcCormickRamaggeSimsSteinberg2021}.  However, the approach of Armstrong et al. provides a different perspective that will no doubt also be useful in future work.

The final thing to point out is that, at the categorical level, the work of Armstrong et al. in \cite{ArmstrongCastroClarkCourtneyLinMcCormickRamaggeSimsSteinberg2021}, just like the original work of Dauns-Hofmann in \cite{DaunsHofmann1966}, only yields a duality with respect to isomorphisms.  In contrast, we obtain a fully-fledged extension of Pierce's duality by considering more general morphisms like those in \cite{Pierce1967}.  In a recent sequel to \cite{ArmstrongCastroClarkCourtneyLinMcCormickRamaggeSimsSteinberg2021}, Steinberg considers a somewhat different class of morphisms in \cite{Steinberg2021} -- it would be an interesting topic of future work to see if this approach could be unified with ours in some way.

\subsection*{Overview}

Our focus in the present paper is on category bundles, which underlie all the structures mentioned in the literature above.  The key point here is that convolution can already be defined on suitable sections of (zero) category bundles.  Naturally one can then ask if the category bundle can be recovered from the multiplicative structure of its sections.  Much of our paper is devoted to showing this can be done, as long as we restrict our attention to ample category bundles.

The basic idea is to first recover the base groupoid as ultrafilters with respect to a certain `domination relation' on sections.  The total category can then be recovered from certain equivalence classes defined by these ultrafilters.  As this process relies only on algebraic data, it even allows us to represent certain abstract algebraic structures as sections of ample category bundles.  This leads to the notion of a Steinberg semigroup, whose algebraic properties are designed to characterise the sections of ample category bundles under convolution, as shown in \autoref{SteinbergBundleCharacterisation}.

As with all algebraic structures, Steinberg semigroups form a category where the morphisms are structure preserving maps.  The next natural question to ask is what these morphisms correspond to in terms of the associated bundles.  Here it turns out we need to consider morphisms that combine maps of the base groupoid with maps of the total category in a compatible way.  The resulting Pierce morphisms allow us to turn the correspondence between Steinberg semigroups and ample category bundles into an equivalence of categories.  The final section is devoted to showing how all this can still be done in the presence of some extra additive structure.

\section{Semimodules}\label{Semimodules}

The semimodules we will focus on here are closely related to the semigroups we considered in \cite{Bice2022}.  These `structured semigroups' were born out of a desire to extend the relevant theory of inverse semigroups so it can be applied to semigroups arising elsewhere, in particular in C*-algebra theory.  Indeed, what we showed in \cite{Bice2022} is that structured semigroups can be represented on \'etale groupoid bundles, much as Kumjian and Renault showed how to represent C*-algebas with Cartan subalgebras on saturated Fell line bundles over their Weyl groupoid.

Here our goal is much the same, except that we want to move even closer to the Kumjian-Renault theory.  For example, we want the representation to apply not just to semigroups but also to larger rings and C*-algebras in which they may sit.  We also want the elements of the semigroup/ring to be represented as total sections rather than the partial slice-sections considered in \cite{Bice2022}.  Last but not least, we want to upgrade the representation to a duality between appropriate categories, something which had not even been achieved in the Kumjian-Renault theory.

Consequently, we will need to work with slightly more restrictive `well-structured semigroups'.  Nevertheless, it is instructive to first review the definition of structured semigroups from \cite{Bice2022}.  For this, we first need to consider the following normality/diagonality conditions from \cite[\S1.3 and Definition 5.5]{Bice2022}.

\begin{dfn}
Let $S$ be a semigroup.  We call a subset $N\subseteq S$ \emph{normal}, \emph{binormal}, \emph{trinormal} or \emph{diagonal} respectively if
\begin{align}
\tag{Normal}s\in S\qquad&\Rightarrow\qquad sN=Ns.\\
\tag{Binormal}s,t\in S\ \text{ and }\ st,ts\in N\qquad&\Rightarrow\qquad sNt\subseteq N.\\
\tag{Trinormal}s,t\in S,\ \ tsn=nts=n\in N\ \text{ and }\ st,ts\in N\qquad&\Rightarrow\qquad snt\in N.\\
\tag{Diagonal}s,t\in S\ \text{ and }\ sn,nt,n\in N\qquad&\Rightarrow\qquad snt\in N.
\end{align}
\end{dfn}

We also denote the \emph{centre} of any $N\subseteq S$ by
\[\mathsf{Z}(N)=\{z\in N:\forall n\in N\ (nz=zn)\}.\]
Now we can define structured semigroups like in \cite[Definition 1.6]{Bice2022} (actually we used a slight weakening of trinormality there but that will not be relevant here).  Let us also adopt the convention that
\begin{center}
\textbf{(sub)semigroups are never allowed to be empty}.
\end{center}

\begin{dfn}
A \emph{structured semigroup} is a triple $(S,Z,N)$ where $S$ is a semigroup with trinormal subsemigroup $N$ and binormal subsemigroup $Z$ contained in $\mathsf{Z}(N)$.
\end{dfn}

For example, any semigroup $S$ with non-empty centre becomes a structured semigroup by taking $N=S$ and $Z=\mathsf{Z}(S)$.  If $S$ is an inverse semigroup, we can get a structured semigroup by instead taking $N=Z=\mathsf{E}(S)$ where $\mathsf{E}(S)$ denotes the idempotents of $S$, i.e.
\[\mathsf{E}(S)=\{e\in S:e=ee\}.\]
If $A$ is C*-algebra with Cartan subalgebra $C$ then the normaliser semigroup
\[S=\{a\in A:aCa^*+a^*Ca\subseteq C\}\]
also forms a structured semigroup when we take $N=Z=C$.

As mentioned above, in this last example we would like to represent not just the semigroup $S$ but also the larger C*-algebra $A$ as sections of some bundle.  To do this, we will need to consider not just the product within $S$, but also the product between elements of $S$ and $A$.  The resulting structure $(A,S)$ is something like a module, but for semigroups rather than rings, and hence we name it as follows.

\begin{dfn}
A \emph{semimodule} is a pair $(A,S)$ with $S\subseteq A$ and an $A$-valued associative product on $S\times A\cup A\times S$, i.e. such that, for all $s,t\in S$ and $a\in A$,
\[(st)a=s(ta),\quad(sa)t=s(at)\quad\text{and}\quad(as)t=a(st).\]
\end{dfn}

Note that, for $(st)a$ and $a(st)$ above to be defined, we must have $st\in S$.  In other words, when $(A,S)$ is semimodule, it is implicit that $SS\subseteq S$ so $S$ is itself a semigroup with respect to the product restricted to $S$.  A semimodule could thus be viewed as a semigroup $S$ acting (from the left and right) on some larger $A\supseteq S$, with the additional requirement that the action restricted to $S$ coincides with the semigroup product within $S$.

To represent a semimodule $(A,S)$ as sections of a bundle, we will also need an appropriate function which brings elements of $A$ back down to $S$, namely a `shiftable expectation' $\Phi:A\rightarrow S$.  Again, the inspiration here comes from operator algebras, where a Cartan subalgebra $C$ of $A$ is, by definition, the range of an expectation, which is shiftable w.r.t. its normaliser semigroup $S$ and also has `bistable' range.

\begin{dfn}
An \emph{expectation} is a map $\Phi:A\rightarrow S$ on a semimodule $(A,S)$ with
\begin{align}
\label{Idempotent}\tag{Idempotent}\Phi(\Phi(a))&=\Phi(a)\quad\text{and}\\
\label{Homogeneous}\tag{Homogeneous}\Phi(a)\Phi(b)&=\Phi(\Phi(a)b)=\Phi(a\Phi(b)),\\
\intertext{for all $a,b\in A$.  We call $\Phi$ \emph{shiftable} if, for all $a\in A$ and $s\in S$,}
\label{Shiftable}\tag{Shiftable}\Phi(sa)s&=s\Phi(as).\\
\intertext{We call $Z\subseteq S$ \emph{bistable} if, for all $s,t\in S$,}
\label{Bistable}\tag{Bistable}st\in Z\qquad&\Rightarrow\qquad \Phi(s)t,s\Phi(t)\in Z.
\end{align}
\end{dfn}

Note that \eqref{Homogeneous} implies $\mathrm{ran}(\Phi)$ is a subsemigroup of $S$.  Also note that the first two defining properties can be rephrased in terms of the range of $\Phi$ as
\begin{align}
\tag{Idempotent$'$}&\mathrm{ran}(\Phi)=\{s\in S:\Phi(s)=s\}.\\
\tag{Homogeneous$'$}\label{Homogeneous'}r\in\ &\mathrm{ran}(\Phi)\quad\Rightarrow\quad\Phi(ar)=\Phi(a)r\text{ and }\Phi(ra)=r\Phi(a).
\end{align}
Further note that $\mathrm{ran}(\Phi)$ is bistable precisely when, for all $a,b\in S$,
\[\tag{$\Phi$-Bistable}ab\in\mathrm{ran}(\Phi)\qquad\Rightarrow\qquad\Phi(a)b=\Phi(a)\Phi(b)=a\Phi(b).\]
Indeed, $\Phi(a)b\in\mathrm{ran}(\Phi)$ iff $\Phi(a)b=\Phi(\Phi(a)b)=\Phi(a)\Phi(b)$ and likewise for $a\Phi(b)$.

Now we are in a position to define the well-structured semigroups and semimodules that will be our primary focus.

\begin{dfn}\label{WSS}
A \emph{well-structured semimodule} is a quadruple $(A,S,Z,\Phi)$ where
\begin{enumerate}
\item $(A,S)$ is a semimodule with shiftable expectation $\Phi:A\rightarrow S$.
\item $Z$ is a binormal bistable subsemigroup of $S$ contained in $\mathsf{Z}(\mathrm{ran}(\Phi))$.
\end{enumerate}
When $A=S$ we call $(S,Z,\Phi)$ a \emph{well-structured semigroup}.
\end{dfn}

First we note that these are, in particular, structured semigroups.

\begin{prp}\label{WellStructuredSemigroups}
If $(S,Z,\Phi)$ is a well-structured semigroup then $\mathrm{ran}(\Phi)$ is both trinormal and diagonal in $S$.  In particular, $(S,Z,\mathrm{ran}(\Phi))$ is a structured semigroup.
\end{prp}

\begin{proof}
For trinormality, take $s,t\in S$ with $tsr=r\in\mathrm{ran}(\Phi)$ and $st\in\mathrm{ran}(\Phi)$.  Then
\[srt=s\Phi(r)t=s\Phi(tsr)t=st\Phi(srt)=\Phi(stsrt)=\Phi(srt)\in\mathrm{ran}(\Phi).\]
For diagonality just note that, as in \cite[Example 4.8]{Bice2022}, $sr,r,rt\in\mathrm{ran}(\Phi)$ implies
\[srt=s\Phi(rt)=sr\Phi(t)=\Phi(srt)\in\mathrm{ran}(\Phi).\qedhere\]
\end{proof}

\subsection{Examples}

An arbitrary semigroup $S$ becomes well-structured by taking $\Phi$ to be the identity $\mathrm{id}_S$ and $Z=\mathsf{Z}(S)$ or even $Z=\mathsf{E}(\mathsf{Z}(S))$.  Indeed, applying the theory we will develop to the case $\Phi=\mathrm{id}_A$ and $Z=\mathsf{E}(\mathsf{Z}(S))$ is all we would need to recover the original duality of Pierce and Dauns-Hofmann.

\begin{prp}
$(S,\mathsf{E}(\mathsf{Z}(S)),\mathrm{id}_S)$ is well-structured, for any semigroup $S$.
\end{prp}

\begin{proof}
First note that $\mathsf{E}(\mathsf{Z}(S))$ is indeed a subsemigroup of $\mathsf{Z}(\mathrm{ran}(\Phi))=\mathsf{Z}(S)$, as $d,e\in\mathsf{E}(\mathsf{Z}(S))$ implies that $dede=ddee=de$.  To see that $\mathsf{E}(\mathsf{Z}(S))$ is binormal in $S$, take $a,b\in S$ and $e\in\mathsf{E}(\mathsf{Z}(S))$.  Then certainly $aeb=abe\in\mathsf{Z}(S)\mathsf{Z}(S)\subseteq\mathsf{Z}(S)$ and, if $ab\in\mathsf{E}(\mathsf{Z}(S))$ too, then $aebaeb=ababee=abe=aeb$, showing that $aeb\in\mathsf{E}(\mathsf{Z}(S))$, as required.  Lastly, just note that the identity on $S$ is certainly a shiftable expectation making any subset bistable.
\end{proof}

Examples with non-identity expectations again come from inverse semigroups and operator algebras.  Indeed, if $C$ is a Cartan subalgebra of $A$ then $(A,S,C,\Phi)$ is a well-structured semimodule, where $S$ is the normaliser semigroup and $\Phi$ is the expectation from $A$ onto $C$.  On the other hand, any inverse $\wedge$-semigroup (i.e. an inverse semigroup which is also a $\wedge$-semilattice w.r.t. the canonical ordering) gives rise to a well-structured semigroup.  Indeed, by \cite[Theorem 1.9]{Leech1995}, an inverse semigroup $S$ is a $\wedge$-semilattice precisely when $a^\geq\cap\mathsf{E}(S)=\{e\in\mathsf{E}(S):e\leq a\}$ has a maximum $\Phi_\mathsf{E}(a)$, for all $a\in S$.

\begin{prp}\label{WellStructuredInverseSemigroup}
If $S$ is an inverse $\wedge$-semigroup then $(S,\mathsf{E}(S),\Phi_\mathsf{E})$ is a well-structured semigroup, where $\Phi_\mathsf{E}$ is the shiftable expectation on $S$ defined by
\[\Phi_\mathsf{E}(a)=\max(a^\geq\cap\mathsf{E}(S))=a\wedge aa^{-1}=a\wedge a^{-1}a.\]
\end{prp}

\begin{proof}
Certainly $\mathsf{E}(S)=\mathrm{ran}(\Phi_\mathsf{E})$ is a commutative subsemigroup of $S$.  Moreover, for any $a\in S$ and $e\in\mathsf{E}(S)$, note $ses^{-1}\in\mathsf{E}(S)$ and $se=ss^{-1}se=ses^{-1}s$, showing that $s\mathsf{E}(S)\subseteq\mathsf{E}(S)s$.  The reverse inclusion follows by a dual argument, showing that $\mathsf{E}(S)$ is normal and hence binormal in $S$.

Next we need to show that $\Phi_\mathsf{E}$ is a shiftable expectation.  As $e\leq e$, for all $e\in\mathsf{E}(S)$, \eqref{Idempotent} is immediate.  By \cite[Theorem 1.11(c)]{Leech1995}, \eqref{Homogeneous} also holds so $\Phi_\mathsf{E}$ is an expectation.  For all $a,s\in S$, \cite[Theorem 1.11(e)]{Leech1995} says that $s^{-1}\Phi_\mathsf{E}(a)s=\Phi_\mathsf{E}(s^{-1}as)$ and hence
\[\Phi_\mathsf{E}(sa)s=\Phi_\mathsf{E}(sa)ss^{-1}s=ss^{-1}\Phi_\mathsf{E}(sa)s=s\Phi_\mathsf{E}(s^{-1}sas)=ss^{-1}s\Phi_\mathsf{E}(as)=s\Phi_\mathsf{E}(as),\]
showing that $\Phi_\mathsf{E}$ is also shiftable.

For all $a,b\in S$, note that $\Phi_\mathsf{E}(a)b,a\Phi_\mathsf{E}(b)\leq ab$.  If $ab\in\mathsf{E}(S)$ then it follows that $\Phi_\mathsf{E}(a)b,a\Phi_\mathsf{E}(b)\in\mathsf{E}(S)$ too, showing that $\mathsf{E}(S)$ is also bistable and hence $(S,\mathsf{E}(S),\Phi_\mathsf{E})$ is a well-structured semigroup.
\end{proof}

To obtain more general examples of well-structured semigroups and semimodules, we now move on to consider category bundles and their continuous sections.

\section{Category Bundles}\label{CategoryBundles}

First we define bundles and make some elementary observations.

As usual, we call a map $\rho$ on a topological space \emph{locally injective} if every point in its domain $\mathrm{dom}(\rho)$ has a neighbourhood on which $\rho$ is an injective map.

\begin{dfn}
A \emph{bundle} is a continuous open surjection $\rho:C\twoheadrightarrow G$ between topological spaces $C$ and $G$.  If $\rho$ is also locally injective then $\rho$ is an \emph{\'etale bundle}.
\end{dfn}

We view any $\rho:C\rightarrow G$ as a subset of $G\times C$ (so functions \emph{are} their graphs).  We call $a:G\rightarrow C$ a \emph{section} of $\rho$ if $a\subseteq\rho^{-1}$, i.e. $\rho(a(g))=g$, for all $g\in G$.

\begin{prp}\label{EtaleBundleOpenRange}
Any continuous section of an \'etale bundle has open range.
\end{prp}

\begin{proof}
If $\rho:C\twoheadrightarrow G$ is an \'etale bundle then, for any $g\in G$, we have an open neighbourhood $O$ of $a(g)$ on which $\rho$ is injective.  As $a$ is continuous, $a^{-1}[O]$ is open and hence so is $O\cap\mathrm{ran}(a)=O\cap\rho^{-1}[a^{-1}[O]]$, showing that $\mathrm{ran}(a)$ is open.
\end{proof}

If $\rho:C\rightarrow G$ is continuous then it is homeomorphic to its domain $C$ and, moreover, any continuous section $a$ of $\rho$ is a homeomorphism onto its range $\mathrm{ran}(a)$ with inverse map $\rho|_{\mathrm{ran}(a)}$.  In particular, $\rho$ is Hausdorff precisely when $C$ is Hausdorff.  As long as $\rho$ has at least one continuous section then $G$ is also Hausdorff.

\begin{prp}\label{ClosedRange}
Any continuous section of a Hausdorff bundle has closed range.
\end{prp}

\begin{proof}
Take a continuous section $a$ of a Hausdorff bundle $\rho:C\rightarrow G$.  Further take a net $(g_\lambda)\subseteq G$ with $a(g_\lambda)\rightarrow c\in C$.  Then $g_\lambda=\rho(a(g_\lambda))\rightarrow\rho(c)$, as $\rho$ is continuous, and hence $a(g_\lambda)\rightarrow a(\rho(c))$, as $a$ is continuous.  If $C$ is Hausdorff then limits are unique and hence $c=a(\rho(c))\in\mathrm{ran}(a)$, showing that $\mathrm{ran}(a)$ is closed.
\end{proof}

Next we recall some basic category theory.  We consider categories as sets of arrows/morphisms with a partial associative product.  The source and range units of any element $a$ of a category $C$ are denoted by $\mathsf{s}(a)$ and $\mathsf{r}(a)$ respectively, so $ab$ is defined iff $\mathsf{s}(a)=\mathsf{r}(b)$.  We denote the units, invertibles and composable pairs by
\begin{align*}
C^0&=\{\mathsf{s}(a):a\in C\}=\{\mathsf{r}(a):a\in C\}.\\
C^\times&=\{a\in C:\exists a^{-1}\in C\ (a^{-1}a=\mathsf{s}(a)\text{ and }aa^{-1}=\mathsf{r}(a))\}.\\
C^2&=\{(a,b):a,b\in C\text{ and }\mathsf{s}(a)=\mathsf{r}(b)\}.
\end{align*}
As usual, we call $C^\times$ the \emph{core} of $C$.  We call $B\subseteq C$ a \emph{bisection} or \emph{slice} if the source and range maps are injective on $B$.  A category is \emph{topological} if it carries a topology making the product, source and range maps continuous.  We call a topological category \emph{\'etale} if both the source and range are also locally injective open maps (and hence local homeomorphisms).  The following result extends \cite[2.$\Rightarrow$3.]{Resende2007} and shows that our \'etale categories are the same as those in \cite{KudryavtsevaLawson2017}.

\begin{prp}
If $C$ is an \'etale category then $C^0$ is an open subset of $C$ and the product on $C^2$ is a locally injective open map and hence a local homeomorphism.
\end{prp}

\begin{proof}
If $C$ is an \'etale category then, in particular, $\mathsf{s}$ and $\mathsf{r}$ are open maps and hence $C^0=\mathsf{s}[C]=\mathsf{r}[C]$ is an open subset.  Now take $(a,b)\in C^2$ and open neighbourhoods $O\ni a$ and $N\ni b$.  As $\mathsf{s}$ and $\mathsf{r}$ are locally injective, we may shrink $O$ and $N$ if necessary to make them slices.  Shrinking them further to $O\cap\mathsf{s}^{-1}[\mathsf{r}[N]]$ and $N\cap\mathsf{r}^{-1}[\mathsf{s}[O]]$ if necessary, we may further assume that $\mathsf{s}[O]=\mathsf{r}[N]$.  Now, for any other $c,c'\in O$ and $d,d'\in N$ with $cd=c'd'$, note that $\mathsf{r}(c)=\mathsf{r}(cd)=\mathsf{r}(c'd')=\mathsf{r}(c')$ and hence $c=c'$.  Likewise $d=d'$, showing that the product is injective on $C^2\cap O\times N$.  Next note that we have a continuous section $\phi$ of $\mathsf{r}|_{\mathsf{r}^{-1}[\mathsf{r}[O]]}$ with range $ON$ defined on $\mathsf{r}[O]$ by $\phi(e)=\mathsf{r}|_{O}^{-1}(e)(\mathsf{r}|_{N}^{-1}\circ\mathsf{s}\circ\mathsf{r}|_{O}^{-1})(e)$.
Thus $ON$ is open, by \autoref{EtaleBundleOpenRange}, showing that the product is a locally injective open map.
\end{proof}

A category $G$ is a \emph{groupoid} if $G=G^\times$ and similar standard definitions apply.  Firstly, a groupoid $G$ is \emph{topological} if it carries a topology making the inverse, product, source and range maps continuous.  Again, a topological groupoid is \emph{\'etale} if both the source and range are are also locally injective open maps.  Equivalently, $G$ is \'etale if it has a basis of open slices forming an inverse semigroup under pointwise products and inverses (see \cite[Theorem 5.18]{Resende2007} and \cite[Proposition 6.6]{BiceStarling2018}).

\begin{prp}\label{OpenInvertibles}
If $C$ is a topological category whose source and range are open maps then its core $C^\times$ is open and hence an \'etale groupoid.
\end{prp}

\begin{proof}
As the product is continuous, for any $g\in C^\times$, we have open $O,N\ni g$ and open $O',N'\ni g^{-1}$ such that $OO'\cup N'N\subseteq C^0(=\mathsf{s}[G^0]=\mathsf{r}[G^0])$.  It follows that
\[g\in O\cap N\cap\mathsf{s}^{-1}[\mathsf{r}[O']]\cap\mathsf{r}^{-1}[\mathsf{s}[N']]\subseteq C^\times,\]
showing that $C^\times$ is open.  Indeed, if $h\in O\cap\mathsf{s}^{-1}[\mathsf{r}[O']]$ then we have $h'\in O'$ with $\mathsf{s}(h)=\mathsf{r}(h')$ and hence $hh'\in OO'\subseteq C^0$.  If $h\in N\cap\mathsf{r}^{-1}[\mathsf{s}[N']]$ too then, likewise, we have $h''\in N'$ with $h''h\in C^0$ and hence $h'=h''hh'=h''=h^{-1}$, so $h\in C^\times$.

Moreover, given any other open $M\ni g^{-1}$,
\[g\in O\cap N\cap\mathsf{s}^{-1}[\mathsf{r}[O'\cap M]]\cap\mathsf{r}^{-1}[\mathsf{s}[N'\cap M]]\subseteq M^{-1},\]
showing that the inverse $g\mapsto g^{-1}$ is also continuous on $C^\times$.  Thus $C^\times$ is a topological groupoid and, as the source is also an open map on $C^\times$, it follows that $C^\times$ is even an \'etale groupoid, by \cite[Theorem 5.18 1.$\Rightarrow$2.]{Resende2007}.
\end{proof}

In particular, for a topology to make a groupoid $C=C^\times$ \'etale, it suffices that the product is continuous and both the source and range are continuous open maps (which could be viewed as an improvement of \cite[Theorem 5.18 1.$\Rightarrow$2.]{Resende2007}).

Before defining the bundles of interest in the present section, recall that a \emph{functor} $\rho:C\rightarrow D$ between categories is a unital product preserving map, i.e. $\rho[C^0]\subseteq D^0$ and $\rho(ab)=\rho(a)\rho(b)$, for all $(a,b)\in C^2$.  In particular, $(\rho(a),\rho(b))\in D^2$ whenever $(a,b)\in C^2$.  As usual, we call $\rho$ an \emph{isocofibration} if the converse also holds, i.e. $(a,b)\in C^2$ whenever $(\rho(a),\rho(b))\in D^2$.  Equivalently, a functor $\rho:C\rightarrow D$ is an isocofibration if and only if it is injective on $C^0$.

\begin{dfn}
A \emph{category bundle} is a bundle $\rho:C\twoheadrightarrow G$ such that
\begin{enumerate}
\item $C$ is a topological category,
\item $G$ is an \'etale groupoid, and
\item $\rho$ is an isocofibration.
\end{enumerate}
If $\rho$ is also locally injective then we call $\rho$ an \emph{\'etale category bundle}.
\end{dfn}

Note that if $\rho:C\twoheadrightarrow G$ is a category bundle then, as $\rho$ is an isocofibration, for every unit $g\in G^0$, we have a unique unit $1_g\in C^0$ with $\rho(1_g)=g$.

As long as the open subbundle $\rho|_{\rho^{-1}[G^0]}$ has at least one continuous section $a$, then $g\mapsto1_g=\mathsf{s}(a(g))$ is also a continuous section and hence $\rho$ restricted to $C^0$ is a homeomorphism onto $G^0$.  In particular, this is true for all \'etale category bundles.

\begin{prp}\label{EtaleCategoryBundle}
If $\rho:C\twoheadrightarrow G$ is an \'etale category bundle then $C$ is an \'etale category and $\rho$ restricted to $C^0$ is a homeomorphism onto $G^0$.
\end{prp}

\begin{proof}
If $\rho:C\twoheadrightarrow G$ is an \'etale category bundle then we first claim that $C^0$ is open in $C$.  To see this, take any $e\in C^0$.  As $\rho$ is locally injective, we have open $O\ni e$ on which $\rho$ is injective.  Replacing $O$ with $O\cap\rho^{-1}[G^0]$ if necessary, we may assume that $O\subseteq\rho^{-1}[G^0]$.  Then $\rho|_O^{-1}$ is a continuous section of the open subbundle $\rho|_{\rho^{-1}[\rho[O]]}$ and hence so too is $\mathsf{s}\circ\rho|_O^{-1}$.  By \autoref{EtaleBundleOpenRange}, $O'=\mathrm{ran}(\mathsf{s}\circ\rho|_O^{-1})$ is open.  As $e\in O'\subseteq C^0$, this proves the claim.

It follows that the restriction $\rho|_{C^0}$ is still an open map.  Certainly, $\rho|_{C^0}$ is also still continuous.  As $\rho$ is an isocofibration, $\rho|_{C^0}$ is also injective hence a homeomorphism onto its range $G^0$.  Noting that $\mathsf{s}(c)=\rho|_{C^0}^{-1}(\mathsf{s}(\rho(c)))$, for all $c\in C$, we see that the source on $C$ is a composition of continuous locally injective open maps.  It follows that the source and, dually, the range is also a (continuous) locally injective open map, i.e. $C$ is an \'etale category.
\end{proof}

In fact, we will be primarily interested in \'etale category bundles $\rho:C\twoheadrightarrow G$ which also have a continuous \emph{zero section} $0:g\mapsto0_g$, i.e. such that
\[a0_g=0_{\rho(a)g}\qquad\text{and}\qquad0_gb=0_{g\rho(b)},\]
for all $a,b\in C$ with $\mathsf{s}(\rho(a))=\mathsf{r}(g)$ and $\mathsf{s}(g)=\mathsf{r}(\rho(b))$.  When $\rho$ has such a continuous zero section, we call $\rho$ a \emph{zero category bundle}.  In this case, we can define the \emph{support} of any section $a$ of $\rho$ by
\[\mathrm{supp}(a)=\{g\in G:a(g)\neq0_g\}.\]

\begin{prp}\label{OpenClosedSupports}
Continuous sections of \'etale zero category bundles have closed supports, while those of Hausdorff zero category bundles have open supports.
\end{prp}

\begin{proof}
For any continuous section $a$ of a zero category bundle $\rho:C\twoheadrightarrow G$, note
\[\mathrm{supp}(a)=a^{-1}[C\setminus\mathrm{ran}(0)].\]
If $\rho$ is \'etale then $\mathrm{ran}(0)$ is open, by \autoref{EtaleBundleOpenRange}, so its complement is closed and hence so is its preimage, by continuity.  Likewise, if $\rho$ is Hausdorff, $\mathrm{ran}(0)$ is closed, by \autoref{ClosedRange}, so its complement is open and hence so is its preimage.
\end{proof}

Sections of zero category bundles yield more examples of well-structured semimodules.  Denote the arbitrary, slice-supported and unit-valued sections of $\rho$ by
\begin{align*}
\mathcal{A}(\rho)&=\{a:a\text{ is a section of }\rho\}.\\
\mathcal{S}(\rho)&=\{a\in\mathcal{A}(\rho):\mathrm{supp}(a)\text{ is a slice}\}.\\
\mathcal{Z}(\rho)&=\{a\in\mathcal{A}(\rho):\mathrm{ran}(a|_{\mathrm{supp}(a)})\subseteq C^0\}.
\end{align*}
We also define $\Phi^\rho:\mathcal{A}(\rho)\rightarrow\mathcal{S}(\rho)$ by
\[\Phi^\rho(a)(g)=\begin{cases}a(g)&\text{if }g\in G^0,\\ 0_g&\text{ otherwise}.\end{cases}\]

\begin{thm}\label{SliceSectionSemigroups}
If $\rho:C\twoheadrightarrow G$ is a zero category bundle then $(\mathcal{A}(\rho),\mathcal{S}(\rho),\mathcal{Z}(\rho),\Phi^\rho)$ is a well-structured semimodule under the product given by
\[ab(f)=\begin{cases}a(g)b(h)&\text{if $f=gh$, for some $g\in\mathrm{supp}(a)$ and }h\in\mathrm{supp}(b),\\0_{f}&\text{ otherwise}.\end{cases}\]
\end{thm}

\begin{proof}
To see that the product well-defined on $\mathcal{S}(\rho)\times\mathcal{A}(\rho)\cup\mathcal{A}(\rho)\times\mathcal{S}(\rho)$, say that $a\in\mathcal{A}(\rho)$, $b\in\mathcal{S}(\rho)$ and $f=gh=ij$, for some $g,i\in\mathrm{supp}(a)$ and $h,j\in\mathrm{supp}(b)$.  Then $\mathsf{s}(h)=\mathsf{s}(gh)=\mathsf{s}(ij)=\mathsf{s}(j)$ so $h=j$, as $\mathrm{supp}(b)$ is a slice, and hence $g=ghh^{-1}=ijj^{-1}=i$.  Moreover, as $f=gh=\rho(a(g))\rho(b(h))$ and $\rho$ is an isocofibration, $a(g)b(h)$ is defined and $\rho(a(g)b(h))=f$.  This shows that $ab$ is a well-defined section and, likewise, so is $ba$.   Further note that $\mathrm{supp}(ab)\subseteq\mathrm{supp}(a)\mathrm{supp}(b)$.  As products of slices are again slices, it follows that $ab\in\mathcal{S}(\rho)$ whenever $a,b\in\mathcal{S}(\rho)$.

Arguing as above we also see that, for any $g\in G$ and any three sections $a_1,a_2,a_3\in\mathcal{A}(\rho)$, two of which are in $\mathcal{S}(\rho)$, there can be at most one triple $g_1,g_2,g_3\in G$ with $g=g_1g_2g_3$ and $g_k\in\mathrm{supp}(a_k)$, for $k=1,2,3$.  In that case, both $((a_1a_2)a_3)(g)$ and $(a_1(a_2a_3))(g)$ equal $a_1(g_1)a_2(g_2)a_3(g_3)$ (which is unambiguous by the associativity of the product in $C$), while they are both zero otherwise, showing that the given product on $\mathcal{S}(\rho)\times\mathcal{A}(\rho)\cup\mathcal{A}(\rho)\times\mathcal{S}(\rho)$ is indeed associative.

For any $z\in\mathcal{Z}(\rho)$, $\mathrm{ran}(z|_{\mathrm{supp}(z)})\subseteq C^0$ implies $\mathrm{supp}(z)\subseteq G^0$ and, for any $e\in\mathrm{supp}(z)$, $z(e)$ is the unit of the monoid fibre $\rho^{-1}\{e\}$.  In particular, $z(e)$ is a central idempotent in $\rho^{-1}\{e\}$ and hence $z$ is a central idempotent in $\mathrm{ran}(\Phi^\rho)$.  As $C^0$ is a subcategory of $C$, it follows that $\mathcal{Z}(\rho)$ is a subsemigroup of $\mathsf{E}(\mathsf{Z}(\mathrm{ran}(\Phi^\rho)))$.  To see that $\mathcal{Z}(\rho)$ is binormal just note that, for any $z\in\mathcal{Z}(\rho)$ and $a,b\in\mathcal{S}(\rho)$ with $ab\in\mathcal{Z}(\rho)$, $\mathrm{ran}((azb)|_{\mathrm{supp}(azb)})\subseteq\mathrm{ran}((ab)|_{\mathrm{supp}(ab)})\subseteq C^0$ and hence $azb\in\mathcal{Z}(\rho)$.

Next, for any $a,b\in\mathcal{A}(\rho)$ and $e\in G^0$, note that
\[(\Phi^\rho(a\Phi^\rho(b)))(e)=(\Phi^\rho(a)\Phi^\rho(b))(e)=(\Phi^\rho(\Phi^\rho(a)b))(e)=a(e)b(e).\]
As all these functions are supported in $G^0$, they are all equal.  Certainly $\Phi^\rho$ is idempotent so $\Phi^\rho$ is an expectation.  To see that $\Phi^\rho$ is also shiftable, take $a\in\mathcal{A}(\rho)$ and $s\in\mathcal{S}(\rho)$.  If $g\in\mathrm{supp}(\Phi^\rho(sa)s))$ then $(\Phi^\rho(sa)s)(g)=sa(\mathsf{r}(g))s(g)=sas(g)$ and hence $g\in\mathrm{supp}(s)\cap\mathrm{supp}(sas)$.  Conversely, if $g\in\mathrm{supp}(s)\cap\mathrm{supp}(sas)$ then $sas(g)=s(g)a(g^{-1})s(g)=sa(\mathsf{r}(g))s(g)=(\Phi^\rho(sa)s)(g)$ and hence $g\in\mathrm{supp}(\Phi^\rho(sa)s))$.  A dual argument shows that $\mathrm{supp}(s\Phi^\rho(as))=\mathrm{supp}(s)\cap\mathrm{supp}(sas)$ too.  As both $\Phi^\rho(sa)s$ and $s\Phi^\rho(as)$ coincide with $sas$ on their supports, $\Phi^\rho(sa)s=s\Phi^\rho(as)$.  This shows that $\Phi^\rho$ is shiftable.

Similarly, for bistability take any $a,b\in\mathcal{A}(\rho)$ with $a$ or $b$ in $\mathcal{S}(\rho)$.  Then
\[\mathrm{supp}(\Phi^\rho(a)b)=(G^0\cap\mathrm{supp}(a))\mathrm{supp}(b)\cap\mathrm{supp}(ab)\]
and $(\Phi_\rho(a)b)(g)=ab(g)$ when $g\in\mathrm{supp}(\Phi_\rho(a)b)$.  In particular, if $\mathrm{supp}(ab)\subseteq G^0$ then $\mathrm{supp}(\Phi_\rho(a)b)\subseteq G^0$ too.  Also, if $ab$ only takes unit values on its support then the same is true of $\Phi_\rho(a)b$.  These and their dual statements show that $\mathrm{ran}(\Phi^\rho)$ and $\mathcal{Z}(\rho)$ are bistable and hence $(\mathcal{A}(\rho),\mathcal{S}(\rho),\mathcal{Z}(\rho),\Phi^\rho)$ is well-structured.
\end{proof}

We get more interesting well-structured semimodules by utilising the topology of $\rho$ and considering compactly supported continuous sections.  Specifically, let
\begin{align*}
\mathcal{C}(\rho)&=\{a\in\mathcal{A}(\rho):a\text{ is continuous}\}.\\
\mathcal{C}_\mathsf{c}(\rho)&=\{a\in\mathcal{C}(\rho):\mathrm{cl}(\mathrm{supp}(a))\text{ is compact}\}.\\
\mathcal{S}_\mathsf{c}(\rho)&=\{a\in\mathcal{C}(\rho):\mathrm{cl}(\mathrm{supp}(a))\text{ is a compact slice}\}.\\
\mathcal{Z}_\mathsf{c}(\rho)&=\mathcal{C}_\mathsf{c}(\rho)\cap\mathcal{Z}(\rho).\\
\Phi^\rho_\mathsf{c}&=\Phi^\rho|_{\mathcal{C}_\mathsf{c}(\rho)}.
\end{align*}
As noted in \autoref{OpenClosedSupports}, when $\rho$ is also an \'etale bundle, supports of continuous sections are already closed and then $\mathrm{cl}(\mathrm{supp}(a))$ could be replaced by $\mathrm{supp}(a)$ above, in which case $\mathcal{S}_\mathsf{c}(\rho)=\mathcal{C}_\mathsf{c}(\rho)\cap\mathcal{S}(\rho)$.  In this case, $\mathcal{Z}_\mathsf{c}(\rho)$ could also be described as the characteristic functions $\mathbf{1}_K:G\rightarrow C$ of compact clopen $K\subseteq G^0$, where
\[\mathbf{1}_K(g)\begin{cases}1_g&\text{if }g\in K\\0_g&\text{otherwise}.\end{cases}\]
Indeed, if $K$ is compact clopen then certainly $\mathbf{1}_K$ is continuous and hence in $\mathcal{Z}_\mathsf{c}(\rho)$.  Conversely, if $a\in\mathcal{Z}_\mathsf{c}(\rho)$ then $a=\mathbf{1}_K$ where both $K=a^{-1}[C^0]$ and its complement $G\setminus K=a^{-1}[\mathrm{ran}(0)]$ are open, by \autoref{EtaleBundleOpenRange} and \autoref{EtaleCategoryBundle}.

\begin{thm}\label{ContinuousSubactions}
If $\rho:C\twoheadrightarrow G$ is a zero category bundle and $G$ is Hausdorff then $(\mathcal{C}_\mathsf{c}(\rho),\mathcal{S}_\mathsf{c}(\rho),\mathcal{Z}_\mathsf{c}(\rho),\Phi_\mathsf{c}^\rho)$ is a well-structured semimodule (with the product above).
\end{thm}

\begin{proof}
Take $s\in\mathcal{S}_\mathsf{c}(\rho)$ and $a\in\mathcal{C}(\rho)$.  As $G$ is Hausdorff, we have an open slice $O$ containing the compact set $K=\mathrm{cl}(\mathrm{supp}(a))$ (see \cite[Proposition 6.3]{BiceStarling2018}).  Note that the function $g\mapsto g_\mathsf{L}=\mathsf{r}|_O^{-1}(\mathsf{r}(g))$ is continuous on $OG=\mathsf{r}^{-1}[\mathsf{r}[O]]$.  Thus $sa$ (defined as in \autoref{SliceSectionSemigroups}) is also continuous on $OG$, as $sa(g)=s(g_\mathsf{L})a(g_\mathsf{L}^{-1}g)$, for all $g\in OG$.  On the other hand, note that $KG=\mathsf{r}^{-1}[\mathsf{r}[K]]$ is closed subset of $OG$.  As $sa$ takes zero values on the open set $G\setminus KG$, it is also continuous there, i.e. $sa\in\mathcal{C}(\rho)$.  Dually, $as\in\mathcal{C}(\rho)$, showing that
\[\mathcal{S}_\mathsf{c}(\rho)\mathcal{C}(\rho)\cup\mathcal{C}(\rho)\mathcal{S}_\mathsf{c}(\rho)\subseteq\mathcal{C}(\rho).\]

As $G$ and hence $G^0$ is Hausdorff, products of compact subsets of $G$ are again compact.  Thus it follows immediately that $\mathcal{S}_\mathsf{c}(\rho)\mathcal{C}_\mathsf{c}(\rho)\cup\mathcal{C}_\mathsf{c}(\rho)\mathcal{S}_\mathsf{c}(\rho)\subseteq\mathcal{C}_\mathsf{c}(\rho)$ and $\mathcal{S}_\mathsf{c}(\rho)\mathcal{S}_\mathsf{c}(\rho)\subseteq\mathcal{S}_\mathsf{c}(\rho)$, i.e. the action of $\mathcal{S}(\rho)$ on $\mathcal{A}(\rho)$ restricts to an action of $\mathcal{S}_\mathsf{c}(\rho)$ on $\mathcal{C}_\mathsf{c}(\rho)$.  As $G$ is Hausdorff, $G^0$ is closed (argue as in \autoref{ClosedRange} but with $\mathsf{s}$ or $\mathsf{r}$).  As $G^0$ is also open (because $G$ is \'etale), it follows that $\Phi^\rho(a)$ is continuous whenever $a$ is and hence $\mathrm{ran}(\Phi^\rho_\mathsf{c})\subseteq\mathcal{S}_\mathsf{c}(\rho)$.  All the required properties of a well-structured semimodule now follow as in the proof of \autoref{SliceSectionSemigroups}.
\end{proof}

The only problem with the above well-structured semimodules is that $\mathcal{Z}_\mathsf{c}(\rho)$ could be rather small.  Indeed, if $G^0$ is connected then $\mathcal{Z}_\mathsf{c}(\rho)$ is trivial, consisting only of the zero section and, if $G^0$ is compact, the one section too.  Less trivial well-structured semimodules arise when $G^0$ and hence $G$ is also compactly based (meaning the topology has a basis of compact open sets).

\begin{dfn}
An \emph{ample groupoid} is a compactly based \'etale groupoid.
\end{dfn}

Note that if $\rho:C\twoheadrightarrow G$ is an \'etale category bundle and $\rho[C^\times]=G$ then $C^\times$ is an ample groupoid if and only if $G$ is.  If $G$ is also Hausdorff and $\rho$ also has a zero section, as above, then we will call $\rho$ an `ample category bundle'.

\begin{dfn}
If $\rho:C\twoheadrightarrow G$ is an \'etale zero category bundle and $G$ is a Hausdorff ample groupoid with $\rho[C^\times]=G$ then we call $\rho$ an \emph{ample category bundle}.
\end{dfn}

So any ample category bundle $\rho:C\twoheadrightarrow G$ has an ample supply of compact clopen characteristic sections, i.e. elements of $\mathcal{Z}_\mathsf{c}(\rho)$.  In fact, ample category bundles also have an ample supply of core-valued continuous sections supported on compact slices.  These are precisely the $\mathcal{Z}_\mathsf{c}(\rho)$-invertible elements of $\mathcal{S}_\mathsf{c}(\rho)$ -- see \autoref{Inverses} below.

\begin{prp}
If $\rho:C\twoheadrightarrow G$ is a ample category bundle then, for every neighbourhood $O$ of every $g\in G$, we have $a\in\mathcal{S}_\mathsf{c}(\rho)$ such that
\[g\in\mathrm{supp}(a)\subseteq O\qquad\text{and}\qquad\mathrm{ran}(a|_{\mathrm{supp}(a)})\subseteq C^\times.\]
\end{prp}

\begin{proof}
Take any open neighbourhood $O$ of any $g\in G$.  As $\rho[C^\times]=G$, we have $c\in C^\times$ with $\rho(c)=g$.  As $C$ is \'etale, we have a neighbourhood $N$ of $c$ on which $\rho$ is locally injective.  As $G$ is ample and $O\cap\rho[N\cap C^\times]$ is neighbourhood of $g$, we have a compact clopen slice $K$ with $g\in K\subseteq O\cap\rho[N\cap C^\times]$.  Thus we can define $a\in\mathcal{S}_\mathsf{c}(\rho)$ with $g\in\mathrm{supp}(a)=K\subseteq O$ and $\mathrm{ran}(a|_{\mathrm{supp}(a)})=\rho^{-1}[K]\cap N\subseteq C^\times$ by
\[a(g)=\begin{cases}\rho|_N^{-1}(g)&\text{if }g\in K\\0_g&\text{otherwise}.\end{cases}\qedhere\]
\end{proof}

To summarise, so far we know that any ample category bundle $\rho:C\twoheadrightarrow G$ gives rise to a well-structured semigroup $(\mathcal{S}_\mathsf{c}(\rho),\mathcal{Z}_\mathsf{c}(\rho),\Phi_\mathsf{c}^\rho)$, thanks to \autoref{ContinuousSubactions}.  Our goal in the following sections will be to reverse this process, starting with a well-structured semigroup which we endeavour to represent on an \'etale category bundle.  The bundle here will be constructed from certain ultrafilters and equivalence classes, which are in turn defined from certain relations on the semigroup.  To ensure that the bundle is ample and that the representation is faithful and surjective, we will need some extra conditions -- this leads to the notion of a Steinberg semigroup, which then provides a purely algebraic characterisation of the well-structured semigroups arising from ample category bundles.

\begin{rmk}\label{TwistedRemark}
While general category bundles do not appear to have been given much attention until now, previous work has dealt with certain \emph{groupoid bundles}, i.e. category bundles $\rho:C\twoheadrightarrow G$ where $C$ is also a groupoid.  For example, in the C*-algebra and Steinberg algebra literature these appear as `twists', which can be viewed as groupoid bundles together with an extra group action.

Specifically, for any topological monoid $T$, let us call $C$ a \emph{$T$-category} if $C$ is a topological category with a continuous action of $T$ on $C$ which is compatible with the product in $C$, i.e. for all $t\in T$ and $b,c\in C$,
\[(tb)c=t(bc)=b(tc).\]
For example, any category $C$ can be considered as a $T$-category where $T$ is the singleton monoid $\{1\}$ and $1c=c$, for all $c\in C$.  By a \emph{category twist} we mean a category bundle $\rho:C\twoheadrightarrow G$ together with a topological monoid $T$ such that $C$ is a $T$-category and $T\times G^0$ is homeomorphic to $\rho^{-1}[G^0]$ via the map
\[(t,g)\mapsto tg.\]
When $C$ is a groupoid, we refer to \emph{$T$-groupoids} and \emph{groupoid twists} respectively.

For example, the discrete twists in \cite[Definition 4.1]{ArmstrongClarkCourtneyLinMcCormickRamagge2022} are precisely the \'etale groupoid twists (i.e. groupoid twists when the bundle map is locally injective) where $T$ is a subgroup of $R^\times$, for some discrete ring $R$.  On the other hand, the twists in \cite[Definition 5.1.1]{Sims2018} are precisely the groupoid twists where $T$ is the unit circle in $\mathbb{C}$, in the usual subspace topology, with the additional requirement that there are local sections going through each point of $C$.

For simplicity, we restrict our attention in the present paper to ($\{1\}$-)categories, although everything could be be extended to $T$-categories without difficulty and then applied to the twists considered elsewhere.
\end{rmk}

\section{Relations}\label{Relations}

In the theory of inverse semigroups, their natural partial order plays a crucial role.  In more general well-structured semigroups, the natural order splits into two different but equally important transitive relations, which we dub `restriction' and `domination'.  Here we investigate their basic properties.

While our primary interest is in well-structured semigroups, throughout the next few sections it will be instructive to work with more general \emph{semigroup inclusions}, namely triples $(S,Z,D)$ where $S$ is a semigroup with subsemigroups $Z$ and $D$ with $Z\subseteq D$.  This will allow us to compare a number of related notions and will also make our theory applicable to the quasi-Cartan pairs we consider later in \autoref{QuasiCartanPairs}.

Accordingly, in this section we make the following standing assumption.
\begin{center}
\textbf{$(S,Z,D)$ is a semigroup inclusion (where $Z\subseteq D\subseteq S$).}
\end{center}
(however, if the reader is not interested in quasi-Cartan pairs or the precise conditions needed for each result then they may just assume throughout that $(S,Z,D)$ is a structured semigroup).
We define \emph{restriction} $\leq$ and \emph{domination} $<$ on $S$ by
\begin{align}
\tag{Restriction}\label{Restriction}a\leq b\qquad&\Leftrightarrow\qquad\exists y,z\in Z\ (yb=ya=a=az=bz).\\
\tag{Domination}a<b\qquad&\Leftrightarrow\qquad\exists s\in S\ (a<_sb)\text{ where}\\
\nonumber a<_sb\qquad&\Leftrightarrow\qquad as,sa\in D,\ bs,sb\in Z\text{ and }bsa=a=asb.
\end{align}

Our motivation here comes from sections of ample category bundles where domination and restriction have natural characterisations in terms of their supports -- see \eqref{AmpleRestriction} and \eqref{AmpleDomination} below.  Special cases of these relations have also appeared in the literature in various contexts.  For example, the Mitsch order (see \cite{Mitsch1986} and \cite{Higgins1994}) is the special case of restriction when $Z=S$, at least when $S$ is a monoid.  On the other hand, if $Z\subseteq\mathsf{E}(S)$ then $a\leq b$ reduces to $a\in bZ\cap Zb$.  In particular, when $Z=\mathsf{E}(S)$, restriction is the Nambooripad order on $S$.  When $D\subseteq\mathsf{E}(S)$, domination is also a subrelation of the Nambooripad order and, when $D=Z=\mathsf{E}(S)$ and $S$ is regular (i.e. $a\in aSa$, for all $a\in S$), all these relations coincide.

\begin{prp}
Assume $S$ is regular and $\mathsf{E}(S)\subseteq Z$.  Then, for all $a,b\in S$,
\[a\leq b\qquad\Leftrightarrow\qquad a\in b\mathsf{E}(S)\cap\mathsf{E}(S)b\qquad\Rightarrow\qquad a<b.\]
Moreover, they are all equivalent when $S$ is regular and $D=Z=\mathsf{E}(S)$.
\end{prp}

\begin{proof}
As $\mathsf{E}(S)\subseteq Z$, we immediately see that $a\in b\mathsf{E}(S)\cap\mathsf{E}(S)b$ implies $a\leq b$.  Conversely, assume $S$ is also regular and $a\leq b$, so we have $y,z\in Z$ with
\[yb=ya=a=az=bz.\]
Taking any $s\in S$ with $asa$, we see that $asyasy=asasy=asy\in\mathsf{E}(S)$ and $asyb=asya=asa=a$.  Likewise, we see that $zsa\in\mathsf{E}(S)$ and $a=azsa=bzsa$, showing that $a\in b\mathsf{E}(S)\cap\mathsf{E}(S)b$.  By assumption, $Z\subseteq D$ and hence $\mathsf{E}(S)\subseteq D$ which implies $a<_{zsy}b$, seeing as $bzsy=azsy=asy\in\mathsf{E}(S)\ni zsa=zsya=zsyb$ and $a=asa=azsyb=bzsya$.  This completes the proof of the first part, while the converse of the last implication is immediate when $D\subseteq\mathsf{E}(S)$, as noted above.
\end{proof}

Now we proceed to examine each of these relations in more detail.

\subsection{Restriction}\label{RestrictionSubsection}

As in \cite[\S8]{Bice2022}, for any $T\subseteq S$ let
\begin{align*}
T^Z&=\{z\in Z:\exists t\in T\ (tz=t)\}.\\
{}^ZT&=\{z\in Z:\exists t\in T\ (zt=t)\}.
\end{align*}
Now \eqref{Restriction} can be rewritten as
\[a\leq b\qquad\Leftrightarrow\qquad a\in{}^Z\hspace{-1pt}ab\cap ba^Z.\]
Using this characterisation, we can show that $\leq$ is almost a partial order relation.

\begin{prp}\label{RestrictionTransitiveAntisymmetric}
The restriction relation $\leq$ is both transitive and antisymmetric.
\end{prp}

\begin{proof}
First note that $a^Z$ is a subsemigroup of $Z$, for any $a\in S$.  Indeed, $y,z\in a^Z$ implies $ayz=az=a$ and hence $yz\in a^Z$.  Likewise, ${}^Z\hspace{-1pt}a$ is a subsemigroup of $Z$, for any $a\in S$.  Also, for any $a,b\in S$, note $b^Z\subseteq(ab)^Z$ and ${}^Z\hspace{-1pt}b\subseteq{}^Z\hspace{-1pt}(ba)$ and hence
\[a\leq b\qquad\Rightarrow\qquad a\in Zb\cap bZ\qquad\Rightarrow\qquad b^Z\subseteq a^Z\quad\text{and}\quad{}^Z\hspace{-1pt}b\subseteq{}^Z\hspace{-1pt}a.\]
It follows that if $a\leq b\leq c$ then
\[a\in{}^Z\hspace{-1pt}ab\cap ba^Z\subseteq{}^Z\hspace{-1pt}a{}^Z\hspace{-1pt}bc\cap cb^Za^Z\subseteq{}^Z\hspace{-1pt}a{}^Z\hspace{-1pt}ac\cap ca^Za^Z\subseteq{}^Z\hspace{-1pt}ac\cap ca^Z,\]
i.e. $a\leq c$, showing that $\leq$ is transitive.

For antisymmetry, just note that if $a\leq b\leq a$ then we have $y,z\in Z$ with $a=ya=yb$ and $b=bz=az=ybz=yb=a$.
\end{proof}

For any $T\subseteq S$ and relation $\sqsubset\ \subseteq S\times S$, let
\[T^\sqsubset=\{s\in S:\exists t\in T\ (t\sqsubset s)\}.\]
For example, $S^\geq$ denotes the elements of $S$ that are the restriction of some other element.  We immediately see that
\[s\in S^\geq\qquad\Rightarrow\qquad s\in Zs\cap sZ\qquad\Rightarrow\qquad s\leq s\qquad\Rightarrow\qquad s\in S^\geq.\]
In other words, $S^\geq$ can also be characterised as those elements of $S$ which have some local left and right units in $Z$, or the elements on which $\leq$ is reflexive:
\[S^\geq=\{s\in S:s\leq s\}=\{s\in S:s\in Zs\cap sZ\}.\]
The following is then immediate from \autoref{RestrictionTransitiveAntisymmetric}.

\begin{prp}\label{RestrictionPoset}
The restriction relation $\leq$ is a partial order on $S^\geq$.
\end{prp}

On the other hand, $S^>$ denotes the dominated elements of $S$.  Note that
\[S^>\subseteq S^\geq,\]
as $a<_{c'}c$ implies that $cc'a=a=ac'c$ and $cc',c'c\in Z$.  In particular, restriction is a partial order on dominated elements.  We also have simpler one-sided characterisations of restriction on $S^>$, at least under suitable extra hypotheses.

\begin{prp}
If $Z$ is binormal and $Z\subseteq\mathsf{Z}(D)$ then, for all $a,b\in S^>$,
\begin{equation}\label{OneSidedRestriction}
a\leq b\qquad\Leftrightarrow\qquad a\in ba^Z\qquad\Leftrightarrow\qquad a\in{}^Z\hspace{-1pt}ab.
\end{equation}
\end{prp}

\begin{proof}
Say $b<_{c'}c$ and $a\in ba^Z$, so we have $z\in Z$ with $a=az=bz$.  If $Z$ is binormal then $czc'\in cZc'\subseteq Z$.  If $Z\subseteq\mathsf{Z}(D)$ too then $czc'b=cc'bz=bz=a$ and also $czc'a=czc'bz=az=a$, showing that $czc'\in{}^Z\hspace{-1pt}a$ so $a=czc'b\in{}^Z\hspace{-1pt}ab$ and hence $a\leq b$.  This the proves the first $\Leftarrow$, while the reverse implication is immediate.  The second equivalence follows by a dual argument.
\end{proof}

If all elements are also dominated, it follows that $\leq$ is invariant under products.

\begin{cor}
If $S=S^>$, $Z$ is binormal and $Z\subseteq\mathsf{Z}(D)$ then, for all $a,b,c\in S$,
\[\tag{Invariance}\label{Invariance}a\leq b\qquad\Rightarrow\qquad ac\leq bc\quad\text{and}\quad ca\leq cb.\]
\end{cor}

\begin{proof}
If $a\leq b$ then $a\in{}^Z\hspace{-1pt}ab$ so $ac\in{}^Z\hspace{-1pt}abc\subseteq{}^Z\hspace{-1pt}(ac)bc$ and hence $ac\leq bc$, by \eqref{OneSidedRestriction}, as long as $S=S^>$, $Z$ is binormal and $Z\subseteq\mathsf{Z}(D)$.  Then, dually, $ca\leq cb$ as well.
\end{proof}

We call $a\in S$ a \emph{normaliser} of $T\subseteq S$ if $aT=Ta$.  Note that normalisers are closed under products, as $aT=Ta$ and $bT=Tb$ implies $abT=aTb=Tab$.  Thus the normalisers of any fixed $T\subseteq S$ form a subsemigroup of $S$ which we denote by
\[\tag{Normalisers}T^\mathsf{N}=\{a\in S:aT=Ta\}.\]
We are primarily concerned with normalisers of $Z$.

\begin{prp}
If $Z\subseteq Z^\mathsf{N}$ then restrictions of normalisers are normalisers, i.e.
\begin{equation}\label{ZNRestriction}
Z\subseteq Z^\mathsf{N}\qquad\Rightarrow\qquad Z^{\mathsf{N}\geq}\subseteq Z^\mathsf{N}.
\end{equation}
\end{prp}

\begin{proof}
If $a\leq b\in Z^\mathsf{N}\supseteq Z$ then, in particular, we have $z\in Z$ with $a=bz$ and hence $aZ=bzZ=bZz=Zbz=Za$, i.e. $a\in Z^\mathsf{N}$.
\end{proof}

We define the \emph{commutant} of any $T\subseteq S$ by
\[\tag{Commutant}T^\mathsf{C}=\bigcap_{t\in T}t^\mathsf{N}.\]
Note that the commutant is contained in the normaliser semigroup and, conversely, any idempotent normaliser is in the commutant, i.e.
\[\mathsf{E}(T^\mathsf{N})\subseteq T^\mathsf{C}\subseteq T^\mathsf{N}.\]
Indeed, if we have $a,b\in S$, $e\in\mathsf{E}(S)$ and $t\in T$ with $et=ae$ and $te=eb$ then
\[et=ae=aee=ete=eeb=eb=te,\]
showing that $e\in T^\mathsf{C}$.

On the commutant of $Z$, we have a stronger characterisation of restriction.

\begin{prp}
For all $a\in Z^\mathsf{C}$ and $b\in S$,
\begin{equation}\label{ZCleq}
a\leq b\qquad\Leftrightarrow\qquad\exists z\in Z\ (zb=za=a=az=bz).
\end{equation}
\end{prp}

\begin{proof}
If $Z^\mathsf{C}\ni a\leq b$, then we have $y,z\in Z$ with $yb=ya=a=az=bz$ and hence
\[zyb=zya=za=az=a=ya=ay=azy=bzy.\]
Thus $zy\in ZZ\subseteq Z$ witnesses \eqref{ZCleq}.
\end{proof}

We also have an even stronger characterisation of restriction on the idempotents of $Z$.  If these commute then they form a meet semilattice.

\begin{prp}
For all $a\in S$ and $z\in\mathsf{E}(Z)$,
\begin{equation}\label{zyyyz}
a\leq z\qquad\Leftrightarrow\qquad az=za=a\in\mathsf{E}(Z).
\end{equation}
If $\mathsf{E}(Z)$ is commutative then products in $\mathsf{E}(Z)$ are meets w.r.t. $\leq$.
\end{prp}

\begin{proof}
If $az=za=a\in\mathsf{E}(Z)$ then $a=aa$ itself witnesses $a\leq z$.  Conversely, if $a\leq z\in\mathsf{E}(Z)$ then we have $w,x\in Z$ with $wz=wa=a=ax=zx\in ZZ\subseteq Z$.  Thus $za=zzx=zx=a=wz=wzz=az$ and $aa=azx=ax=a$, proving \eqref{zyyyz}.

Now take any $y,z\in\mathsf{E}(Z)$.  If $\mathsf{E}(Z)$ is commutative then $yzyz=yyzz=yz$ and hence $yz\in\mathsf{E}(Z)$.  Also $yyz=yz=yzz$ and hence $yz\leq y,z$, by \eqref{zyyyz}.  Moreover, if $x\leq y,z$ then $\mathsf{E}(Z)\ni x=xx=xy=xz$, by what we just proved, so $xyz=xz=x$ and hence $x\leq yz$, by \eqref{zyyyz}.  This shows that $yz$ is indeed the meet of $y$ and $z$ with respect to the restriction relation $\leq$ on $S$.
\end{proof}

If $Z=\mathsf{E}(Z)$ then we immediately see that
\[a\leq b\qquad\Leftrightarrow\qquad a\in Zb\cap bZ.\]
In particular this holds if $Z$ is commutative and up-directed.

\begin{prp}\label{leqUpDirected}
If $Z$ is commutative and up-directed w.r.t. $\leq$ then $Z=\mathsf{E}(Z)$.
\end{prp}

\begin{proof}
If $Z$ is commutative and up-directed w.r.t. $\leq$ then, in particular, $Z\subseteq S^\geq$.  For any $x\in Z$, this means we have $y\in Z$ with $x=xy$.  By up-directedness, we then have $z\in Z$ with $x,y\leq z$, so we have $x',y'\in Z$ with $x=xx'=zx'$ and $y=yy'=zy'$.  Then $x=xy=xyy'=xy'$ too and hence
\[x=xy=xy'z=xz=xx'z=xx.\qedhere\]
\end{proof}

\subsection{Domination}

The first thing to note is that domination is transitive and also in an auxiliary relationship to restriction, at least when $Z$ is commutative.
\begin{prp}
For all $a,b,c,s,t\in S$,
\begin{align}
\tag{Transitivity}\label{Transitivity}a<b<_tc\qquad&\Rightarrow\qquad a<_tc.\\
\tag{Left-Auxiliarity}\label{LeftAuxiliarity}a\leq b<_tc\qquad&\Rightarrow\qquad a<_tc.\\
\intertext{If $Z$ is commutative then, moreover, for all $a,b,c\in S$,}
\tag{Right-Auxiliarity}\label{RightAuxiliarity}a<b\leq c\qquad&\Rightarrow\qquad a<c.
\end{align}
\end{prp}

\begin{proof}
For \eqref{Transitivity} and \eqref{LeftAuxiliarity}, it suffices to show that
\[a\in bD\cap Db\quad\text{and}\quad b<_tc\qquad\Rightarrow\qquad a<_tc.\]
To see this, just note $a\in bD\cap Db$ and $b<_tc$ implies $ta\in tbD\subseteq DD\subseteq D$, $at\in Dbt\subseteq DD\subseteq D$ and, taking $d,e\in D$ with $db=a=be$,
\[cta=ctbe=be=a=db=dbtc=atc.\]

For \eqref{RightAuxiliarity}, assume $Z$ is commutative and $a<_sb\leq c$, so we have $y,z\in Z$ with $zc=zb=b=by=cy$.  We claim that $a<_{ysz}c$.  To see this, note that $yszc=ysb=sby=sb\in Z$, by the commutativity of $Z$.  Likewise, $cysz=bs\in Z$.  Also $aysz=asbysz=asbsz=aszbs=asbs=as\in D$ and $ayszc=asb=a$.  Likewise, $ysza\in D$ and $cysza=a$, proving the claim.
\end{proof}

Under suitable hypotheses, we have a more one-sided characterisation of $<$.

\begin{prp}
If $D$ is trinormal and $Z\subseteq\mathsf{Z}(D)$ then
\begin{equation}\label{OneSidedDomination}
a<_{b'}b\qquad\Leftrightarrow\qquad ab'\in D,\ b'b,bb'\in Z\text{ and }a=ab'b.
\end{equation}
\end{prp}

\begin{proof}
Assume $D$ is trinormal, $Z\subseteq\mathsf{Z}(D)$, $ab'\in D$, $b'b,bb'\in Z$ and $a=ab'b$.  Then $bb'a=bb'ab'b=ab'bb'b=a$ and $ab'=ab'bb'=bb'ab'$, as $b'b\in Z\subseteq\mathsf{Z}(D)$, and hence $b'a=b'ab'b\in D$, as $D$ is trinormal, showing that $a<_{b'}b$.
\end{proof}

Next we observe that domination is invariant under products on the left with $D$ or on the right with $Z$, again under suitable extra hypotheses.

\begin{prp}
If $Z\subseteq\mathsf{Z}(D)$ then, for all $a,b,b'\in S$, $d\in D$ and $z\in Z$,
\begin{align*}
\tag{$D$-Invariance}\label{DInvariance}D\text{ is trinormal and }a<_{b'}b\quad&\Rightarrow\quad ad<_{b'}b.\\
\tag{$Z$-Invariance}\label{ZInvariance}Z\text{ is binormal and }az=a<_{b'}b\quad&\Rightarrow\quad a<_{b'}bz\text{ and }a<_{zb'}b.
\end{align*}
\end{prp}

\begin{proof}
Assume $Z\subseteq\mathsf{Z}(D)$ and $a<_{b'}b$.  If $d\in D$ then $b'ad\in DD\subseteq D$ and $b'bb'ad=b'ad=b'adb'b$, as $b'b\in Z\subseteq\mathsf{Z}(D)$ so $adb'=bb'adb'\in D$, if $D$ is trinormal.  Also $bb'ad=ad=ab'bd=adb'b$, showing that $ad<_{b'}b$.

If $z\in Z\subseteq\mathsf{Z}(D)$ then $zb'b=b'bz\in ZZ\subseteq Z$.  If $Z$ is also binormal then $bzb'\in Z$.  Also $ab'bz=az=bb'az=bzb'a$ so if $a=az$ then $a<_{b'}bz$ and $a<_{zb'}b$.
\end{proof}

In structured semigroups, domination also respects products.

\begin{prp}
If $(S,Z,D)$ a stuctured semigroup and $a,b,b',c,d,d'\in S$ then
\[\tag{Multiplicativity}\label{Multiplicativity}a<_{b'}b\quad\text{and}\quad c<_{d'}d\qquad\Rightarrow\qquad ac<_{d'b'}bd.\]
\end{prp}

\begin{proof}
Assume $(S,Z,D)$ a stuctured semigroup, $a<_{b'}b$ and $c<_{d'}d$.  Then $d'b'bd,bdd'b'\in Z$, as $Z$ is binormal.  Also $b'acd'b'b=b'bb'acd'=b'acd'\in DD\subseteq D$, as $b'b\in Z\subseteq\mathsf{Z}(D)$, and hence $acd'b'=bb'acd'b'\in D$, as $D$ is trinormal.  As $cd'\in D$ commutes with $b'b\in Z$, also $acd'b'bd=ab'bcd'd=ac$ so $ac<_{d'b'}bd$, by \eqref{OneSidedDomination}.
\end{proof}

Unlike restriction, $<$ is not always antisymmetric.  Indeed, antisymmetry is only guaranteed to hold on idempotents, where a formula like \eqref{zyyyz} holds for $<$.

\begin{prp}
$<$ is antisymmetric on $\mathsf{E}(S)$.  Moreover, for all $a\in S$.
\begin{equation}\label{EZdomination}
a<e\in\mathsf{E}(S)\qquad\Rightarrow\qquad ea=a=ae.
\end{equation}
The converse also holds if $a\in D$ and $e\in\mathsf{E}(Z)$.
\end{prp}

\begin{proof}
If $a<_{e'}e\in\mathsf{E}(S)$ then $ea=eee'a=ee'a=a=ae'e=ae'ee=ae$.  In particular, if $d,e\in\mathsf{E}(S)$ and $d<e<d$ then $d=de=e$, i.e. $<$ is antisymmetric on $\mathsf{E}(S)$.  For the converse of \eqref{EZdomination}, just note that if $a\in D$, $e\in\mathsf{E}(Z)$ and $ea=a=ae$ then $eea=ea=a=ae=aee\in D$ and $ee=e\in Z$ so $a<_ee$.
\end{proof}

Lastly let us note that binormality and normality are equivalent conditions on $Z$ when $Z$ is central in $D$ and all elements are dominated, i.e. $S=S^>$.

\begin{prp}\label{BiNormality}
If $Z$ is binormal then
\[D\subseteq Z^\mathsf{N}\qquad\Rightarrow\qquad S^>\subseteq Z^\mathsf{N}.\]
\end{prp}

\begin{proof}
If $Z$ is binormal, $D\subseteq Z^\mathsf{N}$ and $a<_{s'}s$ then $aZ=ss'aZ=sZs'a\subseteq Za$ and, dually, $Za\subseteq aZ$, showing that $a\in Z^\mathsf{N}$.
\end{proof}

Before moving on, let us note that domination allows us to define `duals' which will useful in the next subsection and also key to defining inverses in the (ultra)filter groupoids considered below.  Specifically, the \emph{dual} of any $T\subseteq S$ is given by
\[T^*=\{a\in S:\exists t\in T\ \exists s\in S\ (t<_as)\}.\]
In particular, we let $a^*=\{a\}^*$ denote the dual of any $a\in S$.

\subsection{Orthogonality}

As usual, we call $0\in S$ a \emph{zero} if
\[\tag{Zero}S0=\{0\}=0S.\]
Throughout the rest of this section we assume that $S$ has a (necessarily unique) zero.  This allows us to define the \emph{orthogonality relation} $\perp$ on $S$ by
\begin{align*}
a\perp b\qquad&\Leftrightarrow\qquad\exists y,z\in Z\ (ya=a=az\text{ and }yb=0=bz).\\
&\Leftrightarrow\qquad\exists z\in\!{}^Z\hspace{-1pt}a\,(zb=0)\text{ and }\exists z\in a^Z(bz=0).
\end{align*}
As with the restriction relation $\leq$, the orthogonality relation $\perp$ will play a key role later on when it comes to Steinberg semigroups and Steinberg rings.  For Steinberg semigroups coming from ample category bundles, orthogonality also has a simple characterisation in terms of supports -- see \eqref{AmpleOrthogonality} below.

Here we briefly examine what can be said for the more general semigroup inclusions we are considering in this section.  First we note orthogonality respects the product and has an auxiliary relationship both to restriction and domination.

\begin{prp}
For all $a,b,c,d\in S$,
\begin{align}
\label{abperpcd}a\perp b\quad\text{and}\quad c\perp d\qquad&\Rightarrow\qquad ac\perp bd.\\
\label{perpAux}a\leq b\perp c\geq d\quad\text{or}\quad a<b\perp c>d\qquad&\Rightarrow\qquad a\perp d.
\end{align}
\end{prp}

\begin{proof}
If $a\perp b$ and $c\perp d$ then we have $y\in{}^Z\hspace{-1pt}a\subseteq{}^Z\hspace{-1pt}\{ac\}$ with $0=yb=ybd$ and $z\in c^Z\subseteq\{ac\}^Z$ with $0=dz=bdz$, showing that $ac\perp bd$.  This proves \eqref{abperpcd}.

Similarly, note that if $a\in bS\cap Sb$ then ${}^Z\hspace{-1pt}b\subseteq{}^Z\hspace{-1pt}a$ and $b^Z\subseteq a^Z$, in which case $b\perp c$ immediately implies $a\perp c$.  In particular,
\[a\leq b\perp c\quad\text{or}\quad a<b\perp c\qquad\Rightarrow\qquad a\perp c.\]
Likewise, if $d\in cS\cap Sc$ then $yc=0=cz$ implies $yd=0=dz$, in which case $b\perp c$ immediately implies $b\perp d$.  In particular,
\[b\perp c\geq d\quad\text{or}\quad b\perp c>d\qquad\Rightarrow\qquad b\perp d.\]
Combining these two observations yields \eqref{perpAux}.
\end{proof}

Lastly, we note orthogonality has a simple characterisation in terms of duals and is also invariant under products, at least by dominated elements.

\begin{prp}
If $Z\subseteq\mathsf{Z}(D)$ is binormal then, for all $a,b\in S$,
\begin{align}
\label{DualOrthogonality}S^>\ni a\perp b\qquad&\Leftrightarrow\qquad\exists s\in a^*\ (sb=0=bs).\\
\label{OrthoInvariance}s\in S^>\text{ and }a\perp b\qquad&\Rightarrow\qquad as\perp bs\ \ \text{and}\ \ sa\perp sb.
\end{align}
\end{prp}

\begin{proof}
Take $s\in a^*$, so we have $t\in S$ with $a<_st$.  In particular, $ts\in{}^Z\hspace{-1pt}a$ and $st\in a^Z$ so if $sb=0=bs$ and hence $tsb=0=bst$ then $a\perp b$.  Conversely, if $a<_st$ and $a\perp b$ then we have $y,z\in Z$ with $ya=a=az$ and $yb=0=bz$.  Then $a<_{zsy}t$, by \eqref{ZInvariance}, i.e. $zsy\in a^*$.  As $zsyb=0=bzsy$, this proves \eqref{DualOrthogonality}.

Now say $s<_{t'}t$ and $a\perp b$, so again we have $y,z\in Z$ with $ya=a=az$ and $yb=0=bz$.  Then certainly $yas=as$ and $ybs=0$.  As $Z$ is binormal, $t'zt\in Z$.  As $Z\subseteq\mathsf{Z}(D)$, we see that $as(t'zt)=azst't=as$ and $bs(t'zt)=bzst't=0$, showing that $as\perp bs$.  Dually, $sa\perp sb$, proving \eqref{OrthoInvariance}.
\end{proof}

\section{Inverses}\label{Inverses}

Again let us make the following standing assumption.
\begin{center}
\textbf{$(S,Z,D)$ is a semigroup inclusion}
\end{center}
(although again our primary focus later will be on structured semigroups).

For any $T\subseteq S$ and $a\in S$, we call $s\in S$ a \emph{$T$-inverse} of $a$ if
\[\tag{$T$-Inverse}a=asa,\quad s=sas\quad\text{and}\quad as,sa\in T.\]
Note that $a=asa$ implies that $sa=sasa$ and $as=asas$ so $T$-inverses and $\mathsf{E}(T)$-inverses are the same.  If $S$ has a unit $1$ then, for $s$ to be a $\{1\}$-inverse of $a$, it suffices that $as=1=sa$, i.e. $\{1\}$-inverses are inverses in the usual sense.  On the other hand, for $s$ to be an $S$-inverse of $a$, it suffices that $a=asa$ and $s=sas$.  If every element of $S$ has an $S$-inverse then $S$ is a \emph{regular semigroup} while if these $S$-inverses are also unique then $S$ is an \emph{inverse semigroup} (see \cite{Lawson1998}).

Here we are primarily interested in $Z$-inverses, as these will play an important role when it comes to examining Steinberg semigroups and Steinberg rings.  We denote the set of all $Z$-inverses of all elements of any subset $T\subseteq S$ by
\[T^\dagger=\{s\in S:s\text{ is a $Z$-inverse of some }t\in T\}.\]
In particular, $S^\dagger$ denotes the set of all $Z$-invertible elements, while if $a\in S$ then $a^\dagger=\{a\}^\dagger$ is the set of all $Z$-inverses of $a$.  Equivalently, in terms of domination,
\[a^\dagger=\{s\in S:a<_sa\text{ and }s<_as\}.\]
In particular, $T^\dagger\subseteq T^*$, for all $T\subseteq S$.  Again the situation to keep in mind is when $(S,Z,D)$ arises from an ample category bundle as in \autoref{ContinuousSubactions} (with $D=\mathrm{ran}(\Phi_\mathsf{c}^\rho)$) -- then the characterisation of domination in \eqref{AmpleDomination} implies that the $Z$-invertible sections are precisely those which take only invertible values on their support (and the only $Z$-inverse of such a section is then its pointwise inverse).

\begin{prp}
The domination relation is reflexive precisely on the $Z$-invertible elements, while every restriction of a $Z$-invertible element is again $Z$-invertible, i.e.
\begin{equation}\label{Sdaggergeq}
S^\dagger=\{a\in S:a<a\}=S^{\dagger\geq}.
\end{equation}
\end{prp}

\begin{proof}
If $a\in S^\dagger$ then, in particular, we have $s\in S$ with $a<_sa$.  Conversely, $a<_sa$ implies that $a$ is a $Z$-inverse of $sas$, as $a(sas)a=a$, $(sas)a(sas)=sas$ $(sas)a=sa\in Z$ and $a(sas)=as\in Z$, which shows that $S^\dagger=\{a\in S:a<a\}$.

Now say $a\leq r$ so $yr=ya=a=az=rz$, for some $y,z\in Z$.  If $r<_sr$ then $as=yrs\in ZZ\subseteq Z$, $sa=srz\in ZZ\subseteq Z$ and $asa=yrsrz=yrz=ya=a$.  Thus $a<_sa$ so $a\in S^\dagger$.  Conversely, if $a\in S^\dagger\subseteq S^>\subseteq S^\geq$ then $a\leq a\in S^{\dagger\geq}$.  This shows that $S^\dagger=S^{\dagger\geq}$, which completes the proof of \eqref{Sdaggergeq}.
\end{proof}

However, $Z$-inverses may not be unique or closed under products.  To fix this we restrict further to $\mathsf{E}(Z)^{\mathsf{N}\dagger}$, i.e. to the $Z$-inverses of normalisers of $\mathsf{E}(Z)$.

\begin{prp}
If $\mathsf{E}(Z)$ is commutative then $\mathsf{E}(Z)^{\mathsf{N}\dagger}$ is an inverse semigroup where $\{ab\}^\dagger=b^\dagger a^\dagger$, for all $a,b\in\mathsf{E}(Z)^{\mathsf{N}\dagger}$, with idempotent semilattice given by
\begin{equation}\label{ESdagger}
\mathsf{E}(Z)=\mathsf{E}(\mathsf{E}(Z)^{\mathsf{N}\dagger})=\mathsf{E}(S^{\dagger}).
\end{equation}
\end{prp}

\begin{proof}
Assume $\mathsf{E}(Z)$ is commutative.  If $b,c\in a^\dagger$ then $ab,ac\in\mathsf{E}(Z)$ commute and hence $b=bab=bacab=babac=bac$.  Likewise $b=bab=cab$, showing that $b\leq c$.  Dually, $c\leq b$ and hence $b=c$, by \autoref{RestrictionTransitiveAntisymmetric}, i.e. $Z$-inverses are unique.

If $a\in\mathsf{E}(Z)^\mathsf{N}$ and $s\in a^\dagger$ then
\[s\mathsf{E}(Z)=sas\mathsf{E}(Z)=s\mathsf{E}(Z)as=sa\mathsf{E}(Z)s\subseteq\mathsf{E}(Z)\mathsf{E}(Z)s\subseteq\mathsf{E}(Z)s.\]
Dually, $\mathsf{E}(Z)s\subseteq s\mathsf{E}(Z)$ so $s\in\mathsf{E}(Z)^\mathsf{N}$, showing that $\mathsf{E}(Z)^{\mathsf{N}\dagger}\subseteq\mathsf{E}(Z)^\mathsf{N}$.

Say we also have $b\in\mathsf{E}(Z)^\mathsf{N}$ and $t\in b^\dagger$.  In particular, $sa,bt\in Z$ commute so $abtsab=asabtb=ab$ and, likewise, $tsabts=ts$.  Moreover,
\[abts\in a\mathsf{E}(Z)s=\mathsf{E}(Z)as\subseteq ZZ\subseteq Z\] and, likewise, $tsab\in Z$.  This shows that $\{ab\}^\dagger=\{ts\}=b^\dagger a^\dagger$, showing that $\mathsf{E}(Z)^{\mathsf{N}\dagger}$ is a subsemigroup of $S$ and hence an inverse semigroup, as $\mathsf{E}(Z)^{\mathsf{N}\dagger}\subseteq\mathsf{E}(Z)^\mathsf{N}$.

For all $z\in\mathsf{E}(Z)$, note $z^\dagger=\{z\}$ so $\mathsf{E}(Z)=\mathsf{E}(Z)^\dagger\subseteq\mathsf{E}(Z)^{\mathsf{C}\dagger}\subseteq\mathsf{E}(Z)^{\mathsf{N}\dagger}\subseteq S^{\dagger}$.  Conversely, if $e=ee\in s^\dagger$ then $s=ses=sees=esse\in eS$, as $es,se\in Z$ commute.  Thus $s=es$ so $e=ese=se\in Z$, showing that $\mathsf{E}(S^\dagger)\subseteq Z$.  This proves \eqref{ESdagger}.
\end{proof}

Restricting further to the $Z$-invertible normalisers of $Z$, we obtain another inverse semigroup $Z^{\mathsf{N}\dagger}$ that is important for the quasi-Cartan pairs in \cite{ArmstrongCastroClarkCourtneyLinMcCormickRamaggeSimsSteinberg2021}, which we also consider below in \autoref{QuasiCartanPairs}.  Indeed, the following characterisation shows that $Z^{\mathsf{N}\dagger}$ consists precisely of the normalisers considered in \cite{ArmstrongCastroClarkCourtneyLinMcCormickRamaggeSimsSteinberg2021}.

\begin{prp}\label{ZInvertibleNormalisers}
If $Z\subseteq Z^\mathsf{N}$ then $Z^{\mathsf{N}\dagger}$ is an inverse subsemigroup of $\mathsf{E}(Z)^{\mathsf{N}\dagger}$ with the same idempotent idempotent semilattice, i.e. $\mathsf{E}(Z^{\mathsf{N}\dagger})=\mathsf{E}(Z)$, and
\begin{equation}\label{ZInvertibleNormaliserCharacterisation}
Z^{\mathsf{N}\dagger}=\{a\in S:a\in aZ\cap Za\text{ and }\exists s\in S\ (asa=a\text{ and }aZs\cup sZa\subseteq Z)\}.
\end{equation}
Moreover, $Z^{\mathsf{N}\dagger}$ is also invariant under restrictions, i.e. $Z^{\mathsf{N}\dagger\geq}\subseteq Z^{\mathsf{N}\dagger}$.
\end{prp}

\begin{proof}
If $a\in s^\dagger$ then $as,sa\in Z$ and $a=asa\in aZ\cap Za$.  If, moreover, $s\in Z^\mathsf{N}$ then $aZs\cup sZa=asZ\cup Zsa\subseteq ZZ\subseteq Z$.

Conversely, assume $Z\subseteq Z^\mathsf{N}$, $asa=a\in aZ\cap Za$ and $aZs\cup sZa\subseteq Z$.  It follows that $as\in aZs\subseteq Z$ and $sa\in sZa\subseteq Z$ so $a<_sa$ and hence $\{sas\}^\dagger=\{a\}$.  As $aZsas\cup sasZa\subseteq aZZs\cup sZZa\subseteq Z$, we may replace $s$ with $sas$ if necessary and assume that $s$ and $a$ are $Z$-inverses.  As $as,sa\in Z\subseteq Z^\mathsf{N}$, it follows that $sZ=sasZ=sZas\subseteq Zs=Zsas\subseteq saZs\subseteq sZ$, showing that $s\in Z^\mathsf{N}$.  This completes the proof of \eqref{ZInvertibleNormaliserCharacterisation}.  Likewise, $aZ=Za$, which shows that $Z^{\mathsf{N}\dagger}\subseteq Z^\mathsf{N}$.  As $a\mathsf{E}(Z)s\cup s\mathsf{E}(Z)a\subseteq\mathsf{E}(Z)$ and $as,sa\in\mathsf{E}(Z)\subseteq\mathsf{E}(Z^\mathsf{N})\subseteq Z^\mathsf{C}\subseteq\mathsf{E}(Z)^\mathsf{N}$, we see that
\[a\mathsf{E}(Z)=asa\mathsf{E}(Z)=a\mathsf{E}(Z)sa\subseteq\mathsf{E}(Z)a=\mathsf{E}(Z)asa=as\mathsf{E}(Z)a\subseteq a\mathsf{E}(Z).\]
Thus $Z^{\mathsf{N}\dagger}\subseteq\mathsf{E}(Z)^{\mathsf{N}\dagger}$ and hence
\[Z^{\mathsf{N}\dagger}Z^{\mathsf{N}\dagger}\subseteq Z^\mathsf{N}Z^\mathsf{N}\cap\mathsf{E}(Z)^{\mathsf{N}\dagger}\mathsf{E}(Z)^{\mathsf{N}\dagger}\subseteq Z^\mathsf{N}\cap S^\dagger=Z^{\mathsf{N}\dagger}.\]
This shows that $Z^{\mathsf{N}\dagger}$ is indeed an inverse subsemigroup of $\mathsf{E}(Z)^{\mathsf{N}\dagger}$.

To see that $Z^{\mathsf{N}\dagger}$ and $\mathsf{E}(Z)^{\mathsf{N}\dagger}$ have the same idempotents, note that \eqref{ESdagger} yields
\[\mathsf{E}(Z^{\mathsf{N}\dagger})\subseteq\mathsf{E}(S^\dagger)\subseteq\mathsf{E}(Z)\subseteq\mathsf{E}(Z^{\mathsf{N}\dagger}).\]
To see that $Z^{\mathsf{N}\dagger}$ is invariant under restrictions, note that \eqref{ZNRestriction} and \eqref{Sdaggergeq} yield
\[Z^{\mathsf{N}\dagger\geq}\subseteq S^\dagger\cap Z^{\mathsf{N}\geq}\subseteq S^\dagger\cap Z^\mathsf{N}=Z^{\mathsf{N}\dagger}.\qedhere\]
\end{proof}

When they exist, $Z$-inverses always witness domination, i.e.
\begin{equation}\label{ReflexiveDomination}
a<b\in s^\dagger\qquad\Rightarrow\qquad a<_sb.
\end{equation}
Indeed, if $a<b<_sb$ then $a<_sb$, by \eqref{Transitivity}.  On $S^\dagger$, the restriction and orthogonality relations can also be characterised in terms of $Z$-inverses.

\begin{prp}
If $\mathsf{E}(Z)\subseteq\mathsf{Z}(Z)$ then, for all $a,b,s\in S$ with $a\in s^\dagger$,
\begin{align}
\label{Inversealeqb}a\leq b\qquad&\Leftrightarrow\qquad bsa=a=asb\qquad\Leftrightarrow\qquad a<_sb.\\
\intertext{If $S$ also has a zero then}
\label{Inverseaperpb}a\perp b\qquad&\Leftrightarrow\qquad bs=0=sb.
\end{align}
\end{prp}

\begin{proof}
If $\mathsf{E}(Z)\subseteq\mathsf{Z}(Z)$ and $a<_sa\leq b$ then we have $z\in Z$ with $a=az=bz$ so $sa=saz=zsa$ and hence $a=asa=bzsa=bsa$.  Dually, $a=asb$, proving the first $\Rightarrow$.  Now if $bsa=a=asb$ and $s<_as$ then $sb=sasb=sa\in Z$ and, dually, $bs\in Z$ and hence $a<_sb$, proving the second $\Rightarrow$.  Conversely, if $a<_sb$ then certainly $bsa=a=asb$, which in turn implies $a\leq b$ when $a<_sa$, proving \eqref{Inversealeqb}.

Now say $S$ also has a zero.  If $bs=0=sb$ then $bsa=0=asb$ and hence $a\perp b$, as $a=asa$.  Conversely, if $a\perp b$ then we have $y,z\in Z$ with $ya=a=az$ and $yb=0=bz$ so $bs=bsas=bsazs=bzsas=0=sasyb=syasb=sasb=sb$.
\end{proof}

In particular, restriction implies domination on $S^\dagger$ when $\mathsf{E}(Z)\subseteq\mathsf{Z}(Z)$.  Actually, in general it follows from \eqref{LeftAuxiliarity} that, for all $s\in S$,
\[a\leq b\in S^\dagger\qquad\Rightarrow\qquad a<b.\]
Conversely, if both $a$ and $b$ are restrictions of a single element of $S$, then $a<b$ implies $a\leq b$, in which case $a$ is also necessarily $Z$-invertible.

\begin{prp}
For all $a,b,c,s\in S$,
\begin{equation}\label{CommonRestrictionDomination}
a,b\leq c\quad\text{and}\quad a<_sb\qquad\Rightarrow\qquad a<_sa\leq b.
\end{equation}
\end{prp}

\begin{proof}
Say $a,b\leq c$ and $a<_sb$.  Then we have $y,z\in Z$ such that $a=ay=cy$ and $b=zb=zc$ and hence $by=zcy=za=zbsa=bsa=a$.  In particular, $a\in ba^Z$, $sa=sby\in ZZ\subseteq Z$ and $asa=asby=ay=a$.  Likewise, $a\in{}^Z\hspace{-1pt}ab$ and $as\in Z$ and hence $a<_sa\leq b$.
\end{proof}

We call $<$ \emph{interpolative} if
\[\tag{Interpolative}\label{Interpolative}a<b\qquad\Rightarrow\qquad\exists s\in S\ (a<s<b).\]
For example, $<$ is interpolative is when $Z$ is a Cartan subalgebra of a C*-algebra or part of a more general structured C*-algebra -- see \cite[Proposition 6.11]{Bice2021DHKR}.  When $Z$ consists of idempotents, $<$ is again automatically interpolative and we can even choose the interpolating element to be $Z$-invertible.

\begin{prp}
If $Z\subseteq\mathsf{E}(S)$ then, for all $a,b\in S$,
\begin{equation}\label{ReflexiveInterpolation}
a<b\qquad\Rightarrow\qquad\exists s\in S^\dagger\ (a<s<b).
\end{equation}
\end{prp}

\begin{proof}
Just note that if $a<_{b'}b$ and $b'b,bb'\in\mathsf{E}(S)$ then $a<_{b'}bb'b<_{b'}bb'b<_{b'}b$.
\end{proof}

Lastly, we characterise structured semigroups $(S,\mathsf{E}(S),\mathsf{E}(S))$ arising from inverse semigroups $S$, in terms of reflexivity and antisymmetry of the domination relation.

\begin{prp}\label{InverseStructuredSemigroups}
If $Z\subseteq\mathsf{Z}(D)$ and $D$ is diagonal and trinormal then
\begin{align}
\label{Zposet}(Z,<)\text{ is a poset}\quad&\Leftrightarrow\quad Z\subseteq\mathsf{E}(S).\\
\label{Dposet}(D,<)\text{ is a poset}\quad&\Leftrightarrow\quad D\subseteq\mathsf{E}(S^\dagger)\quad\Leftrightarrow\quad D=\mathsf{E}(Z).\\
\label{Sposet}(S,<)\text{ is a poset}\quad&\Leftrightarrow\quad D\subseteq\mathsf{E}(S^\dagger)\text{ and }S=S^>\\
\nonumber&\Leftrightarrow\quad D=\mathsf{E}(Z)\text{ and }S=S^\dagger.
\end{align}
\end{prp}

\begin{proof}
Take $d\in S^\dagger\cap D$, so we have $s\in S$ with $d<_sd$.  We claim that
\[d<ds<d.\]
Indeed, as $sdsd=sd\in Z$ and $dsdsd=dsd=d\in D$, \eqref{OneSidedDomination} yields $d<_{sd}sd$.  Also $sds\in D$, as $D$ is diagonal, so \eqref{OneSidedDomination} again yields $sd<_sd$.

In particular, if $d\in Z$ and $(Z,<)$ is a poset then we have $s$ as above, as $<$ is reflexive on $Z$, and hence $d=ds$, as $<$ is antisymmetric on $Z$.  It follows that $d=dsd=d^2$, showing that $Z\subseteq\mathsf{E}(S)$.  Conversely, if $Z\subseteq\mathsf{E}(S)$ then, by \eqref{EZdomination},
\[y<z\qquad\Leftrightarrow\qquad y=yz,\]
for all $y,z\in Z$, which is immediately seen to be a partial order, proving \eqref{Zposet}.

If $<$ is a partial order on $D$ then, in particular, $<$ is reflexive on $D$ and hence $D\subseteq S^\dagger$.  Then the claim above again implies that $D\subseteq\mathsf{E}(S)$, proving the first $\Rightarrow$ in \eqref{Dposet}.  The second $\Rightarrow$ is immediate from \eqref{ESdagger}.  To complete the cycle of equivalences, just note that $D=\mathsf{E}(Z)$ implies that $(D,<)=(Z,<)$ is a poset, by \eqref{Zposet}.  This finishes the proof of \eqref{Dposet}.

The first $\Rightarrow$ in \eqref{Sposet} follows from \eqref{Dposet}.  Next assume $D\subseteq\mathsf{E}(S^\dagger)$ and $S=S^>$.  For any $s\in S$, this means we have $t,t'\in S$ with $s<_{t'}t$.  Then $st',t's\in D=\mathsf{E}(Z)$, by \eqref{Dposet}, so $st's=st'st't=st't=s$ and hence $s<_{t'}s$, by \eqref{OneSidedDomination}.
Thus $S=S^\dagger$, proving the second $\Rightarrow$ in \eqref{Sposet}.  Finally note that if $D=\mathsf{E}(Z)$ then domination implies restriction, which is antisymmetric, by \autoref{RestrictionTransitiveAntisymmetric}.  By \eqref{Sdaggergeq}, $S=S^\dagger$ then implies domination is also reflexive and hence a partial order.
\end{proof}

\subsection{Supports}

A \emph{source-support} of $a\in S$ is a minimum $\mathsf{s}(a)$ of $a^Z$, i.e.
\[\tag{Source-Support}\mathsf{s}(a)\in a^Z\subseteq\mathsf{s}(a)^\leq.\]
Likewise, a \emph{range-support} of $a\in S$ is a minimum $\mathsf{r}(a)$ of ${}^Z\hspace{-1pt}a$.  As $\leq$ is antisymmetric, supports are unique, when they exist, and they are also necessarily idempotents.

\begin{prp}\label{Esupport}
Every support is an idempotent, i.e. $\mathsf{s}(a),\mathsf{r}(a)\in\mathsf{E}(Z)$.
\end{prp}

\begin{proof}
Say $a^Z$ has minimum $\mathsf{s}(a)$.  In particular, $\mathsf{s}(a)\in a^Z$, which means $\mathsf{s}(a)\in Z$ and $a\mathsf{s}(a)=a$.  It then follows that $\mathsf{s}(a)^2\in ZZ\subseteq Z$ and $a\mathsf{s}(a)^2=a\mathsf{s}(a)=a$ so $\mathsf{s}(a)^2\in a^Z\subseteq\mathsf{s}(a)^\leq$, i.e. $\mathsf{s}(a)\leq\mathsf{s}(a)^2$.  But this means that we have $z\in Z$ with $\mathsf{s}(a)=\mathsf{s}(a)z=\mathsf{s}(a)^2z=\mathsf{s}(a)^2$, showing that $\mathsf{s}(a)\in\mathsf{E}(Z)$.  A dual argument shows that $\mathsf{r}(a)\in\mathsf{E}(Z)$ as well.
\end{proof}

We call $a\in S$ \emph{bisupported} if it has both a source-support and range-support.  The key thing to note here is that $Z$-invertible elements are bisupported.

\begin{prp}
If $Z$ is commutative then every $a\in S^\dagger$ is bisupported.  Specifically, for any $s\in S$ with $a<_sa$, the supports of $a$ are given by
\[\mathsf{r}(a)=as\qquad\text{and}\qquad\mathsf{s}(a)=sa.\]
\end{prp}

\begin{proof}
If $a<_sa$ then $a=asa$ and $sa\in Z$ and hence $sa\in a^Z$.  For any other $z\in a^Z$, we see that $a=az$ and hence $sa=saz=zsa=sasa$, showing that $sa\leq z$.  This shows that $sa=\mathsf{s}(a)$ and a dual argument yields $\mathsf{r}(a)=as$.
\end{proof}

Also, under suitable conditions, supports on one side yield supports on the other.

\begin{prp}
If $Z$ is binormal, $Z\subseteq\mathsf{Z}(D)$ and $a\in S^>$ then, for $a$ to be bisupported, it suffices that a has either a source-support or range-support.
\end{prp}

\begin{proof}
Say $a<_sb$.  If $Z$ is binormal then $szb\in Z$, for any $z\in{}^Z\hspace{-1pt}a$.  If $Z\subseteq\mathsf{Z}(D)$ then $aszb=zasb=za=a$, i.e. $szb\in a^Z$.  This and a dual argument shows that
\[s({}^Z\hspace{-1pt}a)b\subseteq a^Z\qquad\text{and}\qquad b(a^Z)s\subseteq{}^Z\hspace{-1pt}a.\]
If $a^Z$ has minimum $\mathsf{s}(a)$ then, in particular, $b\mathsf{s}(a)s\in{}^Z\hspace{-1pt}a$.  Moreover, \autoref{Esupport} says that $\mathsf{s}(a)\in\mathsf{E}(Z)$ so $\mathsf{s}(a)\in\mathsf{s}(a)^\dagger$.  As $sb\in a^Z=\mathsf{s}(a)^\leq$, \eqref{Inversealeqb} yields $\mathsf{s}(a)sb=\mathsf{s}(a)$ and hence $b\mathsf{s}(a)sb\mathsf{s}(a)s=b\mathsf{s}(a)\mathsf{s}(a)sbs=b\mathsf{s}(a)s\in\mathsf{E}(Z)$.  Now, for any $z\in{}^Z\hspace{-1pt}a$, we already showed that $szb\in a^Z$ so $\mathsf{s}(a)\leq szb$ and hence $\mathsf{s}(a)szb=\mathsf{s}(a)$, again by \eqref{Inversealeqb}.  Thus $b\mathsf{s}(a)sz=b\mathsf{s}(a)sbsz=b\mathsf{s}(a)szbs=b\mathsf{s}(a)s$ and hence $b\mathsf{s}(a)s\leq z$, again by \eqref{Inversealeqb}.  This shows that $b\mathsf{s}(a)s$ is the range-support of $a$ so $a$ is bisupported.  Likewise, if $a$ has range-support $\mathsf{r}(a)$ then a dual argument shows that $s\mathsf{r}(a)b$ is a source-support of $a$ hence $a$ is again bisupported.
\end{proof}

\section{Expectations}\label{Expectations}

Here we examine what more can be said for semigroup inclusions $(S,Z,D)$ where $D$ is the range of some expectation $\Phi$ on $S$ (like the canonical expectation $\phi^\rho_\mathsf{c}$ on sections of ample category bundles we saw in \autoref{ContinuousSubactions}).  Accordingly, throughout this section we make the following standing assumption.
\begin{center}
\textbf{$(S,Z,\mathrm{ran}(\Phi))$ is a semigroup inclusion, where $\Phi$ is an expectation on $S$.}
\end{center}
In other words, $S$ is a semigroup on which we have an expectation $\Phi$ whose range contains another subsemigroup $Z$ (although again if the reader is not interested in quasi-Cartan pairs or the precise conditions needed for each result below, they may just assume from the outset that $(S,Z,\Phi)$ is a well-structured semigroup).

The first thing to note is that bistability ensures that $\Phi$ is compatible with $<$.

\begin{prp}
If $Z$ is bistable then, for any $a,t,t'\in S$,
\begin{equation}\label{PhiCompatibility}
a<_{b'}b\qquad\Rightarrow\qquad\Phi(a)<_{b'}b,\quad\Phi(a)<_{b'}\Phi(b)\quad\text{and}\quad\Phi(a)<_{\Phi(b')}b.
\end{equation}
\end{prp}

\begin{proof}
Assume $a<_{b'}b$.  In particular, $bb',b'b\in Z$ and hence $b'\Phi(b),\Phi(b)b'\in Z$, as $Z$ is bistable.  Thus $\Phi(a)b'=\Phi(ab'b)b'=ab'\Phi(b)b'\in\mathrm{ran}(\Phi)Z\subseteq\mathrm{ran}(\Phi)$ and, likewise, $b'\Phi(a)\in\mathrm{ran}(\Phi)$.  To see that $\Phi(a)<_{b'}b$ now just note that
\[\Phi(a)b'b=\Phi(ab'b)=\Phi(a)=\Phi(bb'a)=bb'\Phi(a).\]
On the other hand, to see that $\Phi(a)<_{b'}\Phi(b)$ now just note that
\[\Phi(a)b'\Phi(b)=\Phi(ab'\Phi(b))=\Phi(\Phi(ab'b))=\Phi(a)=\Phi(\Phi(bb'a))=\Phi(\Phi(b)b'a)=\Phi(b)b'\Phi(a).\]

Likewise, $bb',b'b\in Z$ implies that $b\Phi(b'),\Phi(b')b\in Z$.  As $\Phi(a)b',b'\Phi(a)\in\mathrm{ran}(\Phi)$, $\Phi(a)b'=\Phi(\Phi(a)b')=\Phi(a)\Phi(b')$ and $b'\Phi(a)=\Phi(b'\Phi(a))=\Phi(b')\Phi(a)$.  As we already know that $\Phi(a)<_{b'}b$, it follows that $\Phi(a)<_{\Phi(b')}b$ too.
\end{proof}

It follows that $a<_sb$ implies $\Phi(a)<_{\Phi(s)}\Phi(b)$ and hence, for any $T\subseteq S$,
\[\Phi[T^*]\subseteq\Phi[T]^*\qquad\text{and}\qquad\Phi[T^\dagger]\subseteq\Phi[T]^\dagger.\]
In other words, $\Phi$ respects duals and $Z$-inverses whenever $\Phi$ is bistable.

In inverse $\wedge$-semigroups, the meet of any $a,b\in S$ can be expressed in terms of the expectation $\Phi_\mathsf{E}$, thanks to \cite[Theorem 1.9]{Leech1995}, specifically
\[a\wedge b=\Phi_\mathsf{E}(ab^{-1})b.\]
While domination is only reflexive on $S^\dagger$, products like those above still behave much like meets even on $S$, as the following results show.

\begin{lem}
If $(S,Z,\Phi)$ is well-structured then, for any $a,b,b',c,c'\in S$,
\begin{equation}\label{Phi(ab')b}
c<_{a'}a\ \text{ and }\ c<_{b'}b\quad\Rightarrow\quad c<_{a'aa'}\Phi(ab')b\ \text{ and }\ c<_{\Phi(a'b)b'}aa'a.
\end{equation}
\end{lem}

\begin{proof}
As $Z$ is binormal, $ab'ba'aa'\in aZa'Z\subseteq Z$ and hence $\Phi(ab')ba'aa'\in Z$, as $Z$ is bistable.  Likewise, $b'aa'aa'b\in b'ZZb\subseteq Z$ and hence
\[a'aa'\Phi(ab')b=a'\Phi(ab')aa'b=a'a\Phi(b'a)a'b=\Phi(b'a)a'aa'b\in Z,\]
as $Z\subseteq\mathsf{Z}(\mathrm{ran}(\Phi))$ is bistable and $\Phi$ is shiftable.  Also $ca'aa'=ca'\in\mathrm{ran}(\Phi)$ and
\[ca'aa'\Phi(ab')b=\Phi(ca'aa'ab')b=\Phi(cb')b=cb'b=c,\]
showing that $c<_{a'aa'}\Phi(ab')b$.

Similarly, we see that $\Phi(a'b)b'aa'a,aa'a\Phi(a'b)b'\in Z$, $\Phi(a'b)b'c=a'c\in\mathrm{ran}(\Phi)$ and $aa'a\Phi(a'b)b'c=aa'aa'c=c$, showing that $c<_{\Phi(a'b)b'}aa'a$.
\end{proof}

\begin{lem}
If $(S,Z,\Phi)$ is well-structured then, for any $a,b,b',c,c'\in S$,
\begin{equation}\label{PhiDirected}
a<_{c'}c\quad\text{and}\quad bb',b'b\in Z\qquad\Rightarrow\qquad\Phi(ab')b<_{c'}c.
\end{equation}
\end{lem}

\begin{proof}
Note $b'cc'b\in b'Zb\subseteq Z$, as $Z$ is binormal.  As $Z$ is bistable, this implies that
\begin{align*}
c'\Phi(ab')b&=c'\Phi(cc'ab')b=c'cc'\Phi(ab')b=c'\Phi(ab')cc'b=c'\Phi(ab'cc')b\\
&=\Phi(c'ab'c)c'b=c'a\Phi(b'c)c'b\in\mathrm{ran}(\Phi)Z\subseteq\mathrm{ran}(\Phi).
\end{align*}
Also $cc'\Phi(ab')b=\Phi(cc'ab')b=\Phi(ab')b$, showing that $\Phi(ab')b<_{c'}c$.
\end{proof}

When $<$ is interpolative, it satisfies the following meet-like condition, which yields an abstract bi-pseudobasis as per \cite[Definition 1.3]{BiceStarling2021}.

\begin{cor}
If $<$ is interpolative and $(S,Z,\Phi)$ is well-structured then
\begin{equation}\label{Pseudomeets}
a<s\quad\text{and}\quad b<t\qquad\Rightarrow\qquad\exists c\in s^>\cap t^>\ (a^>\cap b^>\subseteq c^>).
\end{equation}
If $S$ has a zero, it follows that $(S\setminus\{0\})^<$ is an abstract bi-pseudobasis w.r.t. $<$.
\end{cor}

\begin{proof}
Take $a,b,s,t\in S$ with $a<s$ and $b<t$.  As $<$ is interpolative, we have $u,u'\in S$ with $b<_{u'}u<t$.  Letting $c=\Phi(au')u$, we see that $a^>\cap b^>\subseteq c^>$, by \eqref{Phi(ab')b}, and $c\in s^>\cap t^>$, by \eqref{DInvariance} and \eqref{PhiDirected}.

As $<$ has interpolation, every $s\in(S\setminus\{0\})^<$ dominates some other element of $(S\setminus\{0\})^<$.  Also $<$ is transitive so $(S\setminus\{0\})^<$ is an abstract bi-pseudobasis w.r.t. $<$, according to \cite[Definition 1.3]{BiceStarling2021}.
\end{proof}

Like in \cite[Definition 3.3]{ArmstrongCastroClarkCourtneyLinMcCormickRamaggeSimsSteinberg2021}, let us call $\Phi$ \emph{quasi-Cartan} on $T\subseteq S$ if
\[\tag{Quasi-Cartan}\label{QuasiCartan}t\in T\qquad\Rightarrow\qquad\Phi(t)\leq t.\]
In this case, $\Phi(t)$ is actually the largest restriction of $t$ in $\mathrm{ran}(\Phi)$, i.e.
\[\tag{Leech}\label{Leech}\Phi(t)=\max\{r\in\mathrm{ran}(\Phi):r\leq t\}.\]
Indeed, to see this just note that $\mathrm{ran}(\Phi)\ni r\leq t$ means $r\in\mathrm{ran}(\Phi)\cap{}^Z\hspace{-1pt}rt\cap tr^Z$ and hence $r=\Phi(r)\in{}^Z\hspace{-1pt}r\Phi(t)\cap\Phi(t)r^Z$, i.e. $r\leq\Phi(t)$.

\begin{prp}\label{BistableQC}
If $Z$ is bistable then $\Phi$ is quasi-Cartan on $S^>$.
\end{prp}

\begin{proof}
Take $a\in S^>$, so we have $s,s'\in S$ with $a<_{s'}s$.  In particular, $as'\in\mathrm{ran}(\Phi)$ so $\Phi(a)=\Phi(as's)=as'\Phi(s)$.  Also $s's\in Z$ so $s'\Phi(s)\in Z$, as $Z$ is bistable, and hence $\Phi(a)=\Phi(\Phi(a))=\Phi(as'\Phi(s))=\Phi(a)s'\Phi(s)$.  Dually, $\Phi(s)s'\in Z$ and $\Phi(a)=\Phi(s)s'\Phi(a)=\Phi(s)s'a$, showing that $\Phi(a)\leq a$.
\end{proof}

On the other hand note that, for all $s,t\in S$,
\begin{align}
\tag{Bistability$'$}\label{Bistability'}st\in Z\text{ and }\Phi(s)\leq s\qquad&\Rightarrow\qquad\Phi(s)t\in Z.\\
\tag{Bistability$''$}\label{Bistability''}st\in\mathrm{ran}(\Phi)\text{ and }\Phi(s)\leq s\qquad&\Rightarrow\qquad\Phi(s)t\in\mathrm{ran}(\Phi).
\end{align}
Indeed, if $st\in Z$ and $\Phi(s)\in Zs$ then $\Phi(s)t\in Zst\subseteq ZZ\subseteq Z$, while if $st\in\mathrm{ran}(\Phi)$ and $\Phi(s)\in\mathrm{ran}(\Phi)t$ then $\Phi(s)t\in\mathrm{ran}(\Phi)st\subseteq\mathrm{ran}(\Phi)\mathrm{ran}(\Phi)\subseteq\mathrm{ran}(\Phi)$.  In particular,
\[\Phi\text{ is quasi-Cartan (on all of $S$)}\qquad\Rightarrow\qquad Z\text{ and }\mathrm{ran}(\Phi)\text{ are bistable}.\]

Similarly, we see that
\[\tag{Unistability}Z^\mathsf{C}\ni\Phi(s)\leq s\qquad\Rightarrow\qquad\Phi(s)s=\Phi(s)\Phi(s)=s\Phi(s).\]
Indeed if $Z^\mathsf{C}\ni\Phi(s)\leq s$ then we have $z\in Z$ with $zs=z\Phi(s)=\Phi(s)=\Phi(s)z=sz$, by \eqref{ZCleq}, and hence $\Phi(s)s=\Phi(s)zs=\Phi(s)\Phi(s)=sz\Phi(s)=s\Phi(s)$.  In particular, if $\Phi$ is quasi-Cartan (on all of $S$) and $Z\subseteq\mathsf{Z}(\mathrm{ran}(\Phi))$ then, for all $s\in S$,
\[\Phi(s)s=s\Phi(s)\in\mathrm{ran}(\Phi).\]

We also have a converse of \eqref{PhiCompatibility}, at least when $Z$ is commutative and $\Phi$ is quasi-Cartan.  Specifically, if $Z$ is commutative then, for all $a,s\in S$,
\[a<\Phi(s)\leq s\qquad\Rightarrow\qquad a=\Phi(a)<s.\]
Indeed, if $a<\Phi(s)\leq s$ then $a<s$ is immediate from \eqref{RightAuxiliarity}, while $a<_{r'}r\in\mathrm{ran}(\Phi)$ implies that $a=ar'r\in\mathrm{ran}(\Phi)\mathrm{ran}(\Phi)\subseteq\mathrm{ran}(\Phi)$.

In structured semigroups, we can further show that $\Phi$ is shiftable.

\begin{prp}
If $(S,Z,\mathrm{ran}(\Phi))$ is a structured semigroup and $S=S^>$ then
\begin{equation}\label{BistableQCShiftable}
Z\text{ is bistable}\qquad\Leftrightarrow\qquad\Phi\text{ is quasi-Cartan}\qquad\Rightarrow\qquad\Phi\text{ is shiftable}.
\end{equation}
\end{prp}

\begin{proof}
If $S=S^>$ then \autoref{BistableQC} and the comments after imply that $Z$ is bistable precisely when $\Phi$ is quasi-Cartan.  To see that $\Phi$ is also shiftable when $(S,Z,\mathrm{ran}(\Phi))$ is also structured semigroup, take any $a,s\in S$.  As $\Phi(as)\leq as$, we have $z\in Z$ with $\Phi(as)=\Phi(as)z=asz$.  As $s\in S=S^>$, we have $t,t'\in S$ with $s<_{t'}t$.  Then $s\Phi(as)t'=saszt'=sast'tzt'$ and
\[s\Phi(as)t'=s\Phi(as)zt'=s\Phi(ast't)zt'=s\Phi(as)t'tzt'.\]
As $Z$ is binormal, $tzt'\in Z$ so \eqref{OneSidedRestriction} then yields $s\Phi(as)t'\leq sast'$.  As $Z\subseteq\mathsf{Z}(\mathrm{ran}(\Phi))$, $\Phi(as)=\Phi(ast't)=\Phi(as)t't=t't\Phi(as)$ and hence $t\Phi(as)t'\in\mathrm{ran}(\Phi)$, as $\mathrm{ran}(\Phi)$ is trinormal.  As $st'\in\mathrm{ran}(\Phi)$ too, $s\Phi(as)t'=st't\Phi(as)t'\in\mathrm{ran}(\Phi)$.  As $\Phi$ is quasi-Cartan and hence satisfies \eqref{Leech}, $s\Phi(as)t'\leq\Phi(sast')=\Phi(sa)st'$ and hence $s\Phi(as)=s\Phi(as)t't\leq\Phi(sa)st't=\Phi(sa)s$, by \eqref{Invariance}.  A dual argument yields $\Phi(sa)s\leq s\Phi(as)$ and hence $\Phi(sa)s=s\Phi(as)$, as $\leq$ is antisymmetric, by \autoref{RestrictionTransitiveAntisymmetric}.  This shows that $\Phi$ is indeed shiftable.
\end{proof}

It follows from \autoref{WellStructuredSemigroups} and \eqref{BistableQCShiftable} above that, when $S=S^>$, we can replace the shiftability condition in the definition of well-structured semigroups in \autoref{WSS} with requirement that $\mathrm{ran}(\Phi)$ is trinormal.  In this case, well-structured semigroups $(S,Z,\Phi)$ could be viewed instead as semigroup inclusions $(S,Z,D)$, where $D$ is the range of some shiftable expectation making $Z$ bistable, which is then unique by \autoref{BistableQC} and \eqref{Leech}.  This viewpoint would be more in line with the way Cartan and quasi-Cartan pairs are usually considered in C*-algebra and Steinberg algebra theory.  However, as the expectation plays such an important role, we feel it deserves to considered explicitly as an intrinsic part of any well-structured semigroup.  The expectation will be particularly important for the morphisms we will consider in \autoref{LawsonSteinbergDuality}, which must preserve the expectation itself, not just its range.

\section{Filters}\label{Filters}

First let us make the standing assumption that
\begin{center}
\textbf{$(S,Z,\Phi)$ is a well-structured semigroup.}
\end{center}
As usual, we call $F\subseteq S$ a \emph{filter} if $F$ is a down-directed up-set, i.e.
\[a,b\in F\qquad\Leftrightarrow\qquad\exists f\in F\ (f<a,b).\]
An \emph{ultrafilter} is a maximal proper filter.  These will from the base groupoid of the bundle on which we will represent the semigroup.  Again the motivating situation to keep in mind is when $(S,Z,\Phi)$ arises from an ample category bundle as in \autoref{ContinuousSubactions} -- then each ultrafilter will consist of all sections whose value at some fixed point in the base groupoid is invertible.

At first it will be instructive to consider more general filters and even more general cosets.  Specifically, as in \cite[Definition 6.1]{Bice2022}, we call $C\subseteq S$ a \emph{coset} if
\[\tag{Coset}\label{Coset}CC^*C\subseteq C=C^<,\]
where $C^*=\{s\in S:\exists c,d\in C\ (c<_sd)\}$ is the dual of $C$ defined earlier.  By \cite[Proposition 7.1]{Bice2022}, the non-empty cosets $\mathcal{C}(S)$ form a groupoid where
\[\mathsf{s}(B)=(B^*B)^<,\quad\mathsf{r}(B)=(BB^*)^<,\quad B\cdot C=(BC)^<,\quad\text{and}\quad C^{-1}=C^*,\]
the product being defined of course only when $\mathsf{s}(B)=\mathsf{r}(C)$.

By an \emph{ideal} of a groupoid $G$, we mean a subset $I\subseteq G$ such that $IG\cup GI\subseteq I$.

\begin{prp}\label{FilterIdeal}
The non-empty filters $\mathcal{F}(S)$ form an ideal of $\mathcal{C}(S)$.
\end{prp}

\begin{proof}
By \cite[Proposition 11.5]{Bice2022}, it suffices to show that the filters are precisely the down-directed cosets.  To see this note first that if $F\subseteq S$ is a coset then $F^{**}=F^{<<}\subseteq F$, by \cite[Proposition 5.4]{Bice2022}.  If $F$ is also down-directed, it follows that $F$ is a filter, by \cite[Proposition 11.3]{Bice2022}.  Conversely, if $F\subseteq S$ is a filter then $F$ is a coset, again by \cite[Proposition 11.3]{Bice2022}, as $\mathrm{ran}(\Phi)$ is diagonal, by \autoref{WellStructuredSemigroups}.
\end{proof}

By \cite[Proposition 7.2]{Bice2022}, the unit cosets $\mathcal{C}^0=\mathcal{C}(S)^0$ are precisely those containing an element of $Z$ or, equivalently, $\mathrm{ran}(\Phi)$.  On the other hand, the unit filters $\mathcal{F}^0=\mathcal{F}(S)^0$ can be characterised as the non-empty $\Phi$-invariant cosets.

\begin{prp}\label{UnitFilters}
A non-empty coset $C$ is a unit filter iff $\Phi[C]\subseteq C$, i.e.
\[\mathcal{F}^0=\{C\in\mathcal{C}(S):\Phi[C]\subseteq C\}.\]
\end{prp}

\begin{proof}
By \cite[Proposition 11.6]{Bice2022}, the unit filters are precisely the non-empty cosets $C$ that are generated by their diagonal, meaning that $C=(C\cap\mathrm{ran}(\Phi))^<$.

So if $C$ is a unit filter and $c\in C$ then we have $b\in C\cap\mathrm{ran}(\Phi)$ with $b<c$ and hence $b=\Phi(b)<\Phi(c)$, by \eqref{PhiCompatibility}.  Thus $\Phi(c)\in C^<\subseteq C$, showing that $\Phi[C]\subseteq C$.

Conversely, say $C$ is a coset with $\Phi[C]\subseteq C$.  For any $c\in C$, we have $b\in C$ with $b<c$ and hence $\Phi(b)<c$, again by \eqref{PhiCompatibility}.  As $\Phi(b)\in\Phi[C]\subseteq C\cap\mathrm{ran}(\Phi)$, this shows that $C\subseteq(C\cap\mathrm{ran}(\Phi))^<$ so $C$ is a unit filter.
\end{proof}

General filters can be characterised like in \eqref{Coset} but again with an extra $\Phi$.

\begin{prp}\label{FilterPhi}
A subset $F\subseteq S$ is a filter precisely when $\Phi[FF^*]F\subseteq F=F^<$.
\end{prp}

\begin{proof}
If $F\subseteq S$ is a filter then $\mathsf{r}(F)\supseteq FF^*$ is a unit filter (as long as $F\neq\emptyset$) so $\Phi[\mathsf{r}[F]]\subseteq\mathsf{r}[F]$, by the above result, and hence $\Phi[FF^*]F\subseteq\mathsf{r}[F]F\subseteq F$.

Conversely, say $\Phi[FF^*]F\subseteq F=F^<$.  For $F$ to be a filter it suffices that $F$ is directed, as we are already assuming $F^<\subseteq F$.  Accordingly, take $t,u\in F$.  As $F\subseteq F^<$, we have $a,b,c\in F$ and $c'\in F^*$ with $a<t$ and $b<_{c'}c<u$.  By \eqref{PhiDirected}, $\Phi(ac')c<t$ and, by \cite[Proposition 4.8]{Bice2022}, $\Phi(ac')c<u$.  As
\[\Phi(ac')c\in\Phi[FF^*]F\subseteq F,\]
this shows that $F$ is directed and hence a filter.
\end{proof}

The coset groupoid $\mathcal{C}(S)$ also has a natural topology generated by the slices
\[\mathcal{C}_a=\{C\in\mathcal{C}(S):a\in C\}.\]
In other words, $(\mathcal{C}_a)_{a\in S}$ is a subbasis for this topology.  By \cite[Theorem 7.4]{Bice2022}, this topology makes $\mathcal{C}(S)$ an \'etale groupoid.  As $\mathcal{F}(S)$ is ideal of $\mathcal{C}(S)$, by \autoref{FilterIdeal}, $\mathcal{F}(S)$ is also an \'etale groupoid in the subspace topology (see \cite[Proposition 2.7]{Bice2021}), i.e. the topology with subbasis $(\mathcal{F}_a)_{a\in S}$, where
\[\mathcal{F}_a=\mathcal{C}_a\cap\mathcal{F}(S).\]
As expected, applying $\Phi$ here corresponds to restricting to the diagonal $\mathcal{F}^0$.

\begin{prp}
For any $a\in S^>$,
\[\mathcal{F}_{\Phi(a)}=\mathcal{F}_a\cap\mathcal{F}^0.\]
\end{prp}

\begin{proof}
If $F\in\mathcal{F}_a\cap\mathcal{F}^0$ then we have $e\in\mathrm{ran}(\Phi)\cap F$.  As $F$ is a filter, we have $f\in F$ with $f<a,e$ and hence $f=\Phi(f)<\Phi(a)$, by \eqref{PhiCompatibility}.  Thus $F\in\mathcal{F}_f\subseteq\mathcal{F}_{\Phi(a)}$, showing that $\mathcal{F}_a\cap\mathcal{F}^0\subseteq\mathcal{F}_{\Phi(a)}$.

Conversely, if $F\in\mathcal{F}_{\Phi(a)}$ then certainly $F\in\mathcal{F}^0$, as $\Phi(a)\in\mathrm{ran}(\Phi)$.  Moreover, as $F\subseteq F^<$, we can take $f\in F$ with $f<\Phi(a)$ and hence $f<a$, by \autoref{BistableQC} and \eqref{RightAuxiliarity}, as $a\in S^>$.  This shows that $\mathcal{F}_{\Phi(a)}\subseteq\mathcal{F}_a\cap\mathcal{F}^0$, as required.
\end{proof}

The following observation will be needed later in \autoref{ContinuousSection}.

\begin{lem}
For any $a,b,b'\in S$,
\begin{equation}\label{rFt}
a<_{b'}b\qquad\Rightarrow\qquad\mathsf{r}[\mathcal{F}_a]=\mathcal{F}_{ab'}.
\end{equation}
\end{lem}

\begin{proof}
If $a<_{b'}b$ then $b'\in F^*$ whenever $a\in F$.  It follows that $\mathcal{F}_a^{-1}\subseteq\mathcal{F}_{b'}$ and hence $\mathsf{r}[\mathcal{F}_a]\subseteq\mathcal{F}_a\cdot\mathcal{F}_{b'}\subseteq\mathcal{F}_{ab'}$.  Conversely, if $F\in\mathcal{F}_{ab'}$ then $ab'\in F$, $ab'bb'=ab'$ and $bb',b'b\in Z$.  Letting $T=(Fb)^<$, \cite[Proposition 6.6]{Bice2022} yields $\mathsf{r}(T)=\mathsf{r}(F)=F\in\mathcal{F}^0$, as $ab'\in\mathrm{ran}(\Phi)$, and hence $T\in\mathcal{F}_a$, by \autoref{FilterIdeal} and the fact that $a=ab'b\in Fb\subseteq T$.  This shows that $\mathcal{F}_{ab'}\subseteq\mathsf{r}[\mathcal{F}_a]$.
\end{proof}

When $S$ has a zero, the proper filters are precisely those avoiding $0$.

\begin{prp}
If $S$ has a zero and $0\in F\in\mathcal{F}(S)$ then $0\in Z$ and $F=0^<=S$.
\end{prp}

\begin{proof}
If $0\in F$ then we have $f\in F$ and $s\in S$ with $f<_s0$.  In particular, $0=0s\in Z$ and hence $0<_0s$, for all $s\in S$.  Thus $S=0^<\subseteq F^<=F$.
\end{proof}

As mentioned above, we will eventually focus our attention on the ultrafilters
\[\mathcal{U}(S)=\{U\subseteq S:U\text{ is a maximal proper filter}\}.\]
When $S$ has a zero, these are precisely the maximal filters avoiding zero, by the above result.  In this case, $\mathcal{U}(S)$ is an ideal, by \cite[Theorem 11.7]{Bice2022}, and hence an \'etale subgroupoid of $\mathcal{F}(S)$.  We also have the following observation.

\begin{prp}\label{UltrafilterDichotomy}
If $S$ has a zero then, for any $U\in\mathcal{U}(S)$, we have $u\in U$ with
\[\Phi(u)=0\qquad\text{or}\qquad\Phi(u)=u.\]
\end{prp}

\begin{proof}
Take any $U\in\mathcal{U}(S)$.  In particular, $U$ is directed and hence so is $\Phi[U]$, by \eqref{PhiCompatibility}.  Thus $\Phi[U]^<$ is a filter which contains $U$, again by \eqref{PhiCompatibility}.  If $0\notin\Phi[U]$ then $0\notin\Phi[U]^<$ and hence $U=\Phi[U]^<$, by maximality.  In other words, either $0\in\Phi[U]$ or $\Phi[U]\subseteq U$.  In the latter case, we can take any $u\in U$ and just note that $\Phi(\Phi(u))=\Phi(u)\in U$, as required.
\end{proof}

Let $\mathcal{U}^0=\mathcal{U}(S)^0$ denote the unit ultrafilters.

\begin{cor}\label{HausdorffUltrafilters}
If $S$ has a zero then $\mathcal{U}^0$ is closed and hence $\mathcal{U}(S)$ is Hausdorff.
\end{cor}

\begin{proof}
If $U\in\mathcal{U}(S)\setminus\mathcal{U}^0$ then we must have $u\in U$ with $\Phi(u)=0$, by \autoref{UltrafilterDichotomy}.  It follows that $\mathcal{U}_u\cap\mathcal{U}^0=\emptyset$, showing that $\mathcal{U}(S)\setminus\mathcal{U}^0$ is open and hence $\mathcal{U}^0$ is closed.  As $\mathcal{U}^0$ is Hausdorff (see the proof of \cite[Theorem 11.7]{Bice2022}), it follows that the entirety of $\mathcal{U}(S)$ is also Hausdorff (see \cite[Lemma 2.3.2]{Sims2018}).
\end{proof}

Before moving on, let us note that the $Z$-invertibles in a given ultrafilter can be characterised as follows, which will be useful later on.

\begin{prp}
If $S$ has a zero then, for all $U\in\mathcal{U}(S)$ and $a,s\in S$ with $a\in s^\dagger$,
\begin{equation}\label{UltrafilterReflexives}
a\in U\qquad\Leftrightarrow\qquad0\notin\Phi[Us].
\end{equation}
\end{prp}

\begin{proof}
If we had $a\in U$ and $0\in\Phi[Us]$ then \autoref{FilterPhi} would yield
\[0\in\Phi[Us]a\subseteq\Phi[UU^*]U\subseteq U,\]
a contradiction.  This proves $\Rightarrow$.

Conversely, note $Us$ is directed, by \eqref{Multiplicativity}, as $U$ is directed and $s<s$.  Thus $\Phi[Us]$ is directed, by \eqref{PhiCompatibility}, and hence so is $\Phi[Us]a$, again because $a<a$.  If $0\notin\Phi[Us]=\Phi[Us]as$ then $0\notin\Phi[Us]a$ so $(\Phi[Us]a)^<$ is a proper filter containing $U$, by \eqref{PhiDirected}.  Thus $a\in(\Phi[Us]a)^<=U$, by maximality, proving $\Leftarrow$.
\end{proof}

\section{Equivalence}\label{Equivalence}

So far we have the base of the bundle we will represent our well-structured semigroups and semimodules on, namely the filter groupoid (in fact, later we will even restrict further to the ultrafilters, but we need not worry about that yet).  To construct the rest of the bundle we need to consider certain equivalence relations defined by these filters, which we now proceed to examine.

Here we consider semimodules and thus make the following standing assumption.
\begin{center}
\textbf{$(A,S,Z,\Phi)$ is a well-structured semimodule.}
\end{center}

For any $s\in S$ and $T\subseteq S$, we define relations $\sim_s$ and $\sim_T$ on $A$ by
\begin{align*}
a\sim_sb\qquad&\Leftrightarrow\qquad\Phi(as)=\Phi(bs)\quad\text{and}\quad\Phi(sa)=\Phi(sb).\\
a\sim_Tb\qquad&\Leftrightarrow\qquad\exists s\in T^*\ (a\sim_sb).
\end{align*}
Yet again the motivating situation to keep in mind is when $(S,Z,\Phi)$ arises from an ample category bundle as in \autoref{ContinuousSubactions} and $T$ is some ultrafilter, necessarily corresponding to some point in the base groupoid -- then $a\sim_Tb$ is just saying that $a$ and $b$ take the same value at this point.

Note $s$ above has to be in $T^*$, not $T$.  Actually, for $a\sim_Tb$ to hold, it suffices that either $\Phi(as)=\Phi(bs)$ or $\Phi(sa)=\Phi(sb)$, for some $s\in T^*$.

\begin{prp}
For all $T\subseteq S$,
\[a\sim_Tb\quad\Leftrightarrow\quad\exists s\in T^*\ (\Phi(as)=\Phi(bs))\quad\Leftrightarrow\quad\exists s\in T^*\ (\Phi(sa)=\Phi(sb)).\]
\end{prp}

\begin{proof}
First we claim that, for all $s,t',t\in S$,
\begin{equation}\label{OneSidedSim}
s<_{t'}t\quad\text{and}\quad\Phi(at')=\Phi(bt')\qquad\Rightarrow\qquad a\sim_{t'tt'}b
\end{equation}
Indeed, if the left side holds then $\Phi(at'tt')=\Phi(at')tt'=\Phi(bt')tt'=\Phi(bt'tt')$ and
\[\Phi(t'tt'a)=t't\Phi(t'a)=\Phi(t'a)t't=t'\Phi(at')t=t'\Phi(bt')t=\Phi(t'tt'b).\]

Now if we have $t'\in T^*$ with $\Phi(at')=\Phi(bt')$ then we have $s\in T$ with $s<_{t'}t$ and hence $a\sim_{t'tt'}b$, by \eqref{OneSidedSim}.  Also $s<_{t'tt'}t$ so $t'tt'\in T^*$ and hence $a\sim_Tb$.  This proves the first equivalence and the second follows dually.
\end{proof}

We will use these slightly simpler equivalents of $\sim_T$ without further reference.

Another thing to note is that we can replace $T$ with $T^<$.

\begin{prp}
$\sim_T\ =\ \sim_{T^<}$ when $T\subseteq T^<$, i.e. for all $a,b\in A$,
\begin{equation}\label{Tup}
a\sim_Tb\qquad\Leftrightarrow\qquad a\sim_{T^<}b.
\end{equation}
\end{prp}

\begin{proof}
For any $T\subseteq S$, \cite[Proposition 4.5]{Bice2022} implies that $T^{<*}\subseteq T^*$ and hence $\sim_{T^<}\ \subseteq\ \sim_T$.  If $T\subseteq T^<$ then $T^*\subseteq T^{<*}$ too and hence $\sim_T\ =\ \sim_{T^<}$.
\end{proof}

For the $T$ we really care about, $\sim_T$ will be an equivalence relation.

\begin{prp}\label{rhoUpseudo}
If $T\neq\emptyset$ is directed then $\sim_T$ is an equivalence relation on $A$.
\end{prp}

\begin{proof}
First note that
\begin{equation}\label{simBelow}
s<t\quad\text{and}\quad a\sim_tb\qquad\Rightarrow\qquad a\sim_sb.
\end{equation}
Indeed, if $s<_{t'}t$ and $\Phi(at)=\Phi(bt)$ then
\[\Phi(as)=\Phi(att's)=\Phi(at)t's=\Phi(bt)t's=\Phi(btt's)=\Phi(bs).\]
Likewise, $s<t$ and $\Phi(ta)=\Phi(tb)$ implies $\Phi(sa)=\Phi(sb)$, proving the claim.

We immediately see that $\sim_s$ is an equivalence relation, for any $s\in S$.  It follows immediately that $\sim_T$ is also symmetric, for any $T\subseteq S$, and reflexive, as long as $T\neq\emptyset$.    To see that $\sim_T$ is transitive when $T$ is also directed, just note that then $T^*$ is directed too, by \cite[Proposition 11.2]{Bice2022}.  Thus if $s,t\in T^*$ and $a\sim_sb\sim_tc$ then we have $u\in T^*$ with $u<s,t$ so \eqref{simBelow} yields $a\sim_ub\sim_uc$ and hence $a\sim_uc$.
\end{proof}

Multiplying by elements of $T^Z$ does not change the equivalence class.

\begin{prp}
For any $a\in A$ and $T\subseteq S$ with $T\subseteq T^<$,
\begin{equation}\label{az0}
z\in T^Z\qquad\Rightarrow\qquad a\sim_Taz.
\end{equation}
\end{prp}

\begin{proof}
If $z\in T^Z$, we have $t\in T$ with $tz=t$.  For any $s,t'\in S$ with $T\ni s<_{t'}t$,
\[\Phi(t'tt'az)=t't\Phi(t'a)z=t'tz\Phi(t'a)=t't\Phi(t'a)=\Phi(t'tt'a),\]
showing that $a\sim_Taz$.
\end{proof}

When $T$ is a unit filter, we can characterise $\sim_T$ using $T^Z$ as follows.

\begin{prp}
For any $a,b\in A$ and unit filter $T\in\mathcal{F}^0$,
\begin{equation}\label{rhoUnit}
a\sim_Tb\qquad\Leftrightarrow\qquad\exists z\in T^Z\ (\Phi(a)z=\Phi(b)z).
\end{equation}
\end{prp}

\begin{proof}
If $T\in\mathcal{F}^0$ then $T^Z\subseteq T$ -- see \cite[Proposition 8.3]{Bice2022}.  Thus if we have $z\in T^Z$ with $\Phi(a)z=\Phi(b)z$ then $a\in T=T^*$ and $\Phi(az)=\Phi(bz)$, showing that $a\sim_Tb$.  Conversely, if $a\sim_Tb$ then we have $t\in T^*=T$ with $\Phi(at)=\Phi(bt)$.  As $T\in\mathcal{F}^0$, we have $e\in T\cap\mathrm{ran}(\Phi)$.  Taking $s\in T$ with $s<t,e$, we see that $s\in\mathrm{ran}(\Phi)$ and $\Phi(a)s=\Phi(as)=\Phi(bs)=\Phi(b)s$ (see \eqref{simBelow}).  Taking $r,s'\in T$ with $r<_ss'$, we see that $ss'\in T^Z$ and $\Phi(a)ss'=\Phi(b)ss'$, as required.
\end{proof}

Note this implies that, for all $a\in A$,
\begin{equation}\label{rhoUnit2}
T\in\mathcal{F}^0\qquad\Rightarrow\qquad a\sim_T\Phi(a).
\end{equation}
Indeed, taking any $z\in T^Z$ we immediately see that $\Phi(a)z=\Phi(\Phi(a))z$.

Later we need the following to show that equivalence classes can be multiplied.

\begin{prp}
For any $T,U\subseteq S$, $a,b\in A$ and $s\in T^{<>}$.
\begin{equation}\label{saTUsb}
a\sim_Ub\qquad\Rightarrow\qquad sa\sim_{TU}sb.
\end{equation}
\end{prp}

\begin{proof}
If $s\in T^{<>}$ then we have $r,s',t,t'\in S$ with $s<_{s'}r$ and $T\ni t<_{t'}r$ and hence $t<_{t'rt'}r$.  If we also have $u\in U^*$ with $a\sim_ub$ then $ut'rt'\in U^*T^*\subseteq(TU)^*$ (see \cite[Proposition 4.7]{Bice2022}).  Also $t's=t'rs's\in Z\mathrm{ran}(\Phi)\subseteq\mathrm{ran}(\Phi)$ so
\[\Phi(saut'rt')=rt'\Phi(saut')=r\Phi(t'sau)t'=rt's\Phi(au)t'=rt's\Phi(bu)t'=\Phi(sbut'rt'),\]
showing that $sa\sim_{TU}sb$.
\end{proof}

On $T^>$ we can verify $\sim_T$ without the expectation.  It follows that, on $T\cap S^>$, our $\sim_T$ is the same relation as in \cite[\S8]{Bice2022} (see \cite[Proposition 8.4]{Bice2022}).

\begin{prp}
If $T\subseteq S$ is directed and $a,b\in T^>$ then
\begin{equation}\label{rhoazbz}
a\sim_Tb\qquad\Leftrightarrow\qquad\exists t'\in T^*\ (at'=bt').
\end{equation}
\end{prp}

\begin{proof}
Take $r,r',s,s'\in T$ with $a<_{r'}r$ and $b<_{s'}s$.  Assume $a\sim_{u'}b$ and $u'\in T$.  As $T$ is directed, so is $T^*$ and hence we have $t'\in T^*$ with $t'<r',s',u'$.  So we have $r'',s'',u''\in S$ with $t'<_{r''}r'$, $t'<_{s''}s'$ and $t'<_{u''}u'$, from which we see that $at'=ar'r''t'\in\mathrm{ran}(\Phi)$ and $bt'=bs's''t'\in\mathrm{ran}(\Phi)$.  Thus
\[at'=\Phi(at')=\Phi(au'u''t')=\Phi(au')u''t'=\Phi(bu')u''t'=bt'.\qedhere\]
\end{proof}

Moreover, every equivalence class in $A$ will contain elements of $T^>$.

\begin{prp}\label{as<U0}
If $t\in T\subseteq T^<$ and $a\in A$ then we have $s\in S$ with
\[a\sim_Ts<t.\]
\end{prp}

\begin{proof}
As $T\subseteq T^<$, we have $t',u,u',v,v',w\in S$ with $T\ni w<_{v'}v<_{u'}u<_{t'}t$.  Letting $s=\Phi(au')u$, we see that $s<t$, by \cite[Proposition 4.8]{Bice2022}.  Replacing $v'$ with $v'vv'$ if necessary, we may assume that $v'<_uu'$ (see \cite[Proposition 5.9]{Bice2022}).  Then $\Phi(sv')=\Phi(\Phi(au')uv')=\Phi(\Phi(au'uv'))=\Phi(av')$, showing that $a\sim_Ts$.
\end{proof}

For any $a\in A$ and $T\subseteq S$ let
\[a_T=a^{\sim_T}=\{b\in A:a\sim_Tb\}.\]
So when $T\in\mathcal{F}(S)$, this is just the equivalence class containing $a$.  The result above in \autoref{as<U0} can now be rephrased as follows.

\begin{cor}\label{EquivalenceContainment}
If $T\subseteq T^<$ and $a\in A$ then $T\subseteq a_T^<$.
\end{cor}

When $<$ is interpolative, $a_T^<$ will be a filter whenever $T$ is.

\begin{prp}\label{EquivalenceFilter}
If $<$ is interpolative and $T\in\mathcal{F}(S)$ then $a_T^<\in\mathcal{F}(S)$.
\end{prp}

\begin{proof}
Take any $b,c\in a_T^<$, so we have $r,s\in a_T$ with $r<b$ and $s<c$.  By \autoref{rhoUpseudo}, $r\sim_Ts$ so we have $t'\in T^*$ with $\Phi(rt')=\Phi(st')$.  Taking $t,u\in S$ with $T\ni u<_{t'}t$, we may let $q=\Phi(rt')t=\Phi(st')t$.  As in the previous proof, we see that $q\sim_Tr,s$ and hence $q\in a_T$.  By \eqref{PhiDirected}, $q<b,c$.  As $<$ is interpolative, we have $d,e\in S$ with $q<d<b$ and $q<e<c$.  By \eqref{Pseudomeets}, we then have $f<b,c$ with $q\in d^>\cap e^>\subseteq f^>$, i.e. $q<f$.  This shows that $a_T^<$ is directed and hence a filter.
\end{proof}

Consequently, ultrafilters can be recovered from their non-zero equivalence classes.

\begin{cor}
If $<$ is interpolative and $S$ has a zero then, for any $U\in\mathcal{U}(S)$,
\begin{equation}\label{aUup}
a\not\sim_U0\qquad\Rightarrow\qquad U=a_U^<.
\end{equation}
\end{cor}

\begin{proof}
Assume $<$ is interpolative, $a\not\sim_U0$ and $U\in\mathcal{U}(S)$.  By \autoref{EquivalenceContainment} and \autoref{EquivalenceFilter}, $U\subseteq a_U^<\in\mathcal{F}(U)$.  Also $0\notin a_U$, as $a\not\sim_U0$, so $0\notin a_U^<$ and hence $U=a_U^<$, by maximality.
\end{proof}

\section{Filter Bundles}\label{FilterBundles}

Let us first reiterate our standing assumption that
\begin{center}
\textbf{$(A,S,Z,\Phi)$ is a well-structured semimodule.}
\end{center}
In the present section we will use the equivalence classes considered above to construct a bundle on which to represent $A$.  First, for any $a\in A$ and $T\subseteq S$, let
\[[a,T]=(a_T,T^<).\]
Actually, we will be exclusively interested in the case when $T\neq\emptyset$ is directed.  In this case $T^<\in\mathcal{F}(S)$ and $a_T$ is the equivalence class of $a$ with respect to $\sim_T$, which is the same as $\sim_{T^<}$ (see \eqref{Tup}).  Let
\[\mathcal{F}[A]=\{[a,T]:a\in A\text{ and }T\in\mathcal{F}(S)\}\]

\begin{thm}\label{FilterCategory}
$\mathcal{F}[A]$ is a category under the product
\[[a,T][b,U]=[ab,TU]\]
defined when $\mathsf{s}(U)=\mathsf{r}(V)$, $a\in U^>$ and $b\in V^>$.  In fact, the product above is valid as long as $\mathsf{s}(U)=\mathsf{r}(V)$ and either $a\in T^>$ or $b\in U^>$.  Moreover,
\begin{align}
\label{FS}\mathcal{F}[A]&=\{[a,T]:a\in T^>\text{ and }T\in\mathcal{F}(S)\}.\\
\label{FStimes}\mathcal{F}[A]^\times&=\{[a,T]:a\in T\text{ and }T\in\mathcal{F}(S)\}.\\
\label{FS0}\mathcal{F}[A]^0&=\{[a,T]:a\in T^Z\text{ and }T\in\mathcal{F}^0\}.
\end{align}
\end{thm}

\begin{proof}
First we must show that the product is well defined.  For this, note that if $T,U\in\mathcal{F}(S)$, $a,a'\in T^>$, $b,b'\in U^>$, $a\sim_Ta'$ and $b\sim_Ub'$ then \eqref{saTUsb} yields
\[ab\sim_{TU}ab'\sim_{TU}a'b'.\]
and hence $[ab,TU]=[a'b',TU]$, as required.

Next we must show that the product is associative.  For this, take $T,U,V\in\mathcal{F}(S)$, $a\in T^>$, $b\in U^>$ and $c\in V^>$.  Then $ab\in T^>U^>\subseteq(TU)^>=(T^<U^<)^>\subseteq(TU)^{<>}$ and, likewise, $bc\in(UV)^{<>}$.  If $\mathsf{s}(T)=\mathsf{r}(U)$ and $\mathsf{s}(U)=\mathsf{r}(V)$ then it follows that
\[([a,T][b,U])[c,V]=[ab,TU][c,V]=[abc,TUV]=[a,T][bc,UV]=[a,T]([b,U][c,V]).\]

Next note that \eqref{FS} follows immediately from \autoref{as<U0}, i.e. for any $T\in\mathcal{F}(S)$ and $a\in A$, we can always find $s\in T^>$ with $[a,T]=[s,T]$.  It follows that $[a,T][b,U]$ is defined whenever $\mathsf{s}(U)=\mathsf{r}(T)$ so $\mathcal{F}[A]$ is a semicategory.  To actually compute $[a,T][b,U]$, we may of course need to replace $a$ and $b$ with elements of $T^>$ and $U^>$ respectively.  However, only one replacement is actually necessary, i.e. $[a,T][b,U]=[ab,TU]$ as long as either $a\in T^>$ or $b\in U^>$, again thanks to \eqref{saTUsb}.

For \eqref{FS0}, take $U\in\mathcal{F}^0$ and $z\in U^Z=U^{*Z}={}^ZU$ (see \cite[Proposition 8.3]{Bice2022}).  To see that $[z,U]$ is a unit in $\mathcal{F}[A]$, take $T\in\mathcal{F}(S)$ with $\mathsf{s}(T)=U$ and $a\in T^>$.  Note that $z\in U^Z=\mathsf{s}(T)^Z=T^Z$, by \cite[Proposition 8.3]{Bice2022}.  Then \eqref{az0} yields $a\sim_Taz$ so $[a,T]=[az,T]=[a,T][z,U]$.  Likewise, if we had $\mathsf{r}(T)=U$ instead then $z\in{}^ZU={}^Z\mathsf{r}(T)={}^ZT$ and $[a,T]=[za,T]=[z,U][a,T]$, i.e. $[z,U]\in\mathcal{F}[A]^0$.  It follows that every $[a,T]\in\mathcal{F}[A]$ has a source and range unit, specifically
\begin{align}
\label{source}\mathsf{s}([a,T])&=[z,\mathsf{s}(T)],\text{ for any }z\in T^Z.\\
\label{range}\mathsf{r}([a,T])&=[z,\mathsf{r}(T)],\text{ for any }z\in{}^ZT.
\end{align}
Thus we have shown that $\mathcal{F}[A]$ is a category satisfying \eqref{FS0}.

For \eqref{FStimes}, take any $T\in\mathcal{F}(S)$ and $t\in T$.  Then we have $s,t,t',u'\in S$ with $T\ni v<_{u'}u<_{t'}t$.  Note $u'u\in T^Z={}^Z(T^*)$, by \cite[Proposition 8.3]{Bice2022}, so \eqref{az0} yields $t'\sim_{T^*}u'ut'<_tt'\in T^*$.  Also note $u'ut't\in T^ZT^Z\subseteq T^Z$ and $tu'ut'\in tT^Zt'\subseteq{}^ZT$ (again see \cite[Proposition 8.3]{Bice2022}) and hence
\begin{align*}
[t',T^*][t,T]&=[u'ut',T^*][t,T]=[u'ut't,\mathsf{s}(T)]\in\mathcal{F}[A]^0.\\
[t,T][t',T^*]&=[t,T][u'ut',T^*]=[tu'ut',\mathsf{r}(T)]\in\mathcal{F}[A]^0.
\end{align*}
This shows that if $t\in T\in\mathcal{F}(S)$ then $[t,T]$ is invertible with
\[[t,T]^{-1}=[t',T^*],\text{ for any $t'$ such that }u<_{t'}t,\text{ for some }u\in T.\]

Conversely, say we have $a\in S$ and $T\in\mathcal{F}(S)$ such that $[a,T]$ is invertible.  This means we have $a'\in S$ such that $[a',T^*][a,T]$ and $[a,T][a',T^*]$ are units.  So we have $y\in T^Z$ and $z\in{}^ZT$ such that $[a',T^*][a,T]=[y,\mathsf{s}(T)]$ and $[a,T][a',T^*]=[z,\mathsf{r}(T)]$.  Take $s,t,t',u,u'\in S$ with $T\in s<_{t'}t<_{u'}u$ and $uy=u$.  Replacing $t'$ with $t'tt'$ if necessary, we may assume that $t'<_uu'$ too.  Replacing $a$ and $a'$ with $\Phi(at')t$ and $\Phi(a't)t'$ if necessary (see the proof of \autoref{as<U0}), we may also assume $a<_{u'}u$ and $a'<_uu'$ and hence $aa',a'a\in\mathrm{ran}(\Phi)$, as $uu',u'u\in Z\subseteq\mathrm{ran}(\Phi)$.  As $[a',T^*][a,T]=[a'a,\mathsf{s}(T)]=[y,\mathsf{s}(T)]$, \eqref{rhoUnit} yields $q\in T^Z$ with $q\Phi(a'a)=q\Phi(y)$, i.e. $qa'a=qy$.  Likewise, we have $r\in{}^ZT$ with $aa'r=zr$.  Taking $m\in T^Z$ with $mq=m$ and $n\in{}^ZT$ with $nr=n$, we claim that
\[num<_{qa'r}a\]
Indeed, $qa'a=qy\in Z$ so $qa'ra\in qa'Za\subseteq Z$.  Likewise, $aqa'r\in aZa'r\subseteq Z$.  Also
\[numqa'r=numu'ua'r\in ZuZu'\mathrm{ran}(\Phi)Z\subseteq\mathrm{ran}(\Phi),\]
and, as $umqa'=uma'=umu'ua'\in uZu'\mathrm{ran}(\Phi)\subseteq\mathrm{ran}(\Phi)$ commutes with $r\in Z$,
\[numqa'ra=nrumqa'a=nrumqy=numy=nuym=num.\]
This proves the claim.  So $a>num\in{}^ZTTT^Z\subseteq\mathsf{r}(T)T\mathsf{s}(T)\subseteq T$ and hence $a\in T^<=T$.  This completes the proof of \eqref{FStimes}.
\end{proof}

For any $B\subseteq A$, let
\[\mathcal{F}[B]=\{[b,T]:b\in B\text{ and }T\in\mathcal{F}(S)\}.\]
Note \eqref{FS} implies that $\mathcal{F}[A]=\mathcal{F}[S]$, i.e. every equivalence class of $A$ w.r.t. every filter contains an element of $S$ (and even $S^>$).  In particular, for any $a\in A$, let
\[\mathcal{F}[a]=\mathcal{F}[\{a\}]=\{[a,T]:T\in\mathcal{F}(S)\}.\]
Also let $\mathsf{F}_{\langle A\rangle}:\mathcal{F}[A]\twoheadrightarrow\mathcal{F}(S)$ denote the canonical projection
\[\mathsf{F}_{\langle A\rangle}([a,T])=\mathsf{F}_{\langle A\rangle}(a_T,T^<)=T^<.\]
We consider the topology on $\mathcal{F}[A]$ generated by $(\mathcal{F}[a])_{a\in A}$ and $(\mathsf{F}_{\langle A\rangle}^{-1}[\mathcal{F}_s])_{s\in S}$.

\begin{prp}\label{FSbasis}
The sets $(a_s)_{a,s\in S^>}$ form a basis for $\mathcal{F}[A]$, where
\[a_s=\mathcal{F}[a]\cap\mathsf{F}_{\langle A\rangle}^{-1}[\mathcal{F}_s].\]
Also, if $a\in A$ and $T\subseteq T^<\in\mathcal{F}(S)$ then $(a_t)_{t\in T}$ is neighbourhood base at $[a,T]$.
\end{prp}

\begin{proof}
First observe that if $s,t\in S$ and $[a,T]\in\mathsf{F}_{\langle A\rangle}^{-1}[\mathcal{F}_s]\cap\mathsf{F}_{\langle A\rangle}^{-1}[\mathcal{F}_t]$ then $s,t\in T^<$ so we have $u\in T$ with $u<s,t$ and hence
\[[a,T]\in a_u\subseteq\mathsf{F}_{\langle A\rangle}^{-1}[\mathcal{F}_u]\subseteq\mathsf{F}_{\langle A\rangle}^{-1}[\mathcal{F}_s]\cap\mathsf{F}_{\langle A\rangle}^{-1}[\mathcal{F}_t].\]
To prove the result it thus suffices to show that, whenever $[a,T]\in\mathcal{F}[b]\cap\mathcal{F}[c]$, we can find $d\in T^>$ and $s\in T$ with $[a,T]\in a_s=d_s\subseteq\mathcal{F}[b]\cap\mathcal{F}[c]$.

Accordingly, take $a,b,c\in A$ and $T\subseteq T^<\in\mathcal{F}(S)$ with $[a,T]\in\mathcal{F}[b]\cap\mathcal{F}[c]$.  This means $a\sim_Tb,c$, so we have $b',c'\in T^*$ with $a\sim_{b'}b$ and $a\sim_{c'}c$.  As $T$ is directed, we have $v'\in T^*$ with $v'<b',c'$ and hence $a\sim_{v'}b,c$ (see \eqref{simBelow}).  Take $s,t,t',u,u',v\in S$ with $T\ni s<_{t'}t<_{u'}u<_{v'}v$.  Replacing $t'$ with $t'tt'$ if necessary, we may further assume that $t'<_uu'$.  Letting $d=\Phi(au')u<v$, we see that
\[\Phi(dt')=\Phi(\Phi(au')ut')=\Phi(au'ut')=\Phi(at').\]
Also $a\sim_{t'}b,c$ (again see \eqref{simBelow}).  For any $U\in\mathcal{F}_s$, it follows that $a\sim_Ub,c,d$ and hence $[a,T]\in a_s=d_s\subseteq\mathcal{F}[b]\cap\mathcal{F}[c]$, as required.
\end{proof}

Now we can show that filters and their equivalence classes form the kind of bundle that we are looking for.

\begin{thm}\label{0CB}
$\mathsf{F}_{\langle A\rangle}$ is an \'etale category bundle.
\end{thm}

\begin{proof}
We immediately see that that $\mathsf{F}_{\langle A\rangle}$ is continuous, as the topology on $\mathcal{F}[A]$ includes $\mathsf{F}_{\langle A\rangle}^{-1}[\mathcal{F}_s]$, for all $s\in S$.  On the other hand, $\mathsf{F}_{\langle A\rangle}[a_s]=\mathcal{F}_s$ is open in $\mathcal{F}(S)$, for all $a\in A$ and $s\in S$.  As these sets form a basis for $\mathcal{F}[A]$, by \autoref{FSbasis}, $\mathsf{F}_{\langle A\rangle}$ is also an open map.  As $\mathsf{F}_{\langle A\rangle}$ is injective $\mathcal{F}[a]$, for all $a\in A$, $\mathsf{F}_{\langle A\rangle}$ is also locally injective and hence a local homeomorphism.  From the definition of the product on $\mathcal{F}[A]$, we immediately see that $\mathsf{F}_{\langle A\rangle}$ is also an isocofibration.

Now take any $T,U\in\mathcal{F}(S)$ with $\mathsf{s}(T)=\mathsf{r}(U)$, $a<t\in T$ and $b<u\in U$.  For any $v\in T$ with $v<t$ and $w\in U$ with $w<u$, certainly $a_ub_v\subseteq(ab)_{uv}$.  As $((ab)_{vw})_{T\ni v<t,U\ni w<u}$ forms a neighbourhood base of $[ab,TU]=[a,T][b,U]$, by \autoref{FSbasis}, this shows that the product on $\mathcal{F}[A]$ is continuous.

To see that the source map $\mathsf{s}$ is also continuous on $\mathcal{F}[A]$, take any $a\in A$ and $T\in\mathcal{F}(S)$.  Further take $t\in T$ and $z\in Z$ with $tz=t$.  By \eqref{source}, $\mathsf{s}([a,T])=[z,\mathsf{s}(T)]$.  For any $y\in\mathsf{s}(T)$, we have $s\in T$ with $\mathsf{s}[\mathcal{F}_s]\subseteq\mathcal{F}_y$, as we already know that $\mathsf{s}$ is continuous on $\mathcal{F}(S)$.  Taking any $r\in T$ with $r<s,t$, we see that $\mathsf{s}[a_r]\subseteq z_y$.  As $(z_y)_{y\in\mathsf{s}[T]}$ is a neighbourhood base of $[z,\mathsf{s}(T)]$, by \autoref{FSbasis}, this shows that $\mathsf{s}$ is continuous.  Likewise, $\mathsf{r}$ is continuous so $\mathcal{F}[A]$ is a topological category.  As we already know that $\mathcal{F}(S)$ is an \'etale groupoid, $\mathsf{F}_{\langle A\rangle}$ is an \'etale category bundle.
\end{proof}

Actually, our primary interest will be in the \emph{ultrafilter bundle}, i.e. the subbundle $\mathsf{U}_{\langle A\rangle}=\mathsf{F}_{\langle A\rangle}|_{\mathcal{U}[A]}$ consisting of ultrafilters $\mathcal{U}(S)$ and their equivalence classes
\[\mathcal{U}[A]=\{[a,U]:a\in A\text{ and }U\in\mathcal{U}(S)\}.\]
Again the topology on $\mathcal{U}[A]$ is generated by the subbasis $\mathcal{U}_a=\mathcal{F}_a\cap\mathcal{U}(S)$ and $\mathcal{U}[a]=\{[a,U]:U\in\mathcal{U}(S)\}$, for $a\in A$, i.e. $\mathcal{U}[A]$ is a subspace of $\mathcal{F}[A]$.  As $\mathcal{U}(S)$ is also an ideal of $\mathcal{F}(S)$, \autoref{0CB} immediately yields the following.

Recall that $S$ is the semigroup part of our well-structured semimodule $(A,S,Z,\Phi)$.

\begin{cor}\label{piUEtaleZeroCategoryBundle}
If $S$ has a zero then $\mathsf{U}_{\langle A\rangle}$ is an \'etale zero category bundle.
\end{cor}

For any $a\in A$, we define $\widehat{a}:\mathcal{U}(S)\rightarrow\mathcal{U}[A]$ by
\[\widehat{a}(U)=[a,U].\]
Note that $\widehat{a}$ is a section of $\mathsf{U}_{\langle A\rangle}$, i.e. $\widehat{a}\in\mathcal{A}(\mathsf{U}_{\langle A\rangle})$.  In fact, more can be said.  First recall from \autoref{CategoryBundles} that $\mathcal{C}(\mathsf{U}_{\langle A\rangle})$ and $\mathcal{S}(\mathsf{U}_{\langle A\rangle})$ denote the continuous and slice-supported sections of $\mathsf{U}_{\langle A\rangle}$, while $\Phi^{\mathsf{U}_{\langle A\rangle}}:\mathcal{A}(\mathsf{U}_{\langle A\rangle})\rightarrow\mathcal{S}(\mathsf{U}_{\langle A\rangle})$ is the canonical expectation.

\begin{thm}\label{ContinuousSection}
If $S$ has a zero then, for any $a\in A$,
\begin{equation}\label{ContinuousPhi}
\widehat{a}\in\mathcal{C}(\mathsf{U}_{\langle A\rangle})\qquad\text{and}\qquad\widehat{\Phi(a)}=\Phi^{\mathsf{U}_{\langle A\rangle}}(\widehat{a}).
\end{equation}
If $<$ is also interpolative then, for any $s\in S^>$,
\[\widehat{s}\in\mathcal{S}(\mathsf{U}_{\langle A\rangle}),\quad\widehat{\,as\,}=\widehat{a}\widehat{s}\quad\text{and}\quad\widehat{\,sa\,}=\widehat{s}\widehat{a}.\]
\end{thm}

\begin{proof}
For any $U\in\mathcal{U}(S)$, $(a_u)_{u\in U}$ is a neighbourhood base of $\widehat{a}(U)=[a,U]$, by \autoref{FSbasis}.  As $\widehat{a}^{-1}[a_u]=\mathcal{U}_u$, for all $u\in U$, $\widehat{a}$ is continuous, i.e. $\widehat{a}\in\mathcal{C}(\mathsf{U}_{\langle A\rangle})$.

To see that $\widehat{\Phi(a)}=\Phi^{\mathsf{U}_{\langle A\rangle}}(\widehat{a})$ note, for any unit ultrafilter $U\in\mathcal{U}^0$, \eqref{rhoUnit} yields
\[\widehat{\Phi(a)}(U)=[\Phi(a),U]=[a,U]=\widehat{a}(U)=\Phi^{\mathsf{U}_{\langle A\rangle}}(\widehat{a})(U).\]
On the other hand, if $U\in\mathcal{U}(S)\setminus\mathcal{U}^0$ then $U^*$ is not a unit either so \autoref{UltrafilterDichotomy} yields $u\in U^*$ with $\Phi(u)=0$ and hence $\Phi(\Phi(a)u)=\Phi(a)\Phi(u)=0$.  Thus  $a\sim_U0$ so again $\widehat{\Phi(a)}(U)=[\Phi(a),U]=[0,U]=\Phi^{\mathsf{U}_{\langle A\rangle}}(\widehat{a})(U)$.  This proves $\widehat{\Phi(a)}=\Phi^{\mathsf{U}_{\langle A\rangle}}(\widehat{a})$.

Now assume $<$ is also interpolative and take $t,t',u,u'\in S$ with $s<_{t'}t<_{u'}u$.  If $U\in\mathcal{U}(S)$ and $\widehat{s}_\mathcal{U}(U)\neq0$, i.e. $s\not\sim_U0$, then $t\in s^<\subseteq s_U^<=U$, by \eqref{aUup}, i.e.
\begin{equation}\label{SliceSupport}
s<t\qquad\Rightarrow\qquad\mathrm{supp}(\widehat{s})\subseteq\mathcal{U}_t.
\end{equation}
By \cite[Proposition 7.3]{Bice2022}, $\mathcal{U}_t$ is a slice so $\widehat{s}\in\mathcal{S}(\mathsf{U}_{\langle A\rangle})$.

It follows that $\mathsf{r}[\mathrm{supp}(\widehat{s}\widehat{a})]\subseteq\mathsf{r}[\mathrm{supp}(\widehat{s})]\subseteq\mathsf{r}[\mathcal{U}_t]$.  We claim $\mathsf{r}[\mathrm{supp}(\!\widehat{\,sa\,}\!)]\subseteq\mathsf{r}[\mathcal{U}_t]$ too.  To see the this, first replace $t'$ with $t'tt'$ if necessary to ensure that $t'<u'$.  Now take $W\in\mathcal{U}(S)$ with $\mathsf{r}(W)\notin\mathsf{r}[\mathcal{U}_t]=\mathcal{U}_{tu'}$, by \eqref{rFt}.  Then $st'<tu'$, by \cite[Proposition 4.7]{Bice2022}, so interpolation yields a sequence $(q_n)\subseteq\mathrm{ran}(\Phi)$ such that
\[st'<q_{n+1}<q_n<tu',\]
for all $n\in\mathbb{N}$.  Let $Q=\{q_n:n\in\mathbb{N}\}$ and note that $(Q\mathsf{r}(W))^<$ is a filter containing $tu'$ and $\mathsf{r}(W)$.  As $\mathsf{r}(W)\in\mathcal{U}(S)\setminus\mathcal{U}_{tu'}$, it follows that $(Q\mathsf{r}(W))^<=S$ and hence $qr=0$, for some $q\in Q$ and $r\in\mathsf{r}(W)$.  As $st'<q$ and $\mathsf{r}(W)=(WW^*)^<$, it follows that $st'w=0$, for some $w\in W$.  Taking $w',x\in S$ with $W\ni x<_{w'}w$, we see that
\[\Phi(w'ww'sa)=\Phi(w'ww'st'ta)=\Phi(w'st'ww'ta)=0\]
and hence $\widehat{\,sa\,}(W)=[sa,W]=[0,W]$.  This proves the claim.

Next note that, for any $U\in\mathcal{U}_t$ and $V\in\mathcal{U}(S)$ with $\mathsf{s}(U)=\mathsf{r}(V)$,
\[\widehat{s}\widehat{a}(U\cdot V)=\widehat{s}(U)\widehat{a}(V)=[s,U][a,V]=[sa,UV]=\widehat{\,sa\,}(U\cdot V),\]
as $s<t\in U$ and hence $s\in U^>$.  On the other hand, $\widehat{s}\widehat{a}(W)=[0,W]=\widehat{\,sa\,}(W)$, for any $W\in\mathcal{U}(S)\setminus\mathsf{r}[\mathcal{U}_t]$, by the claim proved above.  This shows that $\widehat{\,sa\,}=\widehat{s}\widehat{a}$, while a dual argument also yields $\widehat{\,as\,}=\widehat{a}\widehat{s}$.
\end{proof}

The above result says that we can represent $A$ more concretely as continuous sections $\widehat{A}=\{\widehat{a}:a\in A\}$ of the ultrafilter bundle $\mathsf{U}_{\langle A\rangle}$, at least if we assume $<$ is interpolative and restrict multiplication to elements of $S^>$.  Accordingly, the map $a\mapsto\widehat{a}$ will be called the \emph{ultrafilter representation}.  To say more about the ultrafilter representation, we must first consider conditions relating to suprema within the semigroup part of our well-structured semimodules, which we now examine.

\section{Steinberg Semigroups}\label{SteinbergSemigroups}

Given a semigroup $S$ with subsemigroup $Z$, we again consider the restriction relation $\leq$ from \autoref{Restriction}.  As usual, we call $a\vee b\in S$ a \emph{supremum} of $a,b\in S$ if
\[a\vee b=\min\{s\in S:a,b\leq s\},\]
i.e. $a,b\leq a\vee b$ and $a\vee b\leq s$ whenever $a,b\leq s$.  As restriction is antisymmetric, suprema are unique, whenever they exist.  As restriction is reflexive on $S^\geq$, a supremum $a\vee b$ can equivalently characertised as an element of $S$ satisfying
\[a^\leq\cap b^\leq=(a\vee b)^\leq\neq\emptyset.\]

If $S$ has a zero and $a\leq b$ then a \emph{$b$-complement} of $a$ is an element $b\setminus a$ such that
\[a(b\setminus a)=0\qquad\text{and}\qquad a\vee(b\setminus a)=b.\]

\begin{dfn}
A \emph{Steinberg semigroup} is a well-structured semigroup $(S,Z,\Phi)$ where $S=S^>$ has a zero, $Z\subseteq\mathsf{E}(S)$ and, for all $a,b\in S$, $r\in\mathrm{ran}(\Phi)$ and $y,z\in Z$,
\begin{align}
\tag{Orthosuprema}\label{Orthosuprema}a\perp b\qquad&\Rightarrow\qquad a\text{ and $b$ have a supremum }a\vee b.\\
\tag{Distributivity}\label{Distributivity}y\perp z\qquad&\Rightarrow\qquad y\vee z\in Z\text{ and }r(y\vee z)=ry\vee rz.\\
\tag{Complements}\label{Complements}y\leq z\qquad&\Rightarrow\qquad y\text{ has a $z$-complement }z\setminus y\in Z.
\end{align}
\end{dfn}

Incidentally, instead of requiring $Z\subseteq\mathsf{E}(S)$, we could have required that $(Z,\leq)$ is up-directed or that $(Z,<)$ is a poset, by \autoref{leqUpDirected} and \autoref{InverseStructuredSemigroups}.  Also note that the zero of $S$ must necessarily lie in $Z$.  Indeed, as semigroups are non-empty by definition, we can take $z\in Z$ and note that $z\leq z$ so we have $z\setminus z\in Z$ with $0=z(z\setminus z)\in ZZ\subseteq Z$

\subsection{Examples}

Recall from \cite{Lawson2012} that an inverse semigroup is \emph{distributive} if compatible pairs have joins/suprema (w.r.t. the canonical order) which are distributive relative to products, i.e. $a(b\vee c)=ab\vee ac$.  A distributive inverse $\wedge$-semigroup is \emph{Boolean} if its idempotents also form a generalised Boolean algebra.

\begin{prp}
If $S$ is a Boolean inverse $\wedge$-semigroup then $\langle S\rangle=(S,\mathsf{E}(S),\Phi_\mathsf{E})$ is a Steinberg semigroup.
\end{prp}

\begin{proof}
We already know that $\langle S\rangle$ is well-structured, by \autoref{WellStructuredInverseSemigroup}.  As $S$ is an inverse semigroup, we also have $S=S^\dagger=S^>$.  As orthogonal elements are compatible and hence have joins, \eqref{Orthosuprema} is satisfied.  As $S$ is distributive, \eqref{Distributivity} holds.  As $\mathsf{E}(S)$ is a generalised Boolean algebra, \eqref{Complements} also holds, so $\langle S\rangle$ satisfies all the required conditions of a Steinberg semigroup.
\end{proof}

In fact, for inverse $\wedge$-semigroups $S$, the converse also holds in that if $(S,\mathsf{E}(S),\Phi_\mathsf{E})$ is a Steinberg semigroup then $S$ is Boolean -- this follows from \autoref{SteinbergSemigroupBooleanZ} and \eqref{SDistributivity} below and the fact that compatible joins can be defined from orthogonal joins whenever the idempotents have relative complements.

More general examples of Steinberg semigroups come from ample category bundles.  Indeed, we will soon see in \autoref{FaithfulUltrafilters} and \autoref{AmpleUltrafilters} below that such examples actually encompass all Steinberg semigroups, at least up to isomorphism.

\begin{thm}\label{AmpleBundleSteinbergSemigroup}
If $\rho:C\twoheadrightarrow G$ is an ample category bundle then $(\mathcal{S}_\mathsf{c}(\rho),\mathcal{Z}_\mathsf{c}(\rho),\Phi_\mathsf{c}^\rho)$ is a Steinberg semigroup such that, for all $a,b\in\mathcal{S}_\mathsf{c}(\rho)$,
\begin{align}
\label{AmpleRestriction}a\leq b\quad\Leftrightarrow\quad&a|_{\mathrm{supp}(a)}=b|_{\mathrm{supp}(a)}.\\
\label{AmpleDomination}a<b\quad\Leftrightarrow\quad&\mathrm{supp}(a)\subseteq b^{-1}[C^\times].\\
\label{AmpleOrthogonality}a\perp b\quad\Leftrightarrow\quad&\mathsf{r}[\mathrm{supp}(a)]\cap\mathsf{r}[\mathrm{supp}(b)]=\emptyset=\mathsf{s}[\mathrm{supp}(a)]\cap\mathsf{s}[\mathrm{supp}(b)]\\
\nonumber\Leftrightarrow\quad&\mathrm{supp}(a)\cap\mathrm{supp}(b)=\emptyset\text{ and }\mathrm{supp}(a)\cup\mathrm{supp}(b)\in\mathcal{B}(G).
\end{align}
Moreover, $a$ is bisupported precisely when $\mathrm{supp}(a)$ is open $($and hence clopen$)$.
\end{thm}

\begin{proof}
Take an ample category bundle $\rho:C\twoheadrightarrow G$, so $(\mathcal{S}_\mathsf{c}(\rho),\mathcal{Z}_\mathsf{c}(\rho),\Phi_\mathsf{c}^\rho)$ is at least a well-structured semigroup, by \autoref{ContinuousSubactions}.

\begin{itemize}
\item[\eqref{AmpleRestriction}]  If $a\leq b$ then we have $z\in\mathcal{Z}_\mathsf{c}(\rho)$ with $a=az=bz$.  For all $g\in\mathrm{supp}(a)$ this means $a(g)=a(g)z(\mathsf{s}(g))=b(g)z(\mathsf{s}(g))$.  As $a(g)\neq0_g$, this implies $z(\mathsf{s}(g))\in C^0$ and hence $a(g)=b(g)$.  This shows that $a|_{\mathrm{supp}(a)}=b|_{\mathrm{supp}(a)}$.

Conversely, say $a|_{\mathrm{supp}(a)}=b|_{\mathrm{supp}(a)}$.  Then $O=\rho(\mathrm{ran}(a)\cap\mathrm{ran}(b))$ is an open subset containing the compact subset $\mathrm{supp}(a)$, by \autoref{EtaleBundleOpenRange} and \autoref{OpenClosedSupports}.  As $G$ is ample, we have compact clopen $K$ with $\mathrm{supp}(a)\subseteq K\subseteq O$.  Then we have $y,z\in\mathcal{Z}_\mathsf{c}(\rho)$ with $\mathrm{supp}(y)=\mathsf{r}[K]$ and $\mathrm{supp}(z)=\mathsf{s}[K]$ and hence $yb=ya=a=az=bz$, showing that $a\leq b$.

\item[\eqref{AmpleDomination}]  Say $a<_sb$.  For all $g\in\mathrm{supp}(a)$, this means $a(g)=a(g)s(h)b(i)$, for some $h,i\in G$ with $g=ghi$.  In particular, $as(gh)=a(g)s(h)\neq0_{gh}$ so $gh\in G^0$, as $as\in\mathrm{ran}(\Phi_\mathsf{c}^\rho)$.  Then $g=i$ so $sb(\mathsf{s}(g))=s(g^{-1})b(g)=s(h)b(i)\in C^0$, as $sb\in\mathcal{Z}_\mathsf{c}(\rho)$.  Likewise, $bs(\mathsf{r}(g))=b(g)s(g^{-1})\in C^0$ and hence $b(g)\in C^\times$.  This shows that $\mathrm{supp}(a)\subseteq b^{-1}[C^\times]$.

Conversely, say $\mathrm{supp}(a)\subseteq b^{-1}[C^\times]$.  Then $\mathrm{supp}(a)$ is compact, again by \autoref{OpenClosedSupports}, and $b^{-1}[C^\times]$ open, by \autoref{OpenInvertibles}.  As $G$ is ample, we have compact clopen $K$ with $\mathrm{supp}(a)\subseteq K\subseteq b^{-1}[C^\times]$.  We can then define $s\in\mathcal{S}_\mathsf{c}(\rho)$ such that $a<_sb$ by
\[s(g)=\begin{cases}b(g^{-1})^{-1}&\text{if }g\in K^{-1}\\ 0_g&\text{otherwise}.\end{cases}\]

\item[\eqref{AmpleOrthogonality}] If $a\perp b$ then we have $y,z\in\mathcal{Z}_\mathsf{c}(\rho)$ with $ya=a=az$ and $yb=0=bz$ so
\begin{align*}
\mathsf{r}[\mathrm{supp}(a)]&\subseteq\mathrm{supp}(y)\subseteq G^0\setminus\mathsf{r}[\mathrm{supp}(b)]\quad\text{and}\\
\mathsf{s}[\mathrm{supp}(a)]&\subseteq\mathrm{supp}(z)\subseteq G^0\setminus\mathsf{s}[\mathrm{supp}(b)].
\end{align*}
In particular, $\mathsf{r}[\mathrm{supp}(a)]\cap\mathsf{r}[\mathrm{supp}(b)]=\emptyset=\mathsf{s}[\mathrm{supp}(a)]\cap\mathsf{s}[\mathrm{supp}(b)]$.

Conversely, say $\mathsf{r}[\mathrm{supp}(a)]\cap\mathsf{r}[\mathrm{supp}(b)]=\emptyset$.  As $\mathsf{r}[\mathrm{supp}(a)]$ and $\mathsf{r}[\mathrm{supp}(b)]$ are compact subsets, again by \autoref{OpenClosedSupports}, we have compact clopen $K$ with $\mathsf{r}[\mathrm{supp}(a)]\subseteq K\subseteq G^0\setminus\mathsf{r}[\mathrm{supp}(b)]$.  Taking $y\in\mathcal{Z}_\mathsf{c}(\rho)$ with $\mathrm{supp}(y)=\mathsf{r}[K]$, it follows that $ya=a$ and $yb=0$.  If $\mathsf{s}[\mathrm{supp}(a)]\cap\mathsf{s}[\mathrm{supp}(b)]=\emptyset$ too then we likewise obtain $z\in\mathcal{Z}_\mathsf{c}(\rho)$ with $az=a$ and $bz=0$ and hence $a\perp b$.

This proves the first equivalence.  Regarding the second, note that if $\mathsf{r}[\mathrm{supp}(a)]\cap\mathsf{r}[\mathrm{supp}(b)]=\emptyset=\mathsf{s}[\mathrm{supp}(a)]\cap\mathsf{s}[\mathrm{supp}(b)]$ then we immediately see that $\mathrm{supp}(a)\cap\mathrm{supp}(b)=\emptyset$ and also $\mathrm{supp}(a)\cup\mathrm{supp}(b)\in\mathcal{B}(G)$, as $\mathrm{supp}(a),\mathrm{supp}(b)\in\mathcal{B}(G)$.  Conversely, take $g\in\mathrm{supp}(a)$ and $h\in\mathrm{supp}(b)$.  If $g=h$ then this witnesses $\mathrm{supp}(a)\cap\mathrm{supp}(b)\neq\emptyset$.  Otherwise $g\neq h$ so if either $\mathsf{r}(g)=\mathsf{r}(h)$ or $\mathsf{s}(g)=\mathsf{s}(h)$ then $\mathrm{supp}(a)\cup\mathrm{supp}(b)\notin\mathcal{B}(G)$.
\end{itemize}

To see that $(\mathcal{S}_\mathsf{c}(\rho),\mathcal{Z}_\mathsf{c}(\rho),\Phi_\mathsf{c}^\rho)$ is a Steinberg semigroup, first note that $\mathcal{Z}_\mathsf{c}(\rho)$ consists of characteristic functions of compact clopen subsets of $G^0$, which are all idempotents.  This includes the zero section, which is indeed a zero for the whole semigroup $\mathcal{S}_\mathsf{c}(\rho)$.  Also, any $a\in\mathcal{S}_\mathsf{c}(\rho)$ has compact support so, as $G=\rho[C^\times]$ is ample, we have $b\in\mathcal{S}_\mathsf{c}(\rho)$ with $\mathrm{ran}(b)\subseteq C^\times$ and $\mathrm{supp}(a)\subseteq\mathrm{supp}(b)$.  Then $a<_{b^{-1}}b$, showing that $\mathcal{S}_\mathsf{c}(\rho)=\mathcal{S}_\mathsf{c}(\rho)^>$.

For \eqref{Orthosuprema}, take $a,b\in\mathcal{S}_\mathsf{c}(\rho)$ with $a\perp b$.  As $\mathrm{supp}(a)\cap\mathrm{supp}(b)=\emptyset$, we can define another function $s$ by
\[s(g)=\begin{cases}a(g)&\text{if }g\in\mathrm{supp}(a)\\b(g)&\text{if }g\in\mathrm{supp}(b)\\0_g&\text{otherwise}.\end{cases}\]
Note we can partition $G$ into compact clopen sets on which $s$ agrees with $a$, $b$, and the zero section respectively.  As these are all continuous sections, so is $s$.  Also $\mathrm{supp}(s)=\mathrm{supp}(a)\cup\mathrm{supp}(b)$ is a compact slice so $s\in\mathcal{S}_\mathsf{c}(\rho)$.  From \eqref{AmpleRestriction}, we immediately see that $s$ is a supremum of $a$ and $b$, i.e. $s=a\vee b$.

In particular, for any orthogonal $y,z\in\mathcal{Z}_\mathsf{c}(\rho)$, we see that
\[(y\vee z)[\mathrm{supp}(y\vee z)]=y[\mathrm{supp}(y)]\cup z[\mathrm{supp}(z)]\subseteq C^0,\]
i.e. $y\vee z\in\mathcal{Z}_\mathsf{c}(\rho)$.  For any $r\in\mathrm{ran}(\Phi_\mathsf{c}^\rho)$, elementary calculations also show that $r(y\vee z)=ry\vee rz$, so \eqref{Distributivity} holds.  Finally, if $y,z\in\mathcal{Z}_\mathsf{c}(\rho)$ and $y\leq z$ then $\mathrm{supp}(y)\subseteq\mathrm{supp}(z)$ and hence $y$ has a $z$-complement given by $z\setminus y=\mathbf{1}_{\mathrm{supp}(z)\setminus\mathrm{supp}(z)}$.  The shows that \eqref{Complements} also holds and hence $(\mathcal{S}_\mathsf{c}(\rho),\mathcal{Z}_\mathsf{c}(\rho),\Phi_\mathsf{c}^\rho)$ is indeed a Steinberg semigroup.

If $\mathrm{supp}(a)$ is open then so is $\mathsf{r}[\mathrm{supp}(a)]$ and hence we have $z\in\mathcal{Z}_\mathsf{c}(\rho)$ with $\mathrm{supp}(z)=\mathsf{r}[\mathrm{supp}(a)]$, which is immediately seen to be a range-support of $a$.  Likewise, if $\mathrm{supp}(a)$ is open then we also have $z\in\mathcal{Z}_\mathsf{c}(\rho)$ with $\mathrm{supp}(z)=\mathsf{s}[\mathrm{supp}(a)]$, which is source-support of $a$, showing that $a$ is bisupported.  Conversely, if $\mathrm{supp}(a)$ is not open then, for any $z\in\mathcal{Z}_\mathsf{c}(\rho)$ with $za=a$, $\mathrm{supp}(z)\setminus\mathsf{r}[\mathrm{supp}(a)]$ is a non-empty open subset of $G^0$.  As $G$ is ample, when then have $y\in\mathcal{Z}_\mathsf{c}(\rho)$ with $\mathrm{supp}(y)\subseteq\mathrm{supp}(z)\setminus\mathsf{r}[\mathrm{supp}(a)]$, which means that $z\setminus y$ is a strictly smaller element of $Z$ with $(z\setminus y)a=a$.  This shows that $a$ has no range-support so $a$ is not bisupported.
\end{proof}

In particular, every $a\in\mathcal{S}_\mathsf{c}(\rho)$ is bisupported when $C$ above is also Hausdorff, thanks to \autoref{OpenClosedSupports}, which means that $\mathcal{S}_\mathsf{c}(\rho)$ is a Boolean restriction semigroup with projections $\mathcal{Z}_\mathsf{c}(\rho)$, as in \cite[Theorem 8.19]{KudryavtsevaLawson2017}.  However, this is also possible even when $C$ is not Hausdorff, as the following example shows (this discrepancy will disappear when we restrict our attention to Steinberg rings -- see \autoref{SectionRings} below).

\begin{xpl}
Let $G=G^0=\{0\}\cup\{1/n:n\in\mathbb{N}\}$, with the usual subspace topology, and consider the unique ample category bundle $\rho:C\twoheadrightarrow G$ such that $\rho^{-1}\{1/n\}=\{0_{1/n},1_{1/n}\}$, for all $n\in\mathbb{N}$, $\rho^{-1}\{0\}=\{0_0,e,1_0\}$, $ee=e$, and $\rho$ is a homeomorphism when restricted to $\{e\}\cup\{1_{1/n}:n\in\mathbb{N}\}$.  Then $e$ and $1_0$ can not be separated by disjoint open sets, however the only continuous sections taking the value $e$ must be non-zero in a neighbourhood of $e$.  It follows that every continuous section of $\rho$ has open support.
\end{xpl}

\subsection{General Theory}

Let us now make the following standing assumption.
\begin{center}
\textbf{$\langle S\rangle=(S,Z,\Phi)$ is a Steinberg semigroup.}
\end{center}
For all $y,z\in Z\subseteq\mathsf{E}(S)$, note $y^\dagger=\{y\}$ so \eqref{Inverseaperpb} says that
\[y\perp z\qquad\Leftrightarrow\qquad yz=0.\]
Likewise, restriction and domination agree on $Z$, by \eqref{zyyyz} and \eqref{EZdomination}, i.e.
\begin{equation}\label{DominationRestrictionUpDirected}
y\leq z\qquad\Leftrightarrow\qquad y<z\qquad\Leftrightarrow\qquad y=yz.
\end{equation}
As $S=S^>$ too, $Z$ is a poset, by \eqref{RestrictionPoset}.  In fact, $Z$ is even a lattice and hence a generalised Boolean algebra, as $Z$ is also distributive and relatively complemented.

\begin{prp}\label{SteinbergSemigroupBooleanZ}
$Z$ is a lattice and hence a generalised Boolean algebra.
\end{prp}

\begin{proof}
First note that $yz\leq y,z$, for any $y,z\in Z$, as $yz=yyz=yzz$.  Moreover, $x\leq yz$ whenever $Z\ni x\leq y,z$, as this implies $x=xz=xyz$.  This shows that $xy$ is the infimum of $x$ and $y$.  In particular, we have an $x$-complement $x\setminus xy$ which is orthogonal to $y$, as $(x\setminus xy)y=(x\setminus xy)xy=0$.  Thus we have a supremum $(x\setminus xy)\vee y\geq(x\setminus xy)\vee xy=x$ and hence $(x\setminus xy)\vee y$ is also the supremum of $x$ and $y$.  This shows that $Z$ is indeed a lattice.
\end{proof}

In particular, relative complements in $Z$ are unique.  Indeed, if $x,y,z,z'\in Z$, $x\leq y$ and both $z$ and $z'$ are $y$-complements of $x$ then $xz=0$ so \eqref{Distributivity} yields $z=yz=(x\vee z')z=xz\vee zz'=zz'$.  Likewise $z'=zz'$ so $z=z'$.

As usual, we extend $\setminus$ to arbitrary $y,z\in Z$ by defining
\[y\setminus z=y\setminus yz=(y\vee z)\setminus z.\]

\begin{lem}
For any $r\in\mathrm{ran}(\Phi)$ and $y,z\in Z$,
\begin{align}
\label{ryz}r<y\setminus z\qquad&\Leftrightarrow\qquad z\perp r<y.\\
\label{yzr}(y\setminus z)\perp r\qquad&\Leftrightarrow\qquad ry<z.
\end{align}
\end{lem}

\begin{proof}
Take $r\in\mathrm{ran}(\Phi)$ and $y,z\in Z$.  If $r<y\setminus z$ then $ry=r(y\setminus z)y=r(y\setminus z)=r$ and $rz=r(y\setminus z)z=0$, i.e. $z\perp r<y$.  Conversely, if $z\perp r<y$ then
\[r=ry=r((y\setminus z)\vee z)=r(y\setminus z)\vee rz=r(y\setminus z)\vee0=r(y\setminus z),\]
i.e. $r<y\setminus z$.  This proves \eqref{ryz}.

If $ry<z$ then $r(y\setminus z)=ry(y\setminus z)=ryz(y\setminus z)=0$, i.e. $y\setminus z\perp r$.  Conversely, if $(y\setminus z)\perp r$ then $ry=r((y\setminus z)\vee yz)=0\vee ryz=ryz$, i.e. $ry<z$, proving \eqref{yzr}.
\end{proof}

Now we can show that $\perp$ has the following more symmetric characterisation.

\begin{prp}\label{perpSymmetric}
For all $a,b\in S$,
\[a\perp b\quad\Leftrightarrow\quad\exists y,y',z,z'\in Z\ (ya=a=ay',\ zb=b=bz'\text{ and }yz=0=y'z').\]
\end{prp}

\begin{proof}
If the right side holds then $yb=yzb=0=bz'y'=b'y'$ so $a\perp b$.

Conversely, if $a\perp b$ then, in particular, we have $y\in Z$ with $ya=a$ and $yb=0$.  Take $s,s'\in S$ with $b<_{s'}s$ and let $z=ss'\setminus y$.  Note $bs'<ss'$ and $ybs'=0$ and hence $bs'<z$, by \eqref{ryz}.  Thus $zb=zbs's=bs's=b$ and $yz=y(ss'\setminus y)=0$.  A dual argument yields $z'\in Z$ with $b=bz'$ and $y'z'=0$, as required.
\end{proof}

Next let us show that distributivity extends to arbitrary elements of $S$.

\begin{prp}
For all $a,b,c\in S$ such that $b$ and $c$ have a supremum $b\vee c$,
\begin{equation}\label{SDistributivity}
a(b\vee c)=ab\vee ac.
\end{equation}
\end{prp}

\begin{proof}
As $b,c\leq b\vee c$, \eqref{Multiplicativity} yields $ab,ac\leq a(b\vee c)$.  Thus to prove $a(b\vee c)=ab\vee ac$, it suffices to show that $a(b\vee c)\leq d$ whenever $ab,ac\leq d$.

First assume $a\in\mathrm{ran}(\Phi)$ and $b,c\in Z$.  If $ab,ac\leq d$ then $a(b\setminus c)\leq ab\leq d$, by \eqref{Multiplicativity}.  Then \eqref{Distributivity} yields
\[a(b\vee c)=a((b\setminus c)\vee c)=a(b\setminus c)\vee ac\leq d.\]

Now assume $a\in S$ and $b,c\in Z$.  Take $s,s'\in S$ with $a<_{s'}s$.  If $ab,ac\leq d$ then $s'ab,s'ac\leq s'd$, by \eqref{Multiplicativity}, and hence what we just proved yields
\[a(b\vee c)=ss'a(b\vee c)=s(s'ab\vee s'ac)\leq ss'd\leq d.\]
This proves $a(b\vee c)=ab\vee ac$ and a dual argument yields $(b\vee c)a=ba\vee ca$.

Now take arbitrary $a,b,c\in S$.  As $b,c\leq b\vee c$, we have $y,z\in Z$ with $b=(b\vee c)y$ and $c=(b\vee c)z$.  Note $b=by=b(y\vee z)\leq(b\vee c)(y\vee z)$ and, likewise, $c\leq(b\vee c)(y\vee z)$ so $b\vee c\leq(b\vee c)(y\vee z)\leq b\vee c$, i.e. $b\vee c=(b\vee c)(y\vee z)$ and hence
\[a(b\vee c)=a(b\vee c)(y\vee z)=a(b\vee c)y\vee a(b\vee c)z=ab\vee ac,\]
by what we just proved.
\end{proof}

It follows that orthogonality respects suprema.

\begin{cor}
If $a,b\in S$ have a supremum $a\vee b$ then, for all $c\in S$,
\begin{equation}\label{veeperp}
a,b\perp c\qquad\Rightarrow\qquad a\vee b\perp c.
\end{equation}
\end{cor}

\begin{proof}
If $a,b\perp c$ then we have $y,y',z,z'\in Z$ with $ya=a=ay'$, $zb=b=bz'$ and $yc=zc=0=cy'=cz'$.  Then \eqref{SDistributivity} yields $(y\vee z)(a\vee b)=a\vee b=(a\vee b)(y'\vee z')$ and $(y\vee z)c=yc\vee zc=0=cy'\vee cz'=c(y'\vee z')$, showing that $a\vee b\perp c$.
\end{proof}

It also follows that $\Phi$ preserves suprema.

\begin{cor}
If $a,b\in S$ have a supremum $a\vee b$ then
\begin{equation}\label{veePhi}
\Phi(a\vee b)=\Phi(a)\vee\Phi(b)
\end{equation}
\end{cor}

\begin{proof}
Take $s,s'\in S$ with $a\vee b<_{s'}s$.  As $a,b\leq a\vee b$, \eqref{LeftAuxiliarity} yields $a,b<_{s'}s$.  Now \eqref{SDistributivity} yields
\[\Phi(a\vee b)=\Phi((a\vee b)s's)=(a\vee b)s'\Phi(s)=as'\Phi(s)\vee bs'\Phi(s)=\Phi(a)\vee\Phi(b).\qedhere\]
\end{proof}

For any $T,U\subseteq S$, let
\[T\vee U=\{t\vee u:t\in T\text{ and }u\in U\}.\]
In particular, $T\vee\emptyset=\emptyset=\emptyset\vee T$.

Now we show that orthogonal suprema also preserve $Z$-inverses.

\begin{cor}\label{veeInverse}
For any orthogonal $a,b\in S$,
\begin{equation}\label{veeInverseEq}
(a\vee b)^\dagger=a^\dagger\vee b^\dagger.
\end{equation}
\end{cor}

\begin{proof}
First note that if $a,b\in S$ have a $Z$-invertible supremum $a\vee b\in S^\dagger$ then $a,b<a\vee b$, by \eqref{LeftAuxiliarity}, and hence $a$ and $b$ are also $Z$-invertible, by \eqref{CommonRestrictionDomination}.  Thus if $a^\dagger=\emptyset$ or $b^\dagger=\emptyset$ then $(a\vee b)^\dagger=\emptyset$.

Conversely, say $a^\dagger=\{a'\}$ and $b^\dagger=\{b'\}$.  If $a\perp b$ then $a'b=0=ba'$, by \eqref{Inverseaperpb}, and hence $ab'=0=b'a$, as $\{ab\}^\dagger=b^\dagger a^\dagger$, so $a'\perp b'$, again by \eqref{Inverseaperpb}.  Thus we have suprema $a\vee b$ and $a'\vee b'$, by \eqref{Orthosuprema}.  Then \eqref{Inversealeqb} yields
\[ab'=aa'(a\vee b)b'bb'=aa'bb'\leq aa',bb'.\]
Thus $(a\vee b)a'=aa'\vee ba'=aa'$ and $(a\vee b)b'=bb'$ and hence
\[(a\vee b)(a'\vee b')=(a\vee b)a'\vee(a\vee b)b'=aa'\vee bb'\in Z.\]
Also $(a\vee b)(a'\vee b')(a\vee b)=(aa'\vee bb')(a\vee b)=aa'(a\vee b)\vee bb'(a\vee b)=a\vee b$, again by \eqref{Inversealeqb}.  Likewise, $(a'\vee b')(a\vee b)=a'a\vee b'b\in Z$ and $(a'\vee b')(a\vee b)(a'\vee b')=a'\vee b'$, showing that $(a\vee b)^\dagger=\{a'\vee b'\}$.
\end{proof}

For convenience, let $a^{-1}$ denote the unique $Z$-inverse of $a$, whenever one exists, i.e. $a^\dagger=\{a^{-1}\}$, for all $a\in S^\dagger$.  So, for all orthogonal $a,b\in S^\dagger$, \eqref{veeInverseEq} says that
\[(a\vee b)^{-1}=a^{-1}\vee b^{-1}.\]

\subsection{Ultrafilters}

The first thing to note about ultrafilters in Steinberg semigroups is that they are determined by their $Z$-invertible elements.  Indeed, thanks to \eqref{ReflexiveInterpolation},
\[U\mapsto U^{\dagger\dagger}(=U\cap S^\dagger)\qquad\text{and}\qquad U\mapsto U^{**}(=U^<)\]
are mutually inverse homeomorphisms between $\mathcal{F}(S)$ and $\mathcal{F}(S^\dagger)(=$ proper filters in $S^\dagger$ w.r.t. $<)$, which then restrict to homeomorphisms between $\mathcal{U}(S)$ and $\mathcal{U}(S^\dagger)$.  We can then show that ultrafilters are precisely those proper non-empty filters that are prime with respect to orthogonal $Z$-invertible elements.

\begin{thm}
$U\in\mathcal{F}(S)\setminus\{S\}$ is an ultrafilter if and only if, for all $a,b\in S^\dagger$,
\begin{equation}\label{Ultra+}
a\perp b\quad\text{and}\quad a\vee b\in U\qquad\Rightarrow\qquad a\in U\quad\text{or}\quad b\in U.
\end{equation}
\end{thm}

\begin{proof}
If $U\in\mathcal{U}(S)$ and $a,b\in S^\dagger\setminus U$ then, by \eqref{UltrafilterReflexives}, we have $u,v\in U$ with $\Phi(ua^{-1})=0=\Phi(vb^{-1})$.  Taking $w\in U$ such that $w<u,v$, it follows that $\Phi(wa^{-1})=0=\Phi(wb^{-1})$.  If $a\perp b$ then \eqref{SDistributivity}, \eqref{veePhi} and \eqref{veeInverseEq} yield
\[\Phi(w(a\vee b)^{-1})=\Phi(wa^{-1}\vee wb^{-1})=\Phi(wa^{-1})\vee\Phi(wb^{-1})=0.\]
Thus $a\vee b\notin U$, again by \eqref{UltrafilterReflexives}.

Conversely, say $T\in\mathcal{F}(S)\setminus\{S\}$ is not an ultrafilter, so $T$ is contained in a strictly larger ultrafilter $U$.  Take $u\in U^{\dagger\dagger}\setminus T$ and $t\in T^{\dagger\dagger}$.  As $U$ is a filter, we can further take $s\in U^{\dagger\dagger}$ with $s<t,u$ and let $a=t(t^{-1}t\setminus s^{-1}s)$ and $b=ts^{-1}s$ so $a\perp b$ and
\[a\vee b=t(t^{-1}t\setminus s^{-1}s)\vee ts^{-1}s=tt^{-1}t=t\in T.\]
Note that $sa^{-1}=ss^{-1}s(t^{-1}t\setminus s^{-1}s)t^{-1}=0$ and hence $a\notin U\supseteq T$, by \eqref{UltrafilterReflexives}.  Also $ts^{-1}\in\mathrm{ran}(\Phi)$ so $b=ts^{-1}s<u\notin T$ and hence $b\notin T$, i.e. \eqref{Ultra+} fails for $T$.
\end{proof}

Next we show that the ultrafilters in a Steinberg semigroup are compactly based.

\begin{prp}\label{RSbasis}
The ultrafilter groupoid $\mathcal{U}(S)$ has compact open basis $(\mathcal{U}_r)_{r\in S^\dagger}$.
\end{prp}

\begin{proof}
First we claim that $(\mathcal{U}_z)_{z\in Z}$ is a compact open basis of $\mathcal{U}^0$ and, moreover,
\[z\mapsto\mathcal{U}_z\]
is a Boolean homomorphism from $Z$ onto the compact open subsets of $\mathcal{U}^0$.  Indeed, this essentially follows from the classic Stone duality.  To see this, let $\mathcal{F}(Z)$ denote the non-empty filters in $Z$ with the topology generated by $(\mathcal{F}'_z)_{z\in Z}$ where
\[\mathcal{F}'_z=\{F\in\mathcal{F}(Z):z\in F\}.\]
Further let $\mathcal{U}(Z)$ denote its subspace of ultrafilters, i.e. with the topology generated by $\mathcal{U}'_z=\mathcal{U}(Z)\cap\mathcal{F}'_z$, for $z\in Z$.  Note that $F^Z\in\mathcal{F}(Z)$, for any $F\in\mathcal{F}(S)$.  If $F\in\mathcal{F}^0$ then $F^Z=F\cap Z$ and $F=F^{Z<}=(F\cap Z)^<$, by \cite[Proposition 11.6]{Bice2022}, and hence $(\mathcal{F}_z)_{z\in Z}$ is a basis for $\mathcal{F}^0$.  On the other hand, if $F\in\mathcal{F}(Z)$ then $F^<\in\mathcal{F}^0$ and $F=F^{<Z}=F^<\cap Z$.  This shows that $F\mapsto F^Z$ and $F\mapsto F^<$ are mutually inverse homeomorphisms between $\mathcal{F}^0$ and $\mathcal{F}(Z)$ and hence restrict to homeomorphisms between $\mathcal{U}^0$ and $\mathcal{U}(Z)$.  As $Z$ is a generalised Boolean algebra, the classic Stone duality (see \cite{Stone1936} and \cite{Stone1938}) says that $z\mapsto\mathcal{U}'_z$ is a Boolean isomorphism from $Z$ onto the compact open subsets of $\mathcal{U}(Z)$.  Thus the same applies to $\mathcal{U}^0$, i.e. $z\mapsto\mathcal{U}_z$ is a Boolean isomorphism from $Z$ onto the compact open subsets of $\mathcal{U}^0$.

Now take $r\in S^\dagger$.  As $\mathcal{U}(S)$ is \'etale, $\mathcal{U}_r$ is homeomorphic to $\mathsf{r}[\mathcal{U}_r]=\mathcal{U}_{rr^{-1}}$, by \eqref{rFt}, which is compact, by the claim just proved.  Also $(\mathcal{U}_r)_{r\in S^\dagger}$ is a basis for $\mathcal{U}(S)$, thanks to \eqref{ReflexiveInterpolation}, which completes the proof.
\end{proof}

Let us denote the image of any $T\subseteq S$ under the ultrafilter representation by
\[\widehat{T}=\{\widehat{t}:t\in T\}.\]
If the ultrafilter representation $a\mapsto\widehat{a}$ is injective on $T$ then we call it \emph{faithful} on $T$.

\begin{thm}\label{FaithfulUltrafilters}
The ultrafilter representation $a\mapsto\widehat{a}$ is faithful on $S$.
\end{thm}

\begin{proof}
First we show that the ultrafilter representation is faithful on $\mathrm{ran}(\Phi)$.  To see this, take distinct $r,s\in\mathrm{ran}(\Phi)$, so either $r\nleq s$ or $s\nleq r$.  Assume that $r\nleq s$ and consider the subset of $Z$ defined by
\[F=\{y\setminus z:y,z\in Z,\ r<y\text{ and }rz=sz\}.\]

We claim $F$ is directed.  To see this, take $y,y',z,z'\in Z$ with $r<y,y'$, $rz=sz$ and $rz'=sz'$.  It follows that $ryy'=ry'=r$, i.e. $r<yy'$, and
\[r(z\vee z')=rz\vee rz'=sz\vee sz'=s(z\vee z').\]
Thus $(y\setminus z)(y'\setminus z')=yy'\setminus(z\vee z')\in F$, as required.

Next note that if we had $y,z\in Z$ with $r<y\leq z$ and $rz=sz$ then $r=rz=sz$, contradicting our assumption that $r\nleq s$.  Thus $y\nleq z$, for all $y,z\in Z$ with $r<y$ and $rz=sz$, and hence $0\notin F$.

Kuratowski-Zorn then yields $U\in\mathcal{U}^0$ extending $F$.  If $r\sim_Us$ then we would have $z\in U^Z$ with $rz=sz$, by \eqref{rhoUnit}.  For any $y>r$, it follows that $y\setminus z\in F\subseteq U$ and hence $0=(y\setminus z)z\in U^ZU\subseteq U$, a contradiction.  Thus $r\not\sim_Us$, which means
\[\widehat{r}(U)=[r,U]\neq[s,U]=\widehat{s}(U)\]
so $\widehat{r}\neq\widehat{s}$.  This shows that the ultrafilter representation is indeed faithful on $\mathrm{ran}(\Phi)$.

In particular, if $r\in\mathrm{ran}(\Phi)\setminus\{0\}$ then $\widehat{r}\neq\widehat{0}$.  But for arbitrary $s\in S\setminus\{0\}$, we can take $t'\in s^*$ and then $st'\in\mathrm{ran}(\Phi)\setminus\{0\}$ and hence $\widehat{s}\widehat{t'}=\widehat{st'}\neq\widehat{0}$ so again $\widehat{s}\neq\widehat{0}$.

We claim this implies that $s\in\mathrm{ran}(\Phi)$ whenever $\widehat{s}\in\widehat{\mathrm{ran}(\Phi)}$.  To see this, take $s\in S\setminus\mathrm{ran}(\Phi)$.  Then we can further take $y,z\in Z$ with $s=sy$ and $\Phi(s)=sz$, as $\Phi$ is quasi-Cartan, by \autoref{BistableQC}.  Letting $a=s(y\setminus z)$, we see that
\[\Phi(a)=\Phi(s)(y\setminus z)=sz(y\setminus z)=0\neq a.\]
Indeed, for the last $\neq$, take $t,t'\in S$ with $s<_{t'}t$ and note $a=0$ would imply $(y\setminus z)\perp t's$ so $t's=t'sy<z$, by \eqref{yzr}, and hence $\Phi(s)=sz=tt'sz=tt's=s$, contradicting $s\notin\mathrm{ran}(\Phi)$.  By what we just proved, $\widehat{a}\neq\widehat{0}=\widehat{\Phi(a)}=\Phi^{\mathsf{U}_{\langle S\rangle}}(\widehat{a})$, by \eqref{ContinuousPhi}, and hence $\widehat{a}\notin\widehat{\mathrm{ran}(\Phi)}$ (because $\Phi^{\mathsf{U}_{\langle S\rangle}}(\widehat{r})=\widehat{\Phi(r)}=\widehat{r}$, for any $r\in\mathrm{ran}(\Phi)$, again by \eqref{ContinuousPhi}).  As $\widehat{a}=\widehat{s}\widehat{y\setminus z}$ and $y\setminus z\in Z\subseteq\mathrm{ran}(\Phi)$, it follows that $\widehat{s}\notin\widehat{\mathrm{ran}(\Phi)}$ too, which finishes the proof of the claim.

Now take any $a,b\in S$ and assume that $\widehat{a}=\widehat{b}$.  Taking $s,s'\in S$ with $a<_{s'}s$, it follows that $\widehat{bs'}=\widehat{b}\widehat{s'}=\widehat{a}\widehat{s'}=\widehat{as'}\in\widehat{\mathrm{ran}(\Phi)}$ and hence $bs'\in\mathrm{ran}(\Phi)$, by what we just proved.  As we already showed that the ultrafilter representation is faithful on $\mathrm{ran}(\Phi)$, it follows that $as'=bs'$ and hence $a=as's=bs's$.  Likewise, taking $t,t'\in S$ with $b<_{t'}t$, we see that $b=at't$.  Thus $a\leq b\leq a$ and hence $a=b$, which shows that the ultrafilter representation is indeed faithful on the whole of $S$.
\end{proof}

Finally we can now show that the ultrafilter bundle is ample and that its sections correspond precisely to the original elements of the Steinberg semigroup.

\begin{thm}\label{AmpleUltrafilters}
The ultrafilter bundle $\mathsf{U}_{\langle S\rangle}:\mathcal{U}[S]\twoheadrightarrow\mathcal{U}(S)$ is ample.  Moreover,
\[\widehat{S}=\mathcal{S}_\mathsf{c}(\mathsf{U}_{\langle S\rangle})\qquad\text{and}\qquad\widehat{Z}=\mathcal{Z}_\mathsf{c}(\mathsf{U}_{\langle S\rangle}).\]
\end{thm}

\begin{proof}
We know that $\mathsf{U}_{\langle S\rangle}$ is an \'etale zero category bundle, by \autoref{piUEtaleZeroCategoryBundle}, with $\mathsf{U}_{\langle S\rangle}(\mathcal{U}[S]^\times)=\mathcal{U}(S)$, by \eqref{FStimes}.  Moreover, the base groupoid $\mathcal{U}(S)$ is Hausdorff, by \autoref{HausdorffUltrafilters}, and compactly based, by \autoref{RSbasis}.  Thus $\mathsf{U}_{\langle S\rangle}$ is ample.

We also saw in the proof of \autoref{RSbasis} that $z\mapsto\mathcal{U}_z$ is a Boolean isomorphism from $Z$ onto the compact open subsets of $\mathcal{U}^0$.  To show that $\widehat{Z}=\mathcal{Z}_\mathsf{c}(\mathsf{U}_{\langle S\rangle})$, it thus suffices to show that $\widehat{z}$ is the characteristic function of $\mathcal{U}_z$, for each $z\in Z$.  For this, just observe that $\mathrm{supp}(\widehat{z})\subseteq\mathcal{U}_z$, by \eqref{SliceSupport}, as $z<z$, and $\widehat{z}(U)=[z,U]\in\mathcal{U}[S]^0$, for all $U\in\mathcal{U}_z$, by \eqref{FS0}.  So $z\mapsto\widehat{z}$ is a Boolean isomorphism from $Z$ onto $\mathcal{Z}_\mathsf{c}(\mathsf{U}_{\langle S\rangle})$.

It follows that $a\mapsto\widehat{a}$ preserves orthogonal suprema, i.e. for all $a,b\in S$,
\begin{equation}\label{Svee}
a\perp b\qquad\Rightarrow\qquad\widehat{a\vee b}=\widehat{a}\vee\widehat{b}.
\end{equation}
Indeed, if $a\perp b$ then \autoref{perpSymmetric} yields $y,z\in Z$ with $a=ay$, $b=bz$ and $yz=0$.  It follows that $(a\vee b)y=ay\vee by=a\vee bzy=a\vee0=a$.  Likewise, $(a\vee b)z=b$ and hence $(a\vee b)(y\vee z)=(a\vee b)y\vee(a\vee b)z=a\vee b$ so
\[\widehat{a\vee b}=\widehat{(a\vee b)}\widehat{(y\vee z)}=\widehat{(a\vee b)}(\widehat{y}\vee\widehat{z})=\widehat{(a\vee b)}\widehat{y}\vee\widehat{(a\vee b)}\widehat{z}=\widehat{a}\vee\widehat{b}.\]

Next we want to show that $\widehat{S}\subseteq\mathcal{S}_\mathsf{c}(\mathsf{U}_{\langle S\rangle})$.  Note that every $s\in S$ is dominated by some $a\in S^\dagger$, which implies $\mathrm{supp}(\widehat{s})\subseteq\mathcal{U}_a$, again by \eqref{SliceSupport}.  Thus it suffices to show that $\mathcal{U}_a$ is compact, for all $a\in S^\dagger$.  But $\mathcal{U}(S)$ is \'etale so $\mathcal{U}_a$, which is an open slice by \cite[Proposition 7.3]{Bice2022}, is homeomorphic to $\mathsf{r}[\mathcal{U}_a]=\mathcal{U}_{aa^{-1}}$, by \eqref{rFt}, which is compact because $aa^{-1}\in Z$.

Now we just have to show that $\mathcal{S}_\mathsf{c}(\mathsf{U}_{\langle S\rangle})\subseteq\widehat{S}$.  Accordingly, take any section $s\in\mathcal{S}_\mathsf{c}(\mathsf{U}_{\langle S\rangle})$.  As $\mathrm{supp}(s)$ is a compact slice, it is contained in some open slice $O$, by \cite[Proposition 6.3]{BiceStarling2018}.  For every $U\in\mathrm{supp}(s)$, we have $a\in S$ such that $s(U)=[a,U]=\widehat{a}(U)$.  As $\mathsf{U}_{\langle S\rangle}$ is ample, $s$ and $\widehat{a}$ agree on some compact open slice $K\subseteq O$ containing the point $U$.  Then $\mathsf{r}[K]=\mathcal{U}_z$, for some $z\in Z$, so we can replace $a$ with $za$ if necessary to ensure that $\widehat{a}=\widehat{z}s$.  As $\mathrm{supp}(s)$ is compact, we can cover it with finitely many such $K$, i.e. we have compact open $K_1,\ldots,K_n\subseteq O$ covering $\mathrm{supp}(s)$ and $a_1,\ldots,a_n\in S$ with $\widehat{a_m}=\widehat{z_m}s$ where $\mathsf{r}[K_m]=\mathcal{U}_{z_m}$, for all $m\leq n$.  If necessary, we may replace each $K_m$ with $K_m\setminus\bigcup_{l<m}K_l$ (and again replace each $a_m$ with $z_ma_m$, where $\mathsf{r}[K_m]=\mathcal{U}_{z_m}$) to ensure that $K_1,\ldots,K_n$ are disjoint.  Taking $y_1,\ldots,y_n\in Z$ with $U_{y_m}=\mathsf{s}[K_m]$, it follows that $\widehat{y_l}\widehat{y_m}=\widehat{0}=\widehat{z_l}\widehat{z_m}$, whenever $l\neq m$, and hence $y_ly_m=0=z_lz_m$, by \autoref{FaithfulUltrafilters}.  Likewise $z_ma=a_m=ay_m$, for all $m\leq n$, and hence $a_l\perp a_m$, for all $l\neq m$.  Thus we may let $a=\bigvee_{m=1}^na_m\in S$, by \eqref{Orthosuprema} and \eqref{veeperp}.  Further let $z=\bigvee_{m=1}^nz_m$ so $\mathcal{U}_z=\bigcup_{m=1}^n\mathsf{r}[K_m]$ covers $\mathsf{r}(\mathrm{supp}(s))$ and hence $s=\widehat{z}s$.  As $a_m=z_ma_m$, for all $m\leq n$, it also follows that $a=za$ so \eqref{Svee} yields $\widehat{a}=\bigvee_{m=1}^n\widehat{a_m}=\bigvee_{m=1}^ns\widehat{z_m}=s\bigvee_{m=1}^n\widehat{z_m}=s\widehat{z}=s$, showing that $s\in\mathcal{S}_\mathsf{c}(\mathsf{U}_{\langle S\rangle})$.  Thus we do indeed have $\widehat{S}=\mathcal{S}_\mathsf{c}(\mathsf{U}_{\langle S\rangle})$, finishing the proof.
\end{proof}

The following corollary summarises what we have proved so far.

\begin{cor}\label{SteinbergBundleCharacterisation}
Steinberg semigroups are, up to isomorphism, precisely the well-structured semigroups $(\mathcal{S}_\mathsf{c}(\rho),\mathcal{Z}_\mathsf{c}(\rho),\Phi_\mathsf{c}^\rho)$ arising from ample category bundles $\rho$.
\end{cor}

\begin{proof}
If $\rho$ is an ample category bundle then $(\mathcal{S}_\mathsf{c}(\rho),\mathcal{Z}_\mathsf{c}(\rho),\Phi_\mathsf{c}^\rho)$ is a Steinberg semigroup, by \autoref{AmpleBundleSteinbergSemigroup}.  Conversely, if $\langle S\rangle=(S,Z,\Phi)$ is a Steinberg semigroup then its ultrafilter bundle $\rho=\mathsf{U}_{\langle S\rangle}$ is ample and the ultrafilter representation $a\mapsto\widehat{a}$ is an isomorphism of $\langle S\rangle$ onto the resulting Steinberg semigroup $(\mathcal{S}_\mathsf{c}(\rho),\mathcal{Z}_\mathsf{c}(\rho),\Phi_\mathsf{c}^\rho)$, by \autoref{ContinuousSection} and \autoref{AmpleUltrafilters}.
\end{proof}

\section{Lawson-Steinberg Duality}\label{LawsonSteinbergDuality}

To beef up \autoref{SteinbergBundleCharacterisation} to an equivalence of categories we must also consider appropriate morphisms.  We start with the morphisms on the algebraic side.

\begin{dfn}\label{SteinbergMorphisms}
Assume $\langle S\rangle=(S,Z,\Phi)$ and $\langle S'\rangle=(S',Z',\Phi')$ are Steinberg semigroups.  A \emph{Steinberg morphism} from $\langle S\rangle$ to $\langle S'\rangle$ is a semigroup homomorphism $\pi:S\rightarrow S'$ such that $\pi(0)=0$, $\pi[Z]\subseteq Z'$ and, for all orthogonal $a,b\in S$,
\[\pi(a\vee b)=\pi(a)\vee\pi(b)\qquad\text{and}\qquad\pi(\Phi(a))=\Phi'(\pi(a)).\]
\end{dfn}

On the topological side, we first consider just the base groupoid of our bundles.  Recall that a map $\phi:C\rightarrow C'$ between categories $C$ and $C'$ is a \emph{functor} if it preserves units and the product, i.e. $\phi[C^0]\subseteq C'^0$ and $\phi(ab)=\phi(a)\phi(b)$ whenever $(a,b)\in C^2$.  Such a functor $\phi$ is \emph{star-bijective} if it maps all stars $Ce$ and $eC$, for $e\in C^0$, bijectively onto the respective stars $C'\phi(e)$ and $\phi(e)C'$.

\begin{dfn}
Assume $G$ and $G'$ are \'etale groupoids.  We call $\phi:G\rightarrow G'$ an \emph{\'etale morphism} from $G$ to $G'$ if $\phi$ is a continuous star-bijective functor.
\end{dfn}

\begin{thm}\label{UltraFunctoriality}
Any Steinberg morphism $\pi:S\rightarrow S'$ from one Steinberg semigroup $(S,Z,\Phi)$ to another $(S',Z',\Phi')$ defines an \'etale morphism $\underline{\pi}$ in the opposite direction, from an open subgroupoid of $\mathcal{U}(S')$ to $\mathcal{U}(S)$, given by
\begin{equation}\label{underlinepi}
\underline{\pi}(U)=\pi^{-1}[U]^<\quad\text{when}\quad\pi^{-1}[U]\neq\emptyset.
\end{equation}
\end{thm}

\begin{proof}
First we must show $\underline{\pi}$ is well-defined, i.e. $T=\pi^{-1}[U]^<$ is an ultrafilter when $U\in\mathcal{U}(S')$ and $\pi^{-1}[U]\neq\emptyset$.  First note $T=T^<$ is immediate from \eqref{ReflexiveInterpolation}, while $0\notin T$ is immediate from $0\notin U$.  Also note that $T^*\subseteq\pi^{-1}[U^*]^<$.  Indeed, if $t'\in T^*$ then we have $a\in\pi^{-1}[U]$ and $s,s',t\in S$ with $a<_{s'}s<_{t'}t$.  Note that $\pi(a)<_{\pi(s')}\pi(s)$, as $\pi$ is a Steinberg morphism, and hence $\pi(s')\in\pi(a)^*\subseteq U^*$.  Replacing $s'$ with $s'ss'$ if necessary, it also follows that $s'<_tt'$, showing that $t'\in\pi^{-1}[U^*]^<$.

Now take $a,b\in T$ and $b'\in T^*$, so we have $c,d\in\pi^{-1}[U]$ and $d'\in\pi^{-1}[U^*]$ with $c<a$, $d<b$ and $d'<b'$.  Then $\Phi(cd')d<\Phi(ab')b$, by \cite[Proposition 4.7]{Bice2022} and \eqref{PhiCompatibility}.  As $\pi$ is a Steinberg morphism,
\[\pi(\Phi(cd')d)=\Phi'(\pi(c)\pi(d'))\pi(d)\in\Phi'[UU^*]U\subseteq U,\]
and hence $\Phi(ab')b\in\pi^{-1}[U]^<$.  This shows that $\Phi[TT^*]T\subseteq T$ and hence $T$ is also a filter, by \autoref{FilterPhi}.

To see that $T$ is an ultrafilter, take orthogonal $a,b\in\mathsf{R}(S)$ with $a\vee b\in T$, so we have $c<a\vee b$ with $\pi(c)\in U$.  As $\pi$ is a Steinberg morphism, $\pi(a)\perp\pi(b)$ and
\[\pi(c)<\pi(a\vee b)=\pi(a)\vee\pi(b).\]
Thus $\pi(a)\vee\pi(b)\in U$ so $\pi(a)\in U$ or $\pi(b)\in U$, by \eqref{Ultra+}, and hence $a\in T$ or $b\in T$, as $a$ and $b$ are $Z$-invertible.  This shows that $T$ is an ultrafilter, again by \eqref{Ultra+}.

Next note that, for any $T\in\mathcal{U}(S)$ and $U\in\mathcal{U}(S')$,
\[\underline{\pi}(U)=T\qquad\Leftrightarrow\qquad\pi[T]\subseteq U.\]
Indeed, if $T=\underline{\pi}(U)=\pi^{-1}[U]^<$ then $\pi[T]=\pi[\pi^{-1}[U]^<]\subseteq U^<\subseteq U$.  Conversely, if $\pi[T]\subseteq U$ then $T=T^<\subseteq\pi^{-1}[U]^<=\underline{\pi}(U)$ and hence $T=\underline{\pi}(U)$, as $T$ and $\underline{\pi}(U)$ are ultrafilters.

Now to see that $\underline{\pi}$ is a functor, say $T=\underline{\pi}(U)$ and $V=\underline{\pi}(W)$, i.e. $\pi[T]\subseteq U$ and $\pi[V]\subseteq W$ and hence $\pi[(TV)^<]\subseteq(UW)^<$.  Thus if $0\notin UW$ then $0\notin TV$, i.e. if $U\cdot W$ is defined then so is $T\cdot V$, in which case $\underline{\pi}(U)\cdot\underline{\pi}(W)=T\cdot V=\underline{\pi}[U\cdot W]$.

To show that $\underline{\pi}$ is star-bijective, take $U\in\mathcal{U}(S')$ and $T\in\mathcal{U}(S)$ with $\underline{\pi}(U)=\mathsf{s}(T)$.  As $\mathsf{s}(T)$ is a unit, it contains some $r\in\mathrm{ran}(\Phi)$ and hence $U$ contains $\pi(r)\in\mathrm{ran}(\Phi')$, so $U$ is also a unit (see \cite[Proposition 7.2]{Bice2022}).  This shows that $\underline\pi$ is star-injective (see \cite[Proposition 3.4]{Bice2022}).  Also, for any $t\in T$, \cite[Proposition 6.6 and Theorem 11.7]{Bice2022} yield an ultrafilter $V=(\pi(t)U)^<\in\mathcal{U}(S')$ with $\mathsf{s}(V)=U$.  By \cite[Proposition 6.5]{Bice2022}, $(t\mathsf{s}(T))^<=T$ and hence
\[\pi[T]=\pi[(t\mathsf{s}(T))^<]\subseteq(\pi(t)U)^<=V.\]
Thus $T=\underline{\pi}(V)$, showing that $\underline{\pi}$ is also star-surjective.

To show that $\underline\pi$ is also proper and continuous, it suffices to show that the preimage of every set in the compact open basis $(\mathcal{U}_r)_{r\in S^\dagger}$ (see \autoref{RSbasis}) is again compact and open.  To see this, note that
\[\underline\pi^{-1}[\mathcal{U}_r]=\{U\in\mathcal{U}(S'):r\in\underline\pi(U)\}=\{U\in\mathcal{U}(S'):\pi(r)\in U\}=\mathcal{U}'_{\pi(r)}.\]
As $r$ is $Z$-invertible and $\pi$ is a Steinberg morphism, $\pi(r)$ is also $Z'$-invertible and hence $\mathcal{U}'_{\pi(r)}$ is compact open, by \autoref{RSbasis}.
\end{proof}

We can immediately turn Steinberg morphisms into a category
\[\mathbf{SS}=\{(\langle S'\rangle,\pi,\langle S\rangle):\pi\text{ is a Steinberg morphism from }\langle S\rangle\text{ to }\langle S'\rangle\},\]
with the product given by composition when the domain and codomain match, i.e.
\[(\langle S''\rangle,\pi',\langle S'\rangle)(\langle S'\rangle,\pi,\langle S\rangle)=(\langle S''\rangle,\pi'\circ\pi,\langle S\rangle).\]
Likewise, we can turn \'etale morphisms into a category although, as in the above result, this time we want the hom-sets to include partial maps.  Accordingly let
\[\mathbf{\acute{E}G}=\{(G',\pi,G):\pi\text{ is an \'etale morphism from an open subgroupoid }O\text{ of }G\text{ to }G'\},\]
with the product again given by composition, i.e.
\[(G'',\pi',G')(G',\pi,G)=(G'',\pi'\circ\pi,G).\]

\begin{cor}\label{SSEG}
We have a contravariant functor from $\mathbf{SS}$ to $\mathbf{\acute{E}G}$ given by
\begin{equation}\label{SStoEG}
(\langle S'\rangle,\pi,\langle S\rangle)\mapsto(\mathcal{U}(S),\underline\pi,\mathcal{U}(S')).
\end{equation}
\end{cor}

\begin{proof}
By \autoref{UltraFunctoriality}, \eqref{SStoEG} is indeed a well-defined function from $\mathbf{SS}$ to $\mathbf{\acute{E}G}$.  All we have to prove is that $\pi\mapsto\underline\pi$ respects composition.  Accordingly, take Steinberg morphisms $\pi:S\rightarrow S'$ and $\pi':S'\rightarrow S''$ between Steinberg semigroups $(S,Z,\Phi)$, $(S',Z',\Phi')$ and $(S'',Z'',\Phi'')$.  For any $U\in\mathcal{U}(S'')$, it follows from \eqref{ReflexiveInterpolation} that
\[\underline{\pi'\circ\pi}(U)=(\pi'\circ\pi)^{-1}[U]^<=\pi^{-1}[\pi'^{-1}[U]]^{<<}\subseteq\pi^{-1}[\pi'^{-1}[U]^<]^<=\underline\pi\circ\underline\pi'(U).\]
But inclusion on ultrafilters implies equality so this yields $\underline{\pi'\circ\pi}=\underline\pi\circ\underline\pi'$.
\end{proof}

Next we want go the other way and turn \'etale morphisms into Steinberg morphisms.  First note that, from any ample category bundle $\rho:C\twoheadrightarrow G'$ and \'etale morphism $\phi:G\rightarrow G'$ between Hausdorff ample groupoids $G$ and $G'$, we obtain another topological category from the pullback
\[\phi^\rho C=\{(g,c)\in G\times C:\phi(g)=\rho(c)\},\]
considered as a topological subspace of $G\times B$ under the product
\[(f,b)(g,c)=(fg,bc),\text{ when }(f,g)\in G^2.\]
Indeed, continuity of the product follows from the continuity of the products in $G$ and $C$ and likewise for the source $\mathsf{s}(g,c)=(\mathsf{s}(g),\mathsf{s}(c))$ and range $\mathsf{r}(g,c)=(\mathsf{r}(g),\mathsf{r}(c))$.  As $\rho$ is locally injective, so is the \emph{pullback bundle} $\rho_\phi:\phi^\rho C\twoheadrightarrow G$ where $\rho_\phi(g,c)=g$.  Thus $\rho_\phi$ is also an ample category bundle.

The pullback-section map $\phi^*:\mathcal{A}(\rho)\rightarrow\mathcal{A}(\rho_\phi)$ is then given by
\[\phi^*(a)(g)=(g,a(\phi(g))).\]
As $\phi$ is proper and continuous, $\phi^*$ restricts to a map from $\mathcal{C}_\mathsf{c}(\rho)$ to $\mathcal{C}_\mathsf{c}(\rho_\phi)$.  Let us also denote the Steinberg semigroup from \autoref{AmpleBundleSteinbergSemigroup} by 
\[\mathsf{S}\rho=(\mathcal{S}_\mathsf{c}(\rho),\mathcal{Z}_\mathsf{c}(\rho),\Phi^\rho_\mathsf{c}).\]

\begin{prp}\label{phiSteinberg}
If $\rho:C\twoheadrightarrow G'$ is a ample category bundle and $\phi:G\rightarrow G'$ is an \'etale morphism between Hausdorff ample groupoids $G$ and $G'$ then the pullback-section map $\phi^*$ on $\mathcal{S}_\mathsf{c}(\rho)$ is a Steinberg morphism from $\mathsf{S}\rho$ to $\mathsf{S}\rho_\phi$.
\end{prp}

\begin{proof}
Take any $a,b\in\mathcal{S}_\mathsf{c}(\rho)$.  Then
\[\mathrm{supp}(\phi^*(ab))=\phi^{-1}[\mathrm{supp}(ab)]\subseteq\phi^{-1}[\mathrm{supp}(a)\mathrm{supp}(b)]\subseteq\phi^{-1}[\mathrm{supp}(a)]\phi^{-1}[\mathrm{supp}(b)],\]
as $\phi$ is star-surjective (see \cite[Proposition 3.5]{Bice2022}).  As $\phi$ is star-injective, both $\phi^{-1}[\mathrm{supp}(a)]$ and $\phi^{-1}[\mathrm{supp}(b)]$ are slices (see \cite[Proposition 3.4]{Bice2022}) and hence so is $\mathrm{supp}(\phi^*(ab))$.  Now, for any $g\in\phi^{-1}[\mathrm{supp}(a)]$ and $h\in\phi^{-1}[\mathrm{supp}(b)]$,
\[\phi^*(ab)(gh)=(gh,ab(\phi(gh)))=(g,a(\phi(g)))(h,b(\phi(h)))=\phi^*(a)(g)\phi^*(b)(h).\]
Thus $\phi^*(ab)=\phi^*(a)\phi^*(b)$, showing $\phi^*$ restricts to a homomorphism from $\mathcal{S}_\mathsf{c}(\rho)$ to $\mathcal{S}_\mathsf{c}(\rho_\phi)$.  We also immediately see that $\phi^*$ maps the zero section of $\rho$ to that of $\rho_\phi$.

Now, for any orthogonal $a,b\in\mathcal{S}_\mathsf{c}(\rho)$, note that
\[\mathrm{supp}(\phi^*(a\vee b))=\phi^{-1}[\mathrm{supp}(a\vee b)]=\phi^{-1}[\mathrm{supp}(a)]\cup\phi^{-1}[\mathrm{supp}(b)].\]
As $a\vee b$ agrees with $a$ on $\mathrm{supp}(a)$ and with $b$ on $\mathrm{supp}(b)$, it follows that $\phi^*(a\vee b)$ agrees with $\phi^*(a)$ on $\phi^{-1}[\mathrm{supp}(a)]$ and with $\phi^*(b)$ on $\phi^{-1}[\mathrm{supp}(b)]$, showing that
\[\phi^*(a\vee b)=\phi^*(a)\vee\phi^*(b).\]

As $\phi$ is star-injective, $\phi^{-1}[G'^0]\subseteq G^0$ (see \cite[Proposition 3.4]{Bice2022}).  This means
\[(\phi^\rho C)^0=\{(g,c)\in\phi^\rho C:c\in C^0\}\]
and hence if $a\in\mathcal{S}_\mathsf{c}(\rho)$ only takes unit values on its support then the same is true of $\phi^*(a)$, i.e. $\phi^*[\mathcal{Z}_\mathsf{c}(\rho)]\subseteq\mathcal{Z}_\mathsf{c}(\rho_\phi)$.  It also means that if $\mathrm{supp}(a)\subseteq G'^0$ then $\mathrm{supp}(\phi^*(a))\subseteq G^0$, in particular $\mathrm{supp}(\phi^*(\Phi^\rho(a)))\subseteq G^0$.  Also, for all $g\in G^0$, we see that $\phi(g)\in G'^0$, as $\phi$ is a functor, and hence
\[\phi^*(\Phi^\rho(a))(g)=(g,\Phi^\rho(a)(\phi(g)))=(g,a(\phi(g)))=\phi^*(a)(g)=\Phi^{\rho_\phi}(\phi^*(a))(g).\]
This shows that $\phi^*(\Phi^\rho(a))=\Phi^{\rho_\phi}(\phi^*(a))$ and hence $\phi^*$ is indeed a Steinberg morphism from $\mathsf{S}\rho$ to $\mathsf{S}\rho_\phi$.
\end{proof}

We can also obtain Steinberg morphisms from functors on the total categories.

\begin{dfn}
Let $\rho:C\rightarrow G$ and $\rho':C'\rightarrow G'$ be ample category bundles where $G$ is an open subgroupoid of $G'$.  A \emph{bundle morphism} from $\rho$ to $\rho'$ is a continuous functor $\beta:C\rightarrow C'$ such that $\rho'\circ\beta=\rho$ and $\beta(0_g)=0_g$, for all $g\in G$.
\end{dfn}

If $\beta:C\rightarrow C'$ is a bundle morphism from $\rho:C\rightarrow G$ to $\rho':C'\rightarrow G'$ then, as $\rho'\circ\beta=\rho$, we obtain a pushforward-section map $\beta_*:\mathcal{A}(\rho)\rightarrow\mathcal{A}(\rho')$ by
\begin{equation}\label{BetaDef}
\beta_*(a)(g)=\begin{cases}\beta(a(g))&\text{if }g\in G\\0_g&\text{if }g\in G'\setminus G.\end{cases}
\end{equation}
Note that $\mathrm{supp}(\beta_*(a))\subseteq\mathrm{supp}(a)$, for all $a\in\mathcal{A}(\rho)$.  As $\beta$ is continuous, it follows that $\beta_*$ restricts to a map from $\mathcal{C}_\mathsf{c}(\rho)$ to $\mathcal{C}_\mathsf{c}(\rho')$ and from $\mathcal{S}_\mathsf{c}(\rho)$ to $\mathcal{S}_\mathsf{c}(\rho')$.

\begin{prp}\label{betaSteinberg}
If $\beta:C\rightarrow C'$ is a bundle morphism from $\rho:C\rightarrow G$ to $\rho':C'\rightarrow G'$ then $\beta_*$ on $\mathcal{S}_\mathsf{c}(\rho)$ is a Steinberg morphism from $\mathsf{S}\rho$ to $\mathsf{S}\rho'$.
\end{prp}

\begin{proof}
Take any $a,b\in\mathcal{S}_\mathsf{c}(\rho)$.  As $\beta$ is a functor, for all $g\in\mathrm{supp}(a)$ and $h\in\mathrm{supp}(b)$,
\[\beta_*(ab)(gh)=\beta(ab(gh))=\beta(a(g)b(h))=\beta(a(g))\beta(b(h))=\beta^*(a)(g)\beta^*(b)(h).\]
Thus $\beta_*(ab)=\beta^*(a)\beta^*(b)$, as $\mathrm{supp}(\beta_*(ab))\subseteq\mathrm{supp}(ab)\subseteq\mathrm{supp}(a)\mathrm{supp}(b)$ and
\[\mathrm{supp}(\beta^*(a)\beta^*(b))\subseteq\mathrm{supp}(\beta^*(a))\mathrm{supp}(\beta^*(b))\subseteq\mathrm{supp}(a)\mathrm{supp}(b).\]
This shows $\beta_*$ restricts to a homomorphism from $\mathcal{S}_\mathsf{c}(\rho)$ to $\mathcal{S}_\mathsf{c}(\rho')$.  For any orthogonal $a,b\in\mathcal{S}_\mathsf{c}(\rho)$, we also immediately see that $\beta_*(a\vee b)=\beta_*(a)\vee\beta_*(b)$.

As $\beta(0_g)=0_g$, for all $g\in G$, $\beta_*$ maps the zero section of $\rho$ to that of $\rho'$ and $\beta_*(\Phi^\rho(a))=\Phi^{\rho'}\!(\beta_*(a))$, for all $a\in\mathcal{A}(\rho)$.  As $\beta$ is a functor, $\beta[C^0]\subseteq C'^0$ so $\beta_*[\mathcal{Z}_\mathsf{c}(\rho)]\subseteq\beta_*[\mathcal{Z}_\mathsf{c}(\rho')]$, i.e. $\beta_*$ is a Steinberg morphism from $\mathsf{S}\rho$ to $\mathsf{S}\rho'$.
\end{proof}

Now we can form another category by combining \'etale morphisms and bundle morphisms between ample category bundles.

\begin{dfn}
Let $\mathbf{ACB}$ denote quadruples $(\rho',\beta,\phi,\rho)$ where
\begin{enumerate}
\item $\rho:C\twoheadrightarrow G$ and $\rho':C'\twoheadrightarrow G'$ are ample category bundles.
\item $\phi$ is an \'etale morphism from an open subgroupoid of $G'$ to $G$.
\item $\beta$ is a bundle morphism from $\rho_\phi$ to $\rho'$.
\end{enumerate}
\end{dfn}

If $(\rho',\beta,\phi,\rho)\in\mathbf{ACB}$, we say the pair $(\beta,\phi)$ is a \emph{Pierce morphism} from $\rho$ to $\rho'$.

\begin{prp}\label{ACBcategory}
$\mathbf{ACB}$ forms a category under the product
\[(\rho'',\beta',\phi',\rho')(\rho',\beta,\phi,\rho)=(\rho'',\beta'\bullet\beta,\phi\circ\phi',\rho)\]
where $\beta'\bullet\beta=\beta'\midscript{\phi'}\beta$ is the bundle morphism from $\rho_{\phi\circ\phi'}$ to $\rho''$ defined by
\[\beta'\bullet\beta(g'',b)=\beta'(g'',\beta(\phi'(g''),b)).\]
\end{prp}

\begin{proof}
For associativity, take $(\rho''',\beta'',\phi'',\rho''),(\rho'',\beta',\phi',\rho'),(\rho',\beta,\phi,\rho)\in\mathbf{ACB}$ and note that, for all $g'''\in\mathrm{dom}(\phi\circ\phi'\circ\phi'')=\phi''^{-1}[\phi'^{-1}[\phi^{-1}[G]]]$ and $c\in C$,
\begin{align*}
(\beta''\bullet(\beta'\bullet\beta))(g''',c)&=\beta''(g''',(\beta'\bullet\beta)(\phi''(g'''),c))\\
&=\beta''(g''',\beta'(\phi''(g'''),\beta(\phi'(\phi''(g''')),c)))\\
&=(\beta''\bullet\beta')(g''',\beta(\phi'\circ\phi''(g'''),c))\\
&=((\beta''\bullet\beta')\bullet\beta)(g''',c).
\end{align*}
This shows that $\beta''\bullet(\beta'\bullet\beta)=(\beta''\bullet\beta')\bullet\beta$, as required.

We immediately see that every $(\rho',\beta,\phi,\rho)\in\mathbf{ACB}$ has a source unit $(\rho,\mathsf{p}_C,\mathrm{id}_G,\rho)$ and range unit $(\rho',\mathsf{p}_{C'},\mathrm{id}_{G'},\rho')$, where $\mathrm{id}_G$ is the identity on $G$ and $\mathsf{p}_C(g,c)=c$, for all $(g,c)\in\mathrm{id}_G^\rho C$.  Thus $\mathbf{ACB}$ forms a category under the given product.
\end{proof}

Combining a couple of the results above yields a functor from $\mathbf{ACB}$ to $\mathbf{SS}$.

\begin{thm}\label{Afunctor}
We have a functor $\mathsf{S}:\mathbf{ACB}\rightarrow\mathbf{SS}$ given by
\[\mathsf{S}(\rho',\beta,\phi,\rho)=(\mathsf{S}\rho',\beta_*\circ\phi^*,\mathsf{S}\rho)\]
\end{thm}

\begin{proof}
By \autoref{phiSteinberg}, $\phi^*$ on $\mathcal{S}_\mathsf{c}(\rho)$ is a Steinberg morphism from $\mathsf{S}\rho$ to $\mathsf{S}\rho_\phi$.  By \autoref{betaSteinberg}, $\beta_*$ on $\mathcal{S}_\mathsf{c}(\rho_\phi)$ is a Steinberg morphism from $\mathsf{S}\rho_\phi$ to $\mathsf{S}\rho'$.   Thus their composition $\beta_*\circ\phi^*$ is a Steinberg morphism from $\mathsf{S}\rho$ to $\mathsf{S}\rho'$.

To see that $\mathsf{S}$ preserves products, take $(\rho'',\beta',\phi',\rho'),(\rho',\beta,\phi,\rho)\in\mathbf{ACB}$ so
\[\mathsf{S}((\rho'',\beta',\phi',\rho')(\rho',\beta,\phi,\rho))=\mathsf{S}(\rho'',\beta'\bullet\beta,\phi\circ\phi',\rho)=(\mathsf{S}\rho'',(\beta'\bullet\beta)_*\circ(\phi\circ\phi')^*,\mathsf{S}\rho).\]
For any $a\in\mathcal{S}_\mathsf{c}(\rho)$ and $g''\in\mathrm{dom}(\phi\circ\phi')=\phi'^{-1}[\phi^{-1}[\mathrm{ran}(\rho)]]$,
\begin{align*}
((\beta'\bullet\beta)_*\circ(\phi\circ\phi')^*)(a)(g'')&=(\beta'\bullet\beta)_*((\phi\circ\phi')^*(a))(g'')\\
&=(\beta'\bullet\beta)(g'',a(\phi\circ\phi'(g'')))\\
&=\beta'(g'',\beta(\phi'(g''),a(\phi(\phi'(g'')))))\\
&=\beta'(g'',\beta(\phi^*(a)(\phi'(g''))))\\
&=\beta'(g'',\beta_*\circ\phi^*(a)(\phi'(g'')))\\
&=\beta'(\phi'^*\circ\beta_*\circ\phi^*(a)(g''))\\
&=(\beta'_*\circ\phi'^*\circ\beta_*\circ\phi^*(a))(g'').
\end{align*}
For any $g''\in\mathrm{ran}(\rho'')\setminus\mathrm{dom}(\phi\circ\phi')$, we also see that
\[((\beta'\bullet\beta)_*\circ(\phi\circ\phi')^*)(a)(g'')=0_{g''}=(\beta'_*\circ\phi'^*\circ\beta_*\circ\phi^*(a))(g'').\]
This shows that $(\beta'\bullet\beta)_*\circ(\phi\circ\phi')^*=\beta'_*\circ\phi'^*\circ\beta_*\circ\phi^*$ and hence
\begin{align*}
\mathsf{S}((\rho'',\beta',\phi',\rho')(\rho',\beta,\phi,\rho))&=(\mathsf{S}\rho'',\beta'_*\circ\phi'^*\circ\beta_*\circ\phi^*,\mathsf{S}\rho)\\
&=(\mathsf{S}\rho'',\beta'^*\circ\phi'^*,\mathsf{S}\rho')(\mathsf{S}\rho',\beta_*\circ\phi^*,\mathsf{S}\rho)\\
&=\mathsf{S}(\rho'',\beta',\phi',\rho')\mathsf{S}(\rho',\beta,\phi,\rho).
\end{align*}
Also $\mathsf{S}$ takes units in $\mathbf{ACB}$ to units in $\mathbf{SS}$ and hence $\mathsf{S}$ is a functor.
\end{proof}

We call $\mathsf{S}$ the \emph{section functor}.  Next we want to show that ultrafilters yield an inverse `bundle functor', at least modulo natural isomorphims.

Recall that that the ultrafilter bundle of any Steinberg semigroup $\langle S\rangle=(S,Z,\Phi)$ is denoted by by $\mathsf{U}_{\langle S\rangle}:\mathcal{U}[S]\rightarrow\mathcal{U}(S)$.

\begin{prp}\label{PullbackBundleMorphism}
Any Steinberg morphism $\pi:S\rightarrow S'$ between Steinberg semigroups $\langle S\rangle=(S,Z,\Phi)$ and $\langle S'\rangle=(S',Z',\Phi')$ defines a bundle morphism $\overline\pi$ on $\underline\pi^{\mathsf{U}_{\langle S\rangle}}\mathcal{U}[S]$ from the pullback bundle $\mathsf{U}_{\langle S\rangle\underline\pi}$ $($defined by $\underline\pi$ from \eqref{underlinepi}$)$ to $\mathsf{U}_{\langle S'\rangle}$ by
\[\overline\pi(U,[a,\underline\pi(U)])=[\pi(a),U].\]
\end{prp}

\begin{proof}
Take any $U\in\mathrm{dom}(\underline\pi)$ so $U\in\mathcal{U}(S')$ and $\pi^{-1}[U]\neq\emptyset$.  For any $a,b\in S$ with $a\sim_{\underline\pi(U)}b$, we have $s\in\underline\pi(U)^*\subseteq\pi^{-1}[U^*]$ with $\Phi(as)=\Phi(bs)$ and hence
\[\Phi'(\pi(a)\pi(s))=\Phi'(\pi(as))=\pi(\Phi(as))=\pi(\Phi(bs))=\Phi'(\pi(b)\pi(s)).\]
This shows that $\pi(a)\sim_U\pi(b)$ so $\overline\pi$ is well-defined.  We also immediately see that $\mathsf{U}_{\langle S\rangle\underline\pi}=\mathsf{U}_{\langle S'\rangle}\circ\overline\pi$ and $\overline\pi(U,[0,\underline\pi(U)])=[\pi(0),U]=[0,U]$, for all $U\in\mathrm{dom}(\underline\pi)$.

To see that $\overline\pi$ preserves products note that, for any $U,V\in\mathrm{dom}(\underline\pi)$ and $a,b\in S$,
\begin{align*}
\overline\pi((U,[a,\underline\pi(U)])(V,[b,\underline\pi(U)]))&=\overline\pi(U\cdot V,[ab,\underline\pi(U\cdot V)])\\
&=[\pi(ab),UV]\\
&=[\pi(a),U][\pi(b),V]\\
&=\overline\pi(U,[a,\underline\pi(U)])\overline\pi(V,[b,\underline\pi(U)])
\end{align*}
To see that $\overline\pi$ also preserves units, take any $(U,[z,\underline\pi(U)])\in(\underline\pi^{\mathsf{U}\langle S\rangle}\mathcal{U}[S])^0$.  Then $U\in\mathcal{U}(S')^0$ so $\underline\pi(U)\in\mathcal{U}(S)^0$ and $[z,\underline\pi(U)]\in\mathcal{U}[S]^0$.  By \eqref{FS0}, we may assume that $z\in\underline\pi(U)^Z$, i.e. $z\in Z$ and we have $s\in\underline\pi(U)$ with $sz=s$.  Then $\pi(z)\in\pi[Z]\subseteq Z'$, $\pi(s)\in U$ and $\pi(s)\pi(z)=\pi(s)$ so $\pi(z)\in U^{Z'}$ and hence $[\pi(z),U]\in\mathcal{U}[S']^0$, again by \eqref{FS0}.  This shows that $\overline\pi$ is a functor.

To see that $\overline\pi$ is continuous, take any $a\in S$ and $U\in\mathrm{dom}(\underline\pi)$.  Any neighbourhood of $[\pi(a),U]$ contains one of the form $\pi(a)_s$, for some $s\in U$, by \autoref{FSbasis}.  Now just note that $(U,[a,\underline\pi(U)])\in\mathcal{U}_s\times\mathcal{U}[a]$ and $\overline\pi[\mathcal{U}_s\times\mathcal{U}[a]]\subseteq\pi(a)_s$.  This completes the proof that $\overline\pi$ is a bundle morphism.
\end{proof}

We can now define the \emph{ultrafilter bundle functor} $\mathsf{U}$ as follows.

\begin{thm}\label{Ufunctor}
We have a functor $\mathsf{U}:\mathbf{SS}\rightarrow\mathbf{ACB}$ given by
\[\mathsf{U}(\langle S'\rangle,\pi,\langle S\rangle)=(\mathsf{U}_{\langle S'\rangle},\overline\pi,\underline\pi,\mathsf{U}_{\langle S\rangle}).\]
\end{thm}

\begin{proof}
We need to show that, whenever $\langle S\rangle=(S,Z,\Phi)$, $\langle S'\rangle=(S',Z',\Phi')$ and $\langle S''\rangle=(S'',Z'',\Phi'')$ are Steinberg semigroups and $\pi:S\rightarrow S'$ and $\pi':S'\rightarrow S''$ are Steinberg morphisms,
\[(\mathsf{U}_{\langle S''\rangle},\overline{\pi'\circ\pi},\underline{\pi'\circ\pi},\mathsf{U}_{\langle S\rangle})=(\mathsf{U}_{\langle S''\rangle},\overline{\pi'}\bullet\overline\pi,\underline\pi\circ\underline{\pi'},\mathsf{U}_{\langle S\rangle})\]
We already know $\underline{\pi'\circ\pi}=\underline\pi\circ\underline{\pi'}$, by \autoref{SSEG}, so we just need to show that
\[\overline{\pi'\circ\pi}=\overline{\pi'}\bullet\overline\pi.\]
To see this, take any $U\in\mathcal{U}(S'')$ and $a\in S$ and note that
\begin{align*}
\overline{\pi'}\bullet\overline\pi(U,[a,\underline{\pi'\circ\pi}(U)])
&=\overline{\pi'}(U,\overline\pi(\underline{\pi'}(U),[a,\underline\pi(\underline{\pi'}(U))]))\\
&=\overline{\pi'}(U,[\pi(a),\underline{\pi'}(U)])\\
&=[\pi'(\pi(a)),U]\\
&=\overline{\pi'\circ\pi}(U,[a,\underline{\pi'\circ\pi}(U)]).
\end{align*}
We also immediately see that $\mathsf{U}$ takes units to units so $\mathsf{U}$ is indeed a functor.
\end{proof}

Now we just have to show that $\mathsf{S}$ and $\mathsf{U}$ are inverse to each other modulo certain natural isomorphisms.  One of these comes from the ultrafilter representation, which we now denote by $\eta_{\langle S\rangle}$, for each Steinberg semigroup $\langle S\rangle$, i.e. for all $a\in S$,
\[\eta_{\langle S\rangle}(a)=\widehat{a}.\]

\begin{thm}\label{EtaIsomorphism}
We have a natural isomorphism $\eta$ from $\mathrm{id}_\mathbf{SS}$ to $\mathsf{S}\circ\mathsf{U}$ given by the ultrafilter representation $\eta_{\langle S\rangle}:S\rightarrow\mathcal{S}_\mathsf{c}(\mathsf{U}_{\langle S\rangle})$, for each Steinberg semigroup $\langle S\rangle$.
\end{thm}

\begin{proof}
By \autoref{FaithfulUltrafilters} and \autoref{AmpleUltrafilters}, the ultrafilter representation $\eta_{\langle S\rangle}$ yields an isomorphism $(\mathsf{SU}_{\langle S\rangle},\eta_{\langle S\rangle},\langle S\rangle)$ in $\mathbf{SS}$, for each Steinberg semigroup $\langle S\rangle$.  It only remains to prove naturality, i.e. we must show that, for any Steinberg morphism $\pi:S\rightarrow S'$ between Steinberg semigroups $\langle S\rangle=(S,Z,\Phi)$ and $\langle S'\rangle=(S',Z',\Phi')$,
\begin{align*}
(\mathsf{SU}_{\langle S'\rangle},\eta_{\langle S'\rangle},\langle S'\rangle)(\langle S'\rangle,\pi,\langle S\rangle)&=\mathsf{S}(\mathsf{U}(\langle S'\rangle,\pi,\langle S\rangle))(\mathsf{SU}_{\langle S\rangle},\eta_{\langle S\rangle},\langle S\rangle).
\end{align*}
This amounts to showing that
\begin{equation}\label{MorphEq}
\eta_{\langle S'\rangle}\circ\pi=\overline\pi_*\circ\underline\pi^*\circ\eta_{\langle S\rangle}.
\end{equation}
To see this note that, for any $a\in S$ and $U\in\mathcal{U}(S')$ with $\pi^{-1}[U]\neq\emptyset$,
\begin{align*}
(\overline\pi_*\circ\underline\pi^*\circ\eta_{\langle S\rangle}(a))(U)&=(\overline\pi_*\circ\underline\pi^*(\widehat{a}))(U)\\
&=\overline\pi(U,\widehat{a}(\underline\pi(U)))\\
&=\overline\pi(U,[a,\underline\pi(U)])\\
&=[\pi(a),U]\\
&=\widehat{\pi(a)}(U)\\
&=(\eta_{\langle S'\rangle}\circ\pi(a))(U).
\end{align*}
On the other hand, if $\pi^{-1}[U]=\emptyset$ then we claim that $[\pi(a),U]=0_U$.  Indeed, if we had $[\pi(a),U]\neq0_U$ then, taking $b\in S$ with $a<b$, it follows that $\pi(a)<\pi(b)$ and hence $\pi(b)\in\pi(a)^<\subseteq U$, by \eqref{aUup}, so $b\in\pi^{-1}[U]$, contradicting $\pi^{-1}[U]=\emptyset$.  This proves the claim and so again
\[(\eta_{\langle A'\rangle}\circ\pi(a))(U)=[\pi(a),U]=0_U=(\overline\pi_*\circ\underline\pi^*\circ\eta_{\langle S\rangle}(a))(U),\]
by the definition of $\overline\pi_*$ (see \eqref{BetaDef}), as $\pi^{-1}[U]=\emptyset$ means $U\notin\mathrm{dom}(\underline\pi)$.  Thus
\[\eta_{\langle S'\rangle}\circ\pi(a)=\overline\pi_*\circ\underline\pi^*\circ\eta_{\langle S\rangle}(a),\]
for all $a\in A$, proving \eqref{MorphEq}.
\end{proof}

For the other natural isomorphism first note, for each ample category bundle $\rho:C\twoheadrightarrow G$, we have an \'etale groupoid isomorphism $\varepsilon_\rho:G\rightarrow\mathcal{U}(\mathcal{S}_\mathsf{c}(\rho))$ given by
\begin{equation}\label{gamma->Sgamma}
\varepsilon_\rho(g)=\{a\in S:a(g)\in C^\times\}.
\end{equation}
Indeed, this follows from essentially the same argument as in \cite[Theorem 5.3]{BiceClark2021}.

\begin{prp}
For any ample category bundle $\rho:C\twoheadrightarrow G$ and $a,b\in\mathcal{S}_\mathsf{c}(\rho)$,
\begin{equation}\label{agbg}
a(g)=b(g)\qquad\Leftrightarrow\qquad a\sim_{\varepsilon_\rho(g)}b.
\end{equation}
\end{prp}

\begin{proof}
If $a(g)\neq b(g)$ then, for all $s\in\varepsilon_\rho(g)^*=\varepsilon_\rho(g^{-1})$,
\[\Phi^\rho(as)(gg^{-1})=a(g)s(g^{-1})\neq b(g)s(g^{-1})=\Phi^\rho(bs)(gg^{-1}).\]
This means $a\not\sim_{\varepsilon_\rho}b$.  Conversely, if $a(g)=b(g)$ then, as $\rho$ is locally injective, $a$ and $b$ agree on some open neighbourhood $O$ of $g$.  Then we can take $s\in\mathcal{S}_\mathsf{c}(\rho)$ such that $s(g^{-1})\in C^\times$ and $\mathrm{supp}(s)\subseteq O^{-1}$, as $\rho$ is an ample category bundle.  Then $s\in\varepsilon_\rho(g^{-1})=\varepsilon_\rho(g)^*$ and $\Phi^\rho(as)=\Phi^\rho(bs)$, showing that $a\sim_{\varepsilon_\rho(g)}b$.
\end{proof}

Now we pull ultrafilter bundles $\mathsf{U}_{\mathsf{S}\rho}:\mathcal{U}[\mathcal{S}_\mathsf{c}(\rho)]\twoheadrightarrow\mathcal{U}(\mathcal{S}_\mathsf{c}(\rho))$ back along $\varepsilon_\rho$.

\begin{prp}\label{EpsilonIsomorphism}
If $\rho:C\twoheadrightarrow G$ is an ample category bundle then we have a bundle isomorphism $\varepsilon^\rho:(\varepsilon_\rho)^{\mathsf{U}_{\mathsf{S}\rho}}\mathcal{U}[\mathcal{S}_\mathsf{c}(\rho)]\rightarrow C$ from $\mathsf{U}_{\mathsf{S}\rho\varepsilon_\rho}$ to $\rho$ given by
\[\varepsilon^\rho(g,[a,\varepsilon_\rho(g)])=a(g).\]
\end{prp}

\begin{proof}
By \eqref{agbg}, $a(g)=b(g)$ precisely when $[a,\varepsilon_\rho(g)]=[b,\varepsilon_\rho(g)]$, which shows that $\varepsilon^\rho$ is both well-defined and injective.  It is also surjective, as $\rho$ is an ample category bundle.  To see that $\varepsilon^\rho$ preserves the product, take $[a,\varepsilon_\rho(g)],[b,\varepsilon_\rho(h)]\in\mathcal{U}[\mathcal{S}_\mathsf{c}(\rho)]$ with $(g,h)\in G^2$.  Replacing $a$ and $b$ by other elements of $a_{\varepsilon_\rho(g)}$ and $b_{\varepsilon_\rho(h)}$ if necessary, we may assume $a\in\varepsilon_\rho(g)^>$ and $b\in\varepsilon_\rho(h)^>$.  In particular, $\mathrm{supp}(a)$ and $\mathrm{supp}(b)$ are contained in slices which themselves contain $g$ and $h$ respectively so
\[\varepsilon^\rho(gh,[ab,\varepsilon_\rho(gh)])=ab(gh)=a(g)b(h)=\varepsilon^\rho(g,[a,\varepsilon_\rho(g)])\varepsilon^\rho(h,[b,\varepsilon_\rho(h)]).\]
Also any unit of $\mathcal{U}[\mathcal{S}_\mathsf{c}(\rho)]$ is of the form $[z,\varepsilon^\rho(g)]$, for some $g\in G^0$ and $z\in\varepsilon^\rho(g)^{\mathcal{Z}_\mathsf{c}(\rho)}$, which means that $\varepsilon^\rho(g,[z,\varepsilon_\rho(g)])=z(g)=1_g$.  This shows that $\varepsilon^\rho$ also preserves units so $\varepsilon^\rho$ is a functor and hence an isomorphism of categories.

To see $\varepsilon^\rho$ is also a homeomorphism note, for every open $O\subseteq G$ and $a,s\in\mathcal{S}_\mathsf{c}(\rho)$,
\[\varepsilon^\rho[O\times a_s]=a[O\cap s^{-1}[C^\times]],\]
which is open in $C$, by \autoref{EtaleBundleOpenRange} and \autoref{OpenInvertibles}.  In fact, as $\rho$ is an ample category bundle, such open sets of form a basis for the topology of $C$.  As $\varepsilon^\rho$ maps one basis to another, it is a homeomorphism and hence bundle isomorphism.
\end{proof}

\begin{thm}\label{VarepsilonIsomorphism}
We have a natural isomorphism $\varepsilon$ from $\mathsf{U}\circ\mathsf{S}$ to $\mathrm{id}_\mathbf{ACB}$ given by the Pierce morphism $(\varepsilon^\rho,\varepsilon_\rho)$, for each ample category bundle $\rho:C\twoheadrightarrow G$.
\end{thm}

\begin{proof}
It only remains to show naturality, i.e. for all $(\rho',\beta,\phi,\rho)\in\mathbf{ACB}$,
\[(\rho',\beta,\phi,\rho)(\rho,\varepsilon^{\rho},\varepsilon_{\rho},\mathsf{U}_{\mathsf{S}\rho})=(\rho',\varepsilon^{\rho'},\varepsilon_{\rho'},\mathsf{U}_{\mathsf{S}\rho'})\mathsf{U}\circ\mathsf{S}(\rho',\beta,\phi,\rho).\]
Note $(\rho',\beta,\phi,\rho)(\rho,\varepsilon^{\rho},\varepsilon_{\rho},\mathsf{U}_{\mathsf{S}\rho})=(\rho',\beta\bullet\varepsilon^{\rho},\varepsilon_{\rho}\circ\phi,\mathsf{U}_{\mathsf{S}\rho})$ while
\begin{align*}
(\rho',\varepsilon^{\rho'},\varepsilon_{\rho'},\mathsf{U}_{\mathsf{S}\rho'})\mathsf{U}\circ\mathsf{S}(\rho',\beta,\phi,\rho)&=(\rho',\varepsilon^{\rho'},\varepsilon_{\rho'},\mathsf{U}_{\mathsf{S}\rho'})\mathsf{U}(\mathsf{S}\rho',\beta_*\circ\phi^*,\mathsf{S}\rho)\\
&=(\rho',\varepsilon^{\rho'},\varepsilon_{\rho'},\mathsf{U}_{\mathsf{S}\rho'})(\mathsf{U}_{\mathsf{S}\rho'},\overline{\beta_*\circ\phi^*},\underline{\beta_*\circ\phi^*},\mathsf{U}_{\mathsf{S}\rho})\\
&=(\rho',\varepsilon^{\rho'}\bullet\overline{\beta_*\circ\phi^*},\underline{\beta_*\circ\phi^*}\circ\varepsilon_{\rho'},\mathsf{U}_{\mathsf{S}\rho}).
\end{align*}
So we need to show that $\beta\bullet\varepsilon^{\rho}=\varepsilon^{\rho'}\bullet\overline{\beta_*\circ\phi^*}$ and $\varepsilon_{\rho}\circ\phi=\underline{\beta_*\circ\phi^*}\circ\varepsilon_{\rho'}$.

For the latter equation note that, for all $g\in G'$,
\begin{align*}
\underline{\beta_*\circ\phi^*}\circ\varepsilon_{\rho'}(g)
&=(\beta_*\circ\phi^*)^{-1}[\varepsilon_{\rho'}(g)]^<\\
&=\{s\in\mathcal{S}_\mathsf{c}(\rho):\beta_*\circ\phi^*(s)\in\varepsilon_{\rho'}(g)\}^<\\
&=\{s\in\mathcal{S}_\mathsf{c}(\rho):\beta_*\circ\phi^*(s)(g)\in\mathrm{dom}(\rho')^\times\}^<\\
&=\{s\in\mathcal{S}_\mathsf{c}(\rho):\beta(g,s(\phi(g)))\in\mathrm{dom}(\rho')^\times\}^<\\
&\supseteq\{s\in\mathcal{S}_\mathsf{c}(\rho):s(\phi(g))\in\mathrm{dom}(\rho)^\times\}^<\\
&=\varepsilon_\rho(\phi(g)).
\end{align*}
As inclusion on ultrafilters implies equality, it follows that
\[\underline{\beta_*\circ\phi^*}\circ\varepsilon_{\rho'}(g)=\varepsilon_\rho(\phi(g))=\varepsilon_\rho\circ\phi(g)\]
and hence $\underline{\beta_*\circ\phi^*}\circ\varepsilon_{\rho'}=\varepsilon_\rho\circ\phi$, as required.

For the former equation note that, for all $g\in G'$ and $a\in\mathcal{S}_\mathsf{c}(\rho)$,
\[\beta\bullet\varepsilon^{\rho}(g,[a,\varepsilon_\rho(\phi(g))])=\beta(g,\varepsilon^\rho(\phi(g),[a,\varepsilon_\rho(\phi(g))]))=\beta(g,a(\phi(g))),\]
while
\begin{align*}
\varepsilon^{\rho'}\bullet\overline{\beta_*\circ\phi^*}(g,[a,\varepsilon_\rho(\phi(g))])&=\varepsilon^{\rho'}(g,\overline{\beta_*\circ\phi^*}(\varepsilon_{\rho'}(g),[a,\varepsilon_\rho(\phi(g))]))\\
&=\varepsilon^{\rho'}(g,[\beta_*\circ\phi^*(a),\varepsilon_{\rho'}(g)])\\
&=\beta_*\circ\phi^*(a)(g)\\
&=\beta(g,a(\phi(g))).
\end{align*}
This shows that $\beta\bullet\varepsilon^{\rho}=\varepsilon^{\rho'}\bullet\overline{\beta_*\circ\phi^*}$ as well, as required.
\end{proof}

Putting all this together yields the main result of this section, namely that the category of Steinberg morphisms between Steinberg semigroups is equivalent to the category of Pierce morphisms between ample category bundles.

\begin{thm}\label{SSACB}
The category $\mathbf{SS}$ of Steinberg morphisms on Steinberg semigroups is equivalent to the category $\mathbf{ACB}$ of Pierce morphisms on ample category bundles.
\end{thm}

\begin{proof}
This is immediate from the fact that the section functor $\mathsf{S}:\mathbf{ACB}\rightarrow\mathbf{SS}$ and ultrafilter bundle functor $\mathsf{U}:\mathbf{SS}\rightarrow\mathbf{ACB}$ are inverse to each other, modulo the natural isomorphisms $\eta$ and $\varepsilon$, thanks to \autoref{EtaIsomorphism} and \autoref{VarepsilonIsomorphism}.
\end{proof}

\section{Steinberg Rings}\label{SteinbergRings}

Now we wish to extend our results with some extra additive structure.  Before defining Steinberg rings, we must first make a few preliminary definitions.

To start with, let us be clear that a \emph{ring} is just an abelian group $A$ with a compatible associative product, i.e. products distribute over sums.  We call a ring $A$ \emph{unital} if it also has a multiplicative unit $1$, i.e. $1a=a=a1$, for all $a\in A$ (so $(A,\cdot)$ is then a monoid, not just a semigroup).  A ring $A$ is \emph{locally unital} if every $a\in A$ has some central idempotent $z\in\mathsf{Z}(\mathsf{E}(A))$ fixing $a$, i.e. such that $za=a=az$.

If $(S,Z,\Phi)$ is a well-structured semigroup, we call any $T\subseteq S$ \emph{orthodirected} when every orthogonal pair in $T$ is dominated by a single element of $T$, i.e. for all $s,t\in T$,
\[\tag{Orthodirected}s\perp t\qquad\Rightarrow\qquad\exists r\in T\ (s,t<r).\]
If $S$ is a (multiplicative) subsemigroup of a ring $A$ then we say $S$ \emph{generates} $A$ when every element of $A$ is a finite sum of elements from $S$.  We also call any $T\subseteq S$ \emph{orthoadditive} or \emph{subtractive} respectively if, for all $s,t\in T$,
\begin{align}
\tag{Orthoadditive}s\perp t\qquad&\Rightarrow\qquad s+t\in T.\\
\tag{Subtractive}s\geq t\qquad&\Rightarrow\qquad s-t\in T.
\end{align}
We call a subtractive orthoadditive subsemigroup $T$ containing $0$ an \emph{orthosubring}.  For example, if $R$ is any commutative subring of $A$ contained in $S$ then its idempotents $\mathsf{E}(R)$ form an orthosubring, as products $st$, orthogonal sums $s+t$ and differences $s-t$, when $s\geq t$, are idempotent whenever $s$ and $t$ are idempotent.

We call a function $\Phi$ on a ring $A$ \emph{additive} if, for all $a,b\in A$,
\[\tag{Additive}\Phi(a+b)=\Phi(a)+\Phi(b).\]
Note that if $\Phi$ is additive on $A$ and shiftable on a subsemigroup $S$ generating $A$, i.e. $\Phi(sa)s=s\Phi(as)$, for all $a,s\in S$, then this automatically extends to $a\in A$ -- just take $s_1,\ldots,s_n\in S$ with $a=\sum_{k=1}^ns_k$ and note that
\[\Phi(sa)s=\Phi\Big(s\sum_{k=1}^ns_k\Big)s=\Phi\Big(\sum_{k=1}^nss_k\Big)s=\sum_{k=1}^n\Phi(ss_k)s=\sum_{k=1}^ns\Phi(s_ks)=s\Phi(as).\]

\begin{dfn}
A \emph{Steinberg ring} is a quadruple $(A,S,Z,\Phi)$ such that
\begin{enumerate}
\item $A$ is ring with orthodirected subsemigroup $S$ generating $A$.
\item $\Phi:A\rightarrow S$ is an additive shiftable expectation.
\item $Z$ is a binormal bistable orthosubring contained in $\mathsf{E}(\mathsf{Z}(\mathrm{ran}(\Phi)))$.
\end{enumerate}
\end{dfn}

To apply all the theory we have developed so far, we first note the following.

\begin{prp}\label{SteinbergRingsSemi}
Any Steinberg ring $(A,S,Z,\Phi)$ is, with respect to the product, a well-structured semimodule such that $(S,Z,\Phi)$ is a Steinberg semigroup.
\end{prp}

\begin{proof}
It is immediate from the definition that $(A,S,Z,\Phi)$ is a well-structured semimodule.  Also note that $S=S^>$, as $S$ is orthodirected and $0\perp s$, for all $s\in S$.  As $Z\subseteq\mathsf{E}(S)$, it follows that $S=S^{\dagger>}$, by \eqref{ReflexiveInterpolation}.

Next we claim that $S$ is also an orthosubring.  To see this, take any orthogonal $s,t\in S$.  As $S$ is orthodirected, we have $r\in S$ with $s,t<r$.  As $S=S^{\dagger>}$, we may assume that $r$ is $Z$-invertible so $s,t<_{r^{-1}}r$ and hence
\[s+t=\Phi(sr^{-1})r+\Phi(tr^{-1})r=\Phi(sr^{-1}+tr^{-1})r\in SS\subseteq S.\]
Likewise, if $s\geq t$, we have $r\in S^\dagger$ with $s<r$.  By \eqref{LeftAuxiliarity}, $t<r$ so
\[s-t=\Phi(sr^{-1})r-\Phi(tr^{-1})r=\Phi(sr^{-1}-tr^{-1})r\in SS\subseteq S.\]
As $S$ is also a subsemigroup containing $Z$ and hence $0$, $S$ is an orthosubring.

Now we show that $s+t=s\vee t$, for all orthogonal $s,t\in S$.  Indeed, if $s\perp t$ then we have $y\in Z$ with $s=sy$ and $ty=0$ and hence $s=sy=(s+t)y$, i.e. $s\leq s+t$.  Taking $z\in Z$ with $tz=t$, we see that $t(z-yz)=t$ and $y(z-yz)=yz-yz=0$.  As $Z$ is an orthosubring, we may thus replace $z$ with $z-yz\in Z$ if necessary to ensure that $yz=0$ so $sz=syz=0$ and hence $t=tz=(s+t)z$, i.e. $t\leq s+t$.  For any $r\geq s,t$, we have $y',z'\in Z$ with $s=sy'=ry'$ and $t=tz'=rz'$.  Replacing $y'$ and $z'$ with $y'y$ and $z'z$ if necessary, we may assume that $y'\leq y\perp z\geq z'$ so $y+z\in Z$ and $r(y'+z')=ry'+rz'=s+t=(s+t)(y'+z')$, i.e. $s+t\leq r$.  This shows that $s+t$ is indeed a supremum of $s$ and $t$.

We already showed that $S$ is orthoadditive so it now follows that $S$ satisfies \eqref{Orthosuprema}.  As $Z$ is also orthoadditive and products distribute over sums, \eqref{Distributivity} also holds.  It also follows that $y-z=y\setminus z$ whenever $y\geq z$.  Indeed, if $y\geq z$ then $z=yz$ and hence $(y-z)z=z-z=0$, i.e. $y-z\perp z$ so $(y-z)\vee z=(y-z)+z=y$, by what we just proved.  As $Z$ is subtractive, it follows that \eqref{Complements} also holds and hence $(S,Z,\Phi)$ is a Steinberg semigroup.
\end{proof}

Another key thing to note is that the expectation is always nondegenerate.

\begin{prp}\label{Nondeg}
If $(A,S,Z,\Phi)$ is a Steinberg ring then, for all $a\in A\setminus\{0\}$,
\[\tag{Nondegenerate}\label{Nondegenerate}\Phi[aS]\neq\{0\}.\]
\end{prp}

\begin{proof}
This is proved by induction on the number of elements of $S$ needed to generate $a$.  First note that, for any $a\in S=S^>$, we have $b,s\in S$ with $a<_sb$ and hence $\Phi(as)b=asb=a\neq0$ so, in particular, $\Phi(as)\neq0$.  Now say we have proved \eqref{Nondegenerate} for all $a\in S_n=\{\sum_{k=1}^ns_n:s_1,\ldots,s_n\in S\}$ and take $a\in S_{n+1}$, so we have $s_1,\ldots,s_{n+1}\in S$ with $a=\sum_{k=1}^ns_k$.  Then again we have $s,b\in S$ with $s_{n+1}<_sb$.  If $\Phi(as)\neq0$ then we are done, otherwise $\Phi(as)=0$ and hence
\[a=a-\Phi(as)b=\sum_{k=1}^{n+1}(s_k-\Phi(s_ks)b)=\sum_{k=1}^n(s_k-\Phi(s_ks)b),\]
as $s_{n+1}-\Phi(s_{n+1}s)b=s_{n+1}-s_{n+1}=0$.  But for each $k\leq n$, we have $y_k,z_k\in Z$ with $s_ky_k=s_k$ and $\Phi(s_ks)=s_ksz_k$ (see \autoref{BistableQC}) and hence
\[s_k-\Phi(s_ks)b=s_k(y_k-y_ksz_kb)\in s_k(Z-ZsZb)\subseteq s_k(Z-ZZ)\subseteq s_k\mathrm{ran}(\Phi)\subseteq S.\]
Thus $a\in S_n$ so we are again done by the inductive hypothesis.
\end{proof}

The first elementary examples come from locally unital rings.

\begin{prp}
Any locally unital ring $A$ yields a Steinberg ring given by
\[\langle A\rangle=(A,A,\mathsf{Z}(\mathsf{E}(A)),\mathrm{id}_A).\]
\end{prp}

\begin{proof}
As $A$ is locally unital, $A$ is directed and hence orthodirected.  Indeed, any $a,b\in A$ have local units $y,z\in\mathsf{Z}(\mathsf{E}(A))$, i.e. satisfying $a=ay$ and $b=bz$, and then $x=y+z-yz\in\mathsf{Z}(\mathsf{E}(A))$ satisfies $a,b<_xx$.  As noted above, $\mathsf{Z}(\mathsf{E}(A))$ also forms an orthosubring.  All the other required properties are immediate.
\end{proof}

\begin{rmk}\label{DaunsHofmannPierce}
The original work of Dauns-Hofmann \cite{DaunsHofmann1966} and Pierce \cite{Pierce1967} dealt only with \emph{biregular} rings, meaning that every principal ideal is generated by an idempotent.  Being locally unital is weaker than this so, even in the commutative case, the duality we obtain will be significantly more general.  On the topological side, this comes from the fact that, in contrast to \cite{DaunsHofmann1966} and Pierce \cite{Pierce1967}, we make no Hausdorff assumption on the total spaces of our bundles.  We also place no simplicity conditions on the fibres, which corresponds to the fact that we could replace $\mathsf{Z}(\mathsf{E}(A))$ above with some other orthosubring.
\end{rmk}

More general examples of Steinberg rings come from ample ringoid bundles.  Indeed, like before, we will soon see in \autoref{AmpleRingoidUltrafilters} below that every Steinberg ring is isomorphic to one arising from an ample ringoid bundle.

First recall that an \emph{abelian group bundle} is a bundle $\pi:C\rightarrow G$ where each fibre $\pi^{-1}\{g\}$ is an abelian group, the inverse $a\mapsto-a$ is continuous on $C$ and the sum $(a,b)\mapsto a+b$ is continuous on $C\times_GC$ (so $C$ is a topological abelian groupoid).

\begin{dfn}
We call $\rho:C\twoheadrightarrow G$ an \emph{ample ringoid bundle} if it is both an ample category bundle and an abelian group bundle such products distribute over sums, i.e. for all $a,b,c,d,e\in C$ with $\rho(a)=\rho(b)$, $\mathsf{s}(b)=\mathsf{r}(c)$, $\mathsf{s}(c)=\mathsf{r}(d)$ and $\rho(d)=\rho(e)$,
\[(a+b)c=ab+ac\qquad\text{and}\qquad c(d+e)=cd+ce.\]
\end{dfn}

If $G$ is a principal groupoid then $C$ will indeed be a ringoid/preadditive category in the usual sense.  This is not so in general, however, as sums are not defined for elements in different fibres, even if they have the same source and range units.

\begin{rmk}
A similar kind of structure is defined in \cite[\S4]{GoncalvesSteinberg2021} called a \emph{$G$-sheaf of unital rings}.  This consists of a local homeomorphism $p:E\twoheadrightarrow G^0$ and a continuous map $\alpha:G\mathop{{}_\mathsf{s}\times_p}E\rightarrow C$ satisfying analogous conditions.  In our terminology, a $G$-sheaf of unital rings is essentially the same as an abelian ringoid bundle $\rho:C\twoheadrightarrow G$ with a distinguished continuous section $a:G\rightarrow C^\times$ such that $a[G^0]=C^0$.  Indeed, given such a bundle and section, we can take $E=C^0$ and $p=\rho|_E$ and define $\alpha:G\mathop{{}_\mathsf{s}\times_p}E\rightarrow E$ by $\alpha(g,e)=a(g)ea(g)^{-1}$, for all $g\in G$ and $e\in E$ with $\mathsf{s}(g)=p(e)$.  Conversely, given a $G$-sheaf of unital rings consisting of $p:E\twoheadrightarrow G^0$ and $\alpha:G\mathop{{}_\mathsf{s}\times_p}E\rightarrow C$, we can define an ample ringoid bundle $\rho:G\mathop{{}_\mathsf{s}\times_p}E\twoheadrightarrow G$ where $\rho(g,e)=g$ and $(g,e)(h,f)=(gh,e\alpha(g,f))$, for all $(g,e),(h,f)\in G\mathop{{}_\mathsf{s}\times_p}E$ with $\mathsf{s}(g)=\mathsf{r}(h)$.  In this case, the corresponding distinguished continuous section is just given by $a(g)=(g,1_{\mathsf{s}(g)})$, where $1_{\mathsf{s}(g)}$ is the unit of the `stalk' at $\mathsf{s}(g)$.
\end{rmk}

\begin{prp}\label{SectionRings}
If $\rho:C\twoheadrightarrow G$ is an ample ringoid bundle, $(\mathcal{C}_\mathsf{c}(\rho),\mathcal{S}_\mathsf{c}(\rho),\mathcal{Z}_\mathsf{c}(\rho),\Phi^\rho_\mathsf{c})$ is a Steinberg ring where, for all $a,b\in\mathcal{C}_\mathsf{c}(\rho)$ and $g\in G$,
\begin{align*}
(a+b)(g)&=a(g)+b(g).\\
ab(g)&=\sum_{g=hi}a(h)b(i).
\end{align*}
Moreover, $C$ is Hausdorff precisely when each element of $\mathcal{S}_\mathsf{c}(\rho)$ is bisupported.
\end{prp}

\begin{proof}
For any $a\in\mathcal{C}_\mathsf{c}(\rho)$, we can cover $\mathrm{supp}(a)$ with finitely many compact open slices.  Taking differences, we can further ensure that these slices are disjoint.  The restriction of $a$ to these slices then yields elements of $\mathcal{S}_\mathsf{c}(\rho)$ whose sum is $a$, showing that $\mathcal{S}_\mathsf{c}(\rho)$ generates $\mathcal{C}_\mathsf{c}(\rho)$.  It also follows that $\mathrm{supp}(a)\cap\mathsf{s}^{-1}\{g\}$ and $\mathrm{supp}(a)\cap\mathsf{r}^{-1}\{g\}$ are finite, for any $g\in G$.  Thus the sum defining $ab$ or $ba$ has only finitely many non-zero terms, for any $b\in\mathcal{C}_\mathsf{c}(\rho)$, i.e. the product is well-defined.  Associativity and distributivity of the bundle then pass to the sections under these operations so $\mathcal{C}_\mathsf{c}(\rho)$ is indeed a ring.  As addition in $\mathcal{C}_\mathsf{c}(\rho)$ is defined fibrewise, $\Phi^\rho_\mathsf{c}$ is also additive.  The remaining properties required of a Steinberg ring follow from the fact that $(\mathcal{S}_\mathsf{c}(\rho),\mathcal{Z}_\mathsf{c}(\rho),\Phi^\rho_\mathsf{c})$ is a Steinberg semigroup, as shown in \autoref{AmpleBundleSteinbergSemigroup}.

If $C$ is Hausdorff then every $a\in\mathcal{S}_\mathsf{c}(\rho)$ is bisupported, by \autoref{OpenClosedSupports} and \autoref{AmpleBundleSteinbergSemigroup}.  Conversely, if $C$ is not Hausdorff then we have some net $(c_\lambda)\subseteq C$ with distinct limits $c$ and $c'$.  It then follows that $(0_{\rho(c_\lambda)})$ has non-zero limit $b=c-c'$.  Take $s\in\mathcal{S}_\mathsf{c}(\rho)$ such that $b\in\mathrm{ran}(s)$.  Restricting to a subnet if necessary, it follows that $(0_{\rho(c_\lambda)})$ also lies in $\mathrm{ran}(b)$, by \autoref{EtaleBundleOpenRange}.  Thus $(\rho(c_\lambda))$ lies outside $\mathrm{supp}(s)$, even though its limit $\rho(b)$ lies in $\mathrm{supp}(s)$.  Thus $\mathrm{supp}(s)$ is not open and hence $s$ is not bisupported, again by \autoref{AmpleBundleSteinbergSemigroup}.
\end{proof}

\subsection{Quasi-Cartan Pairs}\label{QuasiCartanPairs}

Here we will show that the quasi-Cartan pairs in \cite{ArmstrongCastroClarkCourtneyLinMcCormickRamaggeSimsSteinberg2021} also form Steinberg rings.  As the groupoid twists in \cite{ArmstrongCastroClarkCourtneyLinMcCormickRamaggeSimsSteinberg2021} are special kinds of groupoid bundles (see \autoref{TwistedRemark}), the present paper thus provides an alternative path to the theory developed in \cite{ArmstrongCastroClarkCourtneyLinMcCormickRamaggeSimsSteinberg2021}.

First let $A$ be an algebra over a commutative ring $R$.  We call $A$ \emph{torsion-free} if
\[\tag{Torsion-Free}r\in R,\ a\in\mathsf{E}(A)\text{ and }ra=0\quad\Rightarrow\quad r=0\text{ or }a=0.\]
The subset \emph{spanned} by any $B\subseteq A$ will be denoted by
\[\mathrm{span}(B)=\{\sum_{k=1}^nr_kb_k:r_1,\ldots,r_n\in R\text{ and }b_1,\ldots,b_k\in B\}.\]

\begin{dfn}
We call $(A,Z)$ a \emph{quasi-Cartan pair} if
\begin{enumerate}
\item $A$ is an algebra over a commutative ring $R$, which spanned by $Z^{\mathsf{N}\dagger}$.
\item $Z$ is commutative, torsion-free, spanned by $\mathsf{E}(Z)$ and the range of a linear expectation $\Phi$ on $A$ which is quasi-Cartan on $Z^{\mathsf{N}\dagger}$.
\end{enumerate}
\end{dfn}

These quasi-Cartan pairs come from \cite[Definition 3.3]{ArmstrongCastroClarkCourtneyLinMcCormickRamaggeSimsSteinberg2021}, where the expectation $\Phi$ is also required to satisfy \eqref{Nondegenerate}, although this is actually automatic, by \autoref{Nondeg} and \autoref{quasiCartanSteinberg} below.  Moreover, the expectation is unique, by \eqref{Leech}, as also noted in \cite[Proposition 3.7]{ArmstrongCastroClarkCourtneyLinMcCormickRamaggeSimsSteinberg2021}.

If $(A,Z)$ is a quasi-Cartan pair then, in particular, $(A,\mathsf{E}(Z),Z)$ is a semigroup inclusion and hence we may define the restriction, domination, and orthogonality relations on $A$ as in \autoref{Relations}.  We then define the \emph{orthospan} of any $B\subseteq A$ to be the span of all finite orthogonal subsets $F\subseteq B$ (where $a\perp b$, for all distinct $a,b\in F$), i.e.
\begin{align*}
\mathrm{orthospan}(B)&=\bigcup\{\mathrm{span}(F):F\text{ is a finite orthogonal subset of }B\}\\
&=\big\{\sum_{k=1}^nr_kb_k:r_1,\ldots,r_n\in R\text{ and }b_1,\ldots,b_n\in B\text{ are orthogonal}\big\}.
\end{align*}
We are particularly interested in the orthospan of the $Z$-invertible normalisers $\mathrm{orthospan}(Z^{\mathsf{N}\dagger})$.  This gives us precisely the right subsemigroup of $A$ to turn any quasi-Cartan pair into a Steinberg ring.

\begin{thm}\label{quasiCartanSteinberg}
If $(A,Z)$ is a quasi-Cartan pair and $\Phi$ is the unique quasi-Cartan expectation onto $Z$ then $(A,\mathrm{orthospan}(Z^{\mathsf{N}\dagger}),\mathsf{E}(Z),\Phi)$ is a Steinberg ring.
\end{thm}

\begin{proof}
As $Z$ is the range of a linear expectation, $Z$ is a subalgebra of $A$.  As $Z$ is commutative, $\mathsf{E}(Z)$ is an orthosubring and a generalised Boolean algebra, where
\[y\wedge z=yz,\quad y\vee z=y+z-yz\quad\text{and}\quad y\setminus z=y-yz.\]
If $y,z\in\mathsf{E}(Z)$ then $y-yz$, $z-yz$ and $yz$ are orthogonal so, for any $q,r\in R$,
\[qy+rz=q(y-yz)+r(z-yz)+(q+r)yz\in\mathrm{orthospan}(\mathsf{E}(Z)).\]
This extends to arbitrary finite subsets and hence
\[Z=\mathrm{span}(\mathsf{E}(Z))=\mathrm{orthospan}(\mathsf{E}(Z))\subseteq\mathrm{orthospan}(Z^{\mathsf{N}\dagger}).\]

To see that $\mathrm{orthospan}(Z^{\mathsf{N}\dagger})$ is a subsemigroup of $A$, take $a,b\in\mathrm{orthospan}(Z^{\mathsf{N}\dagger})$ so we have orthogonal $a_1,\ldots,a_m\in Z^{\mathsf{N}\dagger}$ and orthogonal $b_1,\ldots,b_n\in Z^{\mathsf{N}\dagger}$ such that $a\in\mathrm{span}(a_k)$ and $b\in\mathrm{span}(b_k)$.  By \eqref{abperpcd}, \eqref{OrthoInvariance} and \autoref{ZInvertibleNormalisers}, $(a_jb_k)_{j\leq m}^{k\leq n}$ are also orthogonal elements of $Z^{\mathsf{N}\dagger}$ and hence $ab\in\mathrm{span}(a_jb_k)\subseteq\mathrm{orthospan}(Z^{\mathsf{N}\dagger})$.

For all $a\in A^\dagger$, note that $a=aa^{-1}a\in a^{\mathsf{E}(Z)}\cap{}^{\mathsf{E}(Z)}\hspace{-1pt}a$, i.e.
\[A^\dagger\subseteq\{a\in A:a\in a^{\mathsf{E}(Z)}\cap{}^{\mathsf{E}(Z)}\hspace{-1pt}a\}.\]
Note this latter set is closed under taking scalar products and sums (e.g. note that if $a=ay$ and $b=bz$, for some $y,z\in\mathsf{E}(Z)$, then $a+b=(a+b)(y\vee z)$).  Thus
\[A=\mathrm{span}(A^\dagger)\subseteq\{a\in A:a\in a^{\mathsf{E}(Z)}\cap{}^{\mathsf{E}(Z)}\hspace{-1pt}a\}.\]
This implies that orthogonality is symmetric, as in the proof of \autoref{perpSymmetric}, i.e. $a\perp b$ implies $b\perp a$.  Also note that the orthogonality relation respects sums, i.e.
\begin{equation}\label{perpSums}
a\perp b,c\qquad\Rightarrow\qquad a\perp b+c.
\end{equation}
Indeed, if $a\perp b,c$ then we have $y,z\in\mathsf{E}(Z)$ with $a=ay=az$ and $by=cz=0$.  Then $yz\in\mathsf{E}(Z)$ satisfies $a=az=ayz$ and $(b+c)yz=byz+czy=0$.  This and a dual argument shows that $a\perp b+c$.

Likewise, we immediately see that $Z^\mathsf{N}$ is closed under taking scalar products.  We claim that $Z^\mathsf{N}$ is also closed under orthogonal sums.  To see this, take $a,b\in Z^\mathsf{N}$ with $a\perp b$.  In particular, we have $x\in a^Z$ and $y\in b^Z$ with $bx=0=ay$ and $z_a,z_b\in Z$ with $za=az_a$ and $zb=bz_b$.  Thus
\[z(a+b)=az_a+bz_b=axz_a+byz_b=(a+b)(xz_a+yz_b)\in(a+b)Z.\]
This and a dual argument yields $Z(a+b)=(a+b)Z$, i.e. $a+b\in Z^\mathsf{N}$.  By \eqref{perpSums}, this also extends to finite orthogonal sums and hence
\[\mathrm{orthospan}(Z^{\mathsf{N}\dagger})\subseteq\mathrm{orthospan}(Z^\mathsf{N})\subseteq Z^\mathsf{N}.\]

Next we claim that $\mathrm{orthospan}(Z^{\mathsf{N}\dagger})$ is orthodirected.  Accordingly, take non-zero $r_1,\ldots,r_m,s_1,\ldots,s_n\in R$ and orthogonal $a_1,\ldots,a_m,b_1,\ldots,b_n\in Z^{\mathsf{N}\dagger}$ such that $a=\sum_{k=1}^mr_ka_k$ and $b=\sum_{k=1}^ns_kb_k$ are orthogonal.  This means we have $z\in a^Z$ with $bz=0$.  For all $k\leq m$, note $r_ka_k^{-1}a_k=a_k^{-1}a=a_k^{-1}az=r_ka_k^{-1}a_kz$.  So $r_k(a_k^{-1}a_k-a_k^{-1}a_kz)=0$ and hence $a_k^{-1}a_k-a_k^{-1}a_kz=0$, as $Z$ is torsion-free.  This yields $a_k=a_ka_k^{-1}a_k=a_ka_k^{-1}a_kz=a_kz$.  Similarly, for all $k\leq n$, we see that $s_kb_k^{-1}b_kz=b_k^{-1}bz=0$ so $b_kz=b_kb_k^{-1}b_kz=0$, as $Z$ is torsion free.  This means $a_j\perp b_k$, for all $j\leq m$ and $k\leq n$.  Thus $r=\sum_{k=1}^ma_k+\sum_{k=1}^nb_k\in Z^{\mathsf{N}\dagger}$, by what we just proved.  Specifically, $r$ has $Z$-inverse $r^{-1}=\sum_{k=1}^ma_k^{-1}+\sum_{k=1}^nb_k^{-1}$ and elementary calculations show that $a,b<_{r^{-1}}r$.  This shows that $\mathrm{orthospan}(Z^{\mathsf{N}\dagger})$ is indeed orthodirected and, moreover,
\begin{equation}\label{OrthoDominated}
\mathrm{orthospan}(Z^{\mathsf{N}\dagger})\subseteq Z^{\mathsf{N}\dagger>}.
\end{equation}

Now we further claim that $\mathrm{orthospan}(Z^{\mathsf{N}\dagger})\subseteq\mathsf{E}(Z)^\mathsf{N}$.  To see this, take any $a\in\mathrm{orthospan}(Z^{\mathsf{N}\dagger})$, so \autoref{ZInvertibleNormalisers} and \eqref{OrthoDominated} yield $b\in Z^{\mathsf{N}\dagger}\subseteq\mathsf{E}(Z)^\mathsf{N}$ with $a<_{b^{-1}}b$.  Then $b^{-1}a\in Z\subseteq\mathsf{E}(Z)^\mathsf{C}\subseteq\mathsf{E}(Z)^\mathsf{N}$ so
\[a\mathsf{E}(Z)=bb^{-1}a\mathsf{E}(Z)=b\mathsf{E}(Z)b^{-1}a=\mathsf{E}(Z)bb^{-1}a=\mathsf{E}(Z)a.\]
This proves the claim, which means both $Z$ and $\mathsf{E}(Z)$ are normal in $\mathrm{orthospan}(Z^{\mathsf{N}\dagger})$, i.e. $(\mathrm{orthospan}(Z^{\mathsf{N}\dagger}),\mathsf{E}(Z),Z)$ is a structured semigroup.

Recall from the proof of \autoref{SteinbergRingsSemi} that orthogonal sums are suprema.  As $\Phi$ is quasi-Cartan on $Z^{\mathsf{N}\dagger}$, for any orthogonal $a,b\in Z^{\mathsf{N}\dagger}$,
\[\Phi(a+b)=\Phi(a)+\Phi(b)=\Phi(a)\vee\Phi(b)\leq a\vee b=a+b.\]
It follows that $\Phi$ is also quasi-Cartan on $\mathrm{orthospan}(Z^{\mathsf{N}\dagger})$.  By \eqref{BistableQCShiftable} and \eqref{OrthoDominated}, $Z$ is then bistable in $\mathrm{orthospan}(Z^{\mathsf{N}\dagger})$ and $\Phi$ is shiftable with respect to $\mathrm{orthospan}(Z^{\mathsf{N}\dagger})$ on $Z^{\mathsf{N}\dagger}$, i.e. $\Phi(sa)s=s\Phi(as)$, for all $a\in Z^{\mathsf{N}\dagger}$ and $s\in\mathrm{orthospan}(Z^{\mathsf{N}\dagger})$, which immediately extends to all $a\in A=\mathrm{span}(Z^{\mathsf{N}\dagger})$.  This completes the proof that $(A,\mathrm{orthospan}(Z^{\mathsf{N}\dagger}),\mathsf{E}(Z),\Phi)$ is a Steinberg ring.
\end{proof}

Finally, let us just note that every element of $Z^{\mathsf{N}\dagger}$ and hence $\mathrm{orthospan}(Z^{\mathsf{N}\dagger})$ is bisupported.  Thus the ultrafilter bundle arising from quasi-Cartan pairs will always have Hausdorff total space $\mathcal{U}[A]$, thanks to \autoref{FaithfulUltrafilters}, \autoref{AmpleUltrafilters}, \autoref{SectionRings} and \autoref{AmpleRingoidUltrafilters} below.

\subsection{Pierce-Steinberg Duality}\label{PierceSteinbergDuality}

Our goal here is to extend the duality in \autoref{LawsonSteinbergDuality} between Steinberg semigroups and ample category bundles to a duality between Steinberg rings and ample ringoid bundles.  Most the of the hard work has already been done, it is just a matter of verifying that all the constructions in \autoref{LawsonSteinbergDuality} also respect the extra additive structure.

First recall that the ultrafilter bundle $\mathsf{U}_{\langle A\rangle}:\mathcal{U}[A]\twoheadrightarrow\mathcal{U}(S)$ from \autoref{FilterBundles} is an ample category bundle whenever $(S,Z,\Phi)$ is a Steinberg semigroup, by \autoref{AmpleUltrafilters}.  In particular, this applies to Steinberg rings, thanks to \autoref{SteinbergRingsSemi}.

\begin{prp}\label{AmpleRingoidUltrafilters}
If $\langle A\rangle=(A,S,Z,\Phi)$ is a Steinberg ring, $\mathsf{U}_{\langle A\rangle}:\mathcal{U}[A]\twoheadrightarrow\mathcal{U}(S)$ is an ample ringoid bundle where the additive structure of the fibres is given by
\[[a,U]+[b,U]=[a+b,U].\]
\end{prp}

\begin{proof}
First we must show that the additive structure of $\mathcal{U}[A]$ above is well-defined.  To see this, say $a\sim_Ua'$ and $b\sim_Ub'$ so we have $u,v\in U^*$ with $\Phi(au)=\Phi(a'u)$ and $\Phi(bv)=\Phi(b'v)$.  Taking $w\in U$ with $w<u,v$, it follows that $\Phi(aw)=\Phi(a'w)$ and $\Phi(bw)=\Phi(b'w)$.  As $\Phi$ is additive, $\Phi((a+b)w)=\Phi((a'+b')w)$ and hence $a+b\sim_Ua'+b'$, showing that $+$ above is indeed well defined.

The abelian group structure of $A$ then immediately passes to the fibres, as does the distributivity of products over sums.  To see that the addition is also continuous just note that, for any basic neighbourhood $(a+b)_u$ of $[a+b,U]$ (see \autoref{FSbasis}), $a_u$ and $b_u$ are neighbourhoods of $[a,U]$ and $[b,U]$ such that $a_u+b_u=(a+b)_u$.  Likewise, if $a_u$ is a basic neighbourhood of $[a,U]$ then $(-a)_u=-(a_u)$ is a basic neighbourhood of $[-a,U]$ so the additive inverse map is also continuous.  It follows that $\mathsf{U}_{\langle A\rangle}$ is an abelian group bundle and hence an ample ringoid bundle, as we already know $\mathsf{U}_{\langle A\rangle}$ is an ample category bundle.
\end{proof}

Next recall the ultrafilter representation $a\mapsto\widehat{a}$ from \autoref{ContinuousSection}.

\begin{thm}\label{IsomorphismUltrafilters}
If $\langle A\rangle=(A,S,Z,\Phi)$ is a Steinberg ring then $a\mapsto\widehat{a}$ is an isomorphism from $A$ onto the ring of all compactly supported continuous sections $\mathcal{C}_\mathsf{c}(\mathsf{U}_{\langle A\rangle})$.
\end{thm}

\begin{proof}
For all $a,b\in A$ and $U\in\mathcal{U}(S)$, the definition of addition in $\mathcal{U}[A]$ yields
\[\widehat{a+b}(U)=[a+b,U]=[a,U]+[b,U]=\widehat{a}(U)+\widehat{b}(U).\]
Thus $\widehat{a+b}=\widehat{a}+\widehat{b}$.  We also already know that $\widehat{\,ab\,}=\widehat{a}\widehat{b}$, as long as $b\in S$, by \autoref{ContinuousSection}.  To extend this to any $b\in A$ just note that, as $S$ generates $A$, we have $b_1,\ldots,b_n\in S$ with $b=\sum_{k=1}^nb_k$ and hence
\[\widehat{\,ab\,}=\sum_{k=1}^n\widehat{\,ab_k}=\sum_{k=1}^n\widehat{a}\widehat{b_k}=\widehat{a}\sum_{k=1}^n\widehat{b_k}=\widehat{a}\widehat{b}.\]
Thus the ultrafilter representation $a\mapsto\widehat{a}$ is a ring homomorphism.

To see that $a\mapsto\widehat{a}$ is injective, take $a\in A\setminus\{0\}$.  By \eqref{Nondegenerate}, we have $s\in S$ with $0\neq\Phi(as)\in\mathrm{ran}(\Phi)\subseteq S$ so \autoref{FaithfulUltrafilters} yields $\widehat{0}\neq\widehat{\Phi(as)}=\Phi^{\mathsf{U}_{\langle A\rangle}}(\widehat{a}\widehat{s})$, by \eqref{ContinuousPhi}, and hence $\widehat{a}\neq\widehat{0}$.

Lastly note that $\mathcal{S}_\mathsf{c}(\rho)$ generates $\mathcal{C}_\mathsf{c}(\rho)$, for any ample ringoid bundle $\rho$.  By \autoref{AmpleUltrafilters}, $\widehat{S}=\mathcal{S}_\mathsf{c}(\mathsf{U}_{\langle A\rangle})$ and hence $\widehat{A}=\mathcal{C}_\mathsf{c}(\mathsf{U}_{\langle A\rangle})$, showing that $a\mapsto\widehat{a}$ is indeed an isomorphism onto $\mathcal{C}_\mathsf{c}(\mathsf{U}_{\langle A\rangle})$.
\end{proof}

Combined with \eqref{ContinuousPhi}, \autoref{AmpleUltrafilters} and \autoref{SectionRings}, it follows that Steinberg rings are precisely the quadruples of the form $(\mathcal{C}_\mathsf{c}(\rho),\mathcal{S}_\mathsf{c}(\rho),\mathcal{Z}_\mathsf{c}(\rho),\Phi^\rho_\mathsf{c})$, for some ample ringoid bundle $\rho$, at least up to isomorphism.

To turn Steinberg rings into a category, we consider the following morphisms.

\begin{dfn}
If $\langle A\rangle=(A,S,Z,\Phi)$ and $\langle A'\rangle=(A',S',Z',\Phi')$ are Steinberg rings, an \emph{additive Stienberg morphism} from $\langle A\rangle$ to $\langle A'\rangle$ is a ring homomorphism $\pi:A\rightarrow A'$ such that $\pi[S]\subseteq S'$, $\pi[Z]\subseteq Z'$ and $\pi(\Phi(a))=\Phi'(\pi(a))$, for all $a\in A$.
\end{dfn}

The category of such morphisms between Steinberg rings is thus given by
\[\mathbf{SR}=\{(\langle A'\rangle,\pi,\langle A\rangle):\pi\text{ is an additive Steinberg morphism from }\langle A\rangle\text{ to }\langle A'\rangle\},\]
the product being composition, again when the domain and codomain match, i.e.
\[(\langle A''\rangle,\pi',\langle A'\rangle)(\langle A'\rangle,\pi,\langle A\rangle)=(\langle A''\rangle,\pi'\circ\pi,\langle A\rangle).\]

Note that if $\pi$ is an additive Steinberg morphism then its restriction to $S$ is indeed a Steinberg morphism, as per \autoref{SteinbergMorphisms}.  Indeed, as $\pi$ is a ring homomorphism, $\pi(0)=0$ and $\pi(a\vee b)=\pi(a+b)=\pi(a)+\pi(b)=\pi(a)\vee\pi(b)$ when $a\perp b$.  In particular, it follows that $(\overline\pi,\underline\pi)=(\overline{\pi|_S},\underline{\pi|_S})$ is a Pierce morphism between the ultrafilter bundles $\mathsf{U}_{\langle A\rangle}$ and $\mathsf{U}_{\langle A'\rangle}$, by \autoref{UltraFunctoriality} and \autoref{PullbackBundleMorphism}.

Next note that a pullback of an ample ringoid bundle $\rho:C\twoheadrightarrow G'$ along an \'etale morphism $\phi:G\rightarrow G'$ between Hausdorff ample groupoids $G$ and $G'$ is again an ample ringoid bundle, where the extra additive structure of $\phi^\rho C$ comes from that of $C$, i.e. for all $g\in G$ and $c,d\in\rho^{-1}\{\phi(g)\}$,
\[(g,c)+(g,d)=(g,c+d).\]
We call a Pierce morphism $(\beta,\phi)$ between ample ringoid bundles \emph{additive} if
\[\beta(c+d)=\beta(c)+\beta(d),\]
for all $c$ and $d$ in the same fibre of the pullback bundle defined by $\phi$.  The category of additive Pierce morphisms between ample ringoid bundles is thus given by
\[\mathbf{ARB}=\{(\rho',\beta,\phi,\rho):(\beta,\phi)\text{ is an additive Pierce morphism from $\rho$ to $\rho'$}\},\]
with the product given as in \autoref{ACBcategory} by
\[(\rho'',\beta',\phi',\rho')(\rho',\beta,\phi,\rho)=(\rho'',\beta'\bullet\beta,\phi\circ\phi',\rho),\]
where $\beta'\bullet\beta(g'',b)=\beta'\midscript{\phi'}\beta(g'',b)=\beta'(g'',\beta(\phi'(g''),b))$.

\begin{prp}
If $\pi:A\rightarrow A'$ is an additive Steinberg morphism then $(\overline\pi,\underline\pi)$ is an additive Pierce morphism.
\end{prp}

\begin{proof}
Just note that, for any $a,b\in A$ and $U\in\mathcal{U}(S')$,
\begin{align*}
\overline\pi((U,[a,\underline\pi(U)])+(U,[b,\underline\pi(U)]))&=\overline\pi(U,[a+b,\underline\pi(U)])\\
&=[\pi(a+b),U]\\
&=[\pi(a),U]+[\pi(b),U]\\
&=\overline\pi(U,[a,\underline\pi(U)])+\overline\pi(U,[b,\underline\pi(U)]).\qedhere
\end{align*}
\end{proof}

As in \autoref{Ufunctor}, we thus have a functor $\mathsf{U}:\mathbf{SR}\rightarrow\mathbf{ARB}$ given by
\[\mathsf{U}(\langle A'\rangle,\pi,\langle A\rangle)=(\mathsf{U}_{\langle A'\rangle},\overline\pi,\underline\pi,\mathsf{U}_{\langle A\rangle}).\]

Given an ample ringoid bundle $\rho:C\twoheadrightarrow G$, let us denote the corresponding Steinberg ring from \autoref{SectionRings} by 
\[\mathsf{C}\rho=(\mathcal{C}_\mathsf{c}(\rho),\mathcal{S}_\mathsf{c}(\rho),\mathcal{Z}_\mathsf{c}(\rho),\Phi^\rho_\mathsf{c}).\]

\begin{prp}
If $(\beta,\phi)$ is an additive Pierce morphism from $\rho$ to $\rho'$ then $\beta_*\circ\phi^*$ is an additive Steinberg morphism from $\mathsf{C}\rho$ to $\mathsf{C}\rho'$.
\end{prp}

\begin{proof}
By \autoref{phiSteinberg} and \autoref{betaSteinberg}, we already know that $\beta_*\circ\phi^*$ is a Steinberg morphism from $\mathsf{S}\rho$ to $\mathsf{S}\rho'$.  As in the proof of \autoref{IsomorphismUltrafilters}, it thus suffices to show that $\beta_*\circ\phi^*$ is an additive group homomorphism from $\mathcal{C}_\mathsf{c}(\rho)$ to $\mathcal{C}_\mathsf{c}(\rho')$.  To see this just note that, for all $a,b\in\mathcal{C}_\mathsf{c}(\rho)$ and $g'\in\mathrm{ran}(\rho')$,
\begin{align*}
\beta_*\circ\phi^*(a+b)(g')&=\beta(g',(a+b)(\phi(g')))\\
&=\beta(g',a(\phi(g')))+\beta(g',b(\phi(g')))\\
&=\beta_*\circ\phi^*(a)(g')+\beta_*\circ\phi^*(b)(g'),
\end{align*}
and hence $\beta_*\circ\phi^*(a+b)=\beta_*\circ\phi^*(a)+\beta_*\circ\phi^*(b)$, as required.
\end{proof}

As in \autoref{Afunctor}, we thus have a functor $\mathsf{C}:\mathbf{ARB}\rightarrow\mathbf{SR}$ given by
\[\mathsf{C}(\rho',\beta,\phi,\rho)=(\mathsf{C}\rho',\beta_*\circ\phi^*,\mathsf{C}\rho).\]
As in \autoref{EtaIsomorphism}, we also then have a natural isomorphism $\eta$ from $\mathrm{id}_\mathbf{SR}$ to $\mathsf{C}\circ\mathsf{U}$ given by the ultrafilter representation $\eta_{\langle A\rangle}:A\rightarrow\mathcal{C}_\mathsf{c}(\mathsf{U}_{\langle A\rangle})$, i.e. for all $a\in A$,
\[\eta_{\langle A\rangle}(a)=\widehat{a}.\]
Elementary arguments like those above also show that the Pierce morphism $(\varepsilon^\rho,\varepsilon_\rho)$ from \eqref{gamma->Sgamma} and \autoref{EpsilonIsomorphism} is also additive and thus yields another natural isomorphism $\varepsilon$ from $\mathsf{U}\circ\mathsf{C}$ to $\mathrm{id}_\mathbf{ARB}$.

Putting all this together yields the final result we have been working towards.

\begin{cor}\label{SRARB}
The category $\mathbf{SR}$ of additive Steinberg morphisms between Steinberg rings is equivalent to the category $\mathbf{ARB}$ of additive Pierce morphisms between ample ringoid bundles.
\end{cor}

\begin{proof}
As in the proof of \autoref{SSACB}, this is immediate from the fact that the functors $\mathsf{C}:\mathbf{ARB}\rightarrow\mathbf{SR}$ and $\mathsf{U}:\mathbf{SR}\rightarrow\mathbf{ARB}$ are inverse to each other, modulo the natural isomorphisms $\eta$ and $\varepsilon$.
\end{proof}

\bibliography{maths}{}
\bibliographystyle{alphaurl}

\end{document}